\documentclass[]{article}
\usepackage[utf8]{inputenc}
\usepackage[a4paper,top=2.5cm,bottom=2.5cm,left=2cm,right=2cm]{geometry}
\usepackage{graphicx} 
\usepackage{float}
\usepackage{textcomp}
\usepackage{xcolor}
\usepackage{amsmath}
\usepackage{amssymb}
\usepackage{amsthm}
\usepackage{amsfonts}
\usepackage[english]{babel}
\usepackage{verbatim}
\usepackage{array}
\usepackage{hyperref}
\usepackage{cases}
\usepackage{esint}
\usepackage{multirow}
\usepackage{enumerate}
\usepackage{cite}
\usepackage{afterpage}
\usepackage{url}    
\usepackage{esint}    
\usepackage{caption}
\usepackage{subcaption} 
\usepackage{dsfont}                                                                                      
\newcommand{\authoraddress}[2]{%
	\textsc{#1} \textit{E-mail address:} \protect\url{#2}%
}                                                                                                                     
\newcommand{\AuthorAddressone}{                                                                                          
	\authoraddress{Institute for Applied Mathematics, University of Bonn, 53115 Bonn, Germany.}{dematte@iam.uni-bonn.de}%
} 
\newcommand{\AuthorAddresstwo}{                                                                                          
	\authoraddress{Institute for Applied Mathematics, University of Bonn, 53115 Bonn, Germany.}{velazquez@iam.uni-bonn.de}%
} 
\numberwithin{equation}{section}
\newtheorem{theorem}{Theorem}[section]
\newtheorem{lemma}{Lemma}[section]
\newtheorem{corollary}{Corollary}[section]
\newtheorem{prop}{Proposition}[section]

\theoremstyle{definition}
\newtheorem{definition}{Problem}[section]

\theoremstyle{remark}
\newtheorem*{remark}{Remark}

\newcommand{\Ss}{\mathbb{S}}
\newcommand{\R}{\mathbb{R}}

\newcommand{\eps}{\varepsilon}

\DeclareMathOperator*{\Div}{div}
\DeclareMathOperator*{\osc}{osc}
\DeclareMathOperator*{\artanh}{artanh}
\title{Traveling waves for a two-phase Stefan problem with radiation}
\author{Elena Demattè\thanks{\AuthorAddressone}, Juan J.L. Velázquez\thanks{\AuthorAddresstwo}}

\begin{document}
	
	\maketitle
	\begin{abstract}
		In this paper we study the existence of traveling wave solutions for a free-boundary problem modeling the phase transition of a material where the heat is transported by both conduction and radiation. Specifically, we consider a one-dimensional two-phase Stefan problem with an additional non-local non-linear integral term describing the situation in which the heat is transferred in the solid phase also by radiation, while the liquid phase is completely transparent, not interacting with radiation.	We will prove that there are traveling wave solutions for the considered model, differently from the case of the classical Stefan problem in which only self-similar solutions with the parabolic scale $ x\sim \sqrt{t} $ exist. In particular we will show that there exist traveling waves for which the solid expands. The properties of these solutions will be studied using maximum-principle methods, blow-up limits and Liouville-type Theorems for non-linear integral-differential equations. 
	\end{abstract}
\textbf{Acknowledgments:} The authors gratefully acknowledge the financial support of Bonn International Graduate School of Mathematics (BIGS) at the Hausdorff Center for Mathematics founded through the Deutsche Forschungsgemeinschaft (DFG, German Research Foundation) under Germany’s Excellence Strategy – EXC-2047/1 – 390685813. \\


\noindent\textbf{Statements and Declarations:} The authors have no relevant financial or non-financial interests to disclose.\\

\noindent\textbf{Data availability:} Data sharing not applicable to this article as no datasets were generated or analysed during the current study.
	\tableofcontents
	\section{Introduction}
	In this paper we continue the study of the free-boundary problem presented in \cite{StefanRad1} considering a one-dimensional Stefan-like problem which describes the melting of ice (resp. the solidification of water) under the assumption that the heat is transported by conduction in both phases of the material and additionally by radiation in the solid. To be more precise, we are studying the problem in which $ \R^3  $ is divided in two regions, one liquid region with a temperature $ T $ greater then the melting temperature $ T_M $ and one solid region with $ 0<T<T_M $. At the contact surface between the two phases the temperature satisfies $ T=T_M $. This surface moves as the solid melts or the liquid solidifies and it is thus called moving interface. Analogously to the classical Stefan problem, the heat is transferred by conduction in both the liquid and the solid phase. In addition we assume that the heat is transported also by radiation only in the solid. Equivalently, we assume the liquid to be perfectly transparent not allowing any interaction with radiation. 
	
	At the initial time $ t=0 $, the liquid is considered to fill the half-space $ \R_-^3=\{x\in\R^3: x_1<0\} $ and the solid to fill $ \R^3_+=\{x\in\R^3: x_1>0\} $. Thus, the interface is initially the plane $ \{0\}\times\R^2 $. Furthermore, we assume the temperature to depend only on the first variable, i.e. $ T(t,x)=T(t,x_1) $. This implies that the interface is described by the plane $ \{s(t)\}\times \R^2 $ for all $ t>\geq0 $ and the problem reduces to the study of a one-dimensional model.
	\begin{figure}[H]
		\begin{center}
			\includegraphics[height=5cm]{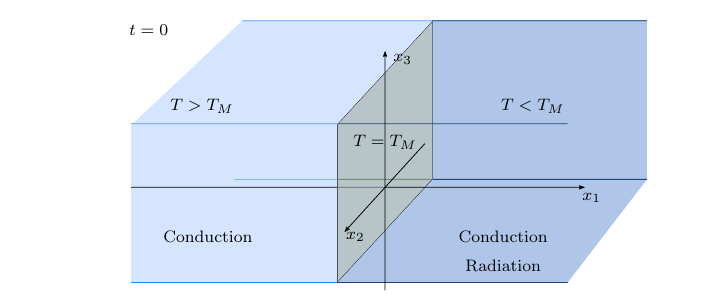}\caption{Illustration of the considered model at initial time $ t=0 $.}
		\end{center}	
	\end{figure}
	We also assume the material to satisfy local thermal equilibrium (i.e. there exists a well-defined temperature for all $ t>0 $, $ x\in\R^3 $) and we consider the case in which the scattering process is negligible. Hence, the interaction of photons with matter is described in the solid phase by the (stationary) radiative transfer equation, which under the further assumption of constant Grey approximation (i.e. $ \alpha\equiv $const.) writes
		\begin{equation}\label{RTE}
		n\cdot \nabla_x I_\nu(t,x,n)=\alpha \left(B_\nu(T(t,x)) -I_\nu(t,x,n)\right)\quad  \nu>0,\;x_1>s(t),\;n\in\Ss^2,\;t>0 ,
	\end{equation}
where $ I_\nu  $ is the radiation intensity , i.e. the energy of photons with frequency $ \nu>0 $, at position $ x\in\Omega$, moving in direction $ n\in\Ss^2 $ at time $ t>0 $, and $ B_\nu(T)= \frac{2h\nu^3}{c^2}\frac{1}{e^{\frac{h\nu}{kT}}-1} $ is the Planck distribution of a black body.

In the transport term of equation \eqref{RTE} the term containing the time derivative of $ I_\nu $, i.e. $ 	\frac{1}{c}\partial_t I_\nu(t,x,n) $ has been neglected since the characteristic time scale required in order to obtain significant changes of the temperature is much larger than the time scale in which the radiation intensity becomes stable. This is due to the fact that photons travel with the speed of light. 

In this paper it is assumed also the absence of external sources of radiation. Thus, since the photons do not interact with the liquid phase, at the interface the radiation intensity has to satisfy
	\begin{equation}\label{bnd.interface}
		I_\nu\left(t,x,n\right)=0\;\;\text{ if }x_1=s(t),\; n_1>0.
	\end{equation}
We remark that the transparency of the liquid implies that the radiation escaping the solid (i.e. traveling with direction $ n_1<0 $) passes through the liquid without interacting with it and hence without any possibility to return in the solid phase. Thus, radiation helps the solid to cool faster.

Under all these assumptions, the two-phase free boundary problem that we study in this paper is given by
	\begin{equation}\label{syst.1}
		\begin{cases}
			C_L \partial_t T(t,x_1)=K_L \partial_{x_1}^2 T(t,x_1)& x_1<s(t),\\
			C_S \partial_t T(t,x_1)=K_S \partial_{x_1}^2 T(t,x_1)-\Div\left(\int_0^\infty d\nu\int_{\Ss} dn nI_\nu(t,x,n)\right)& x_1>s(t),\\
			n\cdot \nabla_x I_\nu(t,x,n)=\alpha \left(B_\nu(T(t,x_1))-I_\nu(t,x,n)\right)& x_1>s(t),\\
			I_\nu(t,x,n)=0&x_1=s(t),\;n_1>0,\\
			T(t,s(t))=T_M&x_1=s(t),\\
			T(0,x)=T_0(x)& x_1\in\R,\\
			\dot{s}(t)=\frac{1}{L}\left(K_S \partial_{x_1} T(t,s(t)^+)-K_L\partial_{x_1} T(t,s(t)^-)\right),
		\end{cases}
	\end{equation}
where $ C_S,\;C_L $ are the volumetric heat capacities of the solid and liquid, $ K_S,\;K_L $ the conductivities of the two phases and $ L $ is the latent heat. Notice that for simplicity we are assuming that the two phases have the same constant density. For a more detailed explanation of the derivation of \eqref{syst.1} and in particular of the Stefan condition for the moving interface we refer to \cite{StefanRad1}. As remarked there, the main feature of system \eqref{syst.1} is that there is no external source of radiation and only the solid is emitting radiation. The addition of a non-trivial external source of radiation heating the solid from far away is another very interesting problem that could be studied, not only developing a well-posedness theory but also examining the possible existence of traveling waves. In this case we would consider as boundary condition
\[	\left(t,(s(t),x_2,x_3),n\right)=g_\nu(n)>0\;\;\text{ if } n_1>0.\] Moreover, in our previous article \cite{StefanRad1} we showed that reducing the radiative transfer equation to a non-local non-linear integral operator for $ T^4 $ and performing suitable rescalings, the system \eqref{syst.1} is equivalent to 
\begin{equation}\label{syst.3}
	\begin{cases}
		\partial_t T(t,x)=\kappa \partial_{x}^2 T(t,x)& x<s(t),\\
		\partial_t T(t,x)=\partial_{x}^2 T(t,x)-I_\alpha[T](t,x)& x>s(t),\\
		T(t,s(t))=T_M\\
		T(0,x)=T_0(x)& x\in\R,\\
		\dot{s}(t)=\frac{1}{L}\left( \partial_{x} T(t,s(t)^+)-K\partial_{x} T(t,s(t)^-)\right),
	\end{cases}
\end{equation}
where \[ I_\alpha[T](t,x_1)=T^4(t,x_1)-\int_{s(t)}^\infty d\eta \frac{\alpha E_1(\alpha|x_1-\eta|)}{2}T^4(t,\eta)\] for $ E_1(x)=\int_{|x|}^\infty \frac{e^-t}{t} $ being the exponential integral. Given a solution $ T $ of \eqref{syst.3}, the intensity of radiation is obtained solving by characteristics the radiative transfer equation \eqref{RTE} with boundary conditions \eqref{bnd.interface} as
\[ I_\nu(t,x,n)=\int_0^{d(x,n)}d\tau \alpha \exp\left(-\alpha \tau\right)B_\nu(T(t,x_1-\tau n_1)) \quad \text{ for }x_1>0,\]
where $ d(x,n) $ is the distance of $ x\in\R^3 $ to the interface $ \{s(t)\}\times \R^2 $ in direction $ -n\in\Ss^2 $ and it is possibly infinity. In \cite{StefanRad1} a local and global well-posedness theory for \eqref{syst.3} has been developed. Thus, a natural question that arises concerns the asymptotic behavior of the solutions to \eqref{syst.3} as $ t\to\infty $. In this paper we construct traveling waves of \eqref{syst.3} and we study their properties. Therefore, considering solutions of the form $ T(t,x)=T(x-s(t)) $ and $ s(t)=-ct $ for $ c\in\R $, in this article we study the system
\begin{equation}\label{trav.wave.1}
	\begin{cases}
		c\partial_y T_1(y)=\kappa \partial_y^2 T_1(y) &y<0\\
		c\partial_y T_2(y)=\partial_y^2 T_2(y)-T_2^4(y)+\int_0^\infty \alpha \frac{E_1(\alpha(y-\eta))}{2}T_2^4(\eta)d\eta&y>0\\
		T_2(0)=T_1(0)=T_M\\
		c=\frac{1}{L}\left(K\partial_y T_1(0^-)-\partial_yT_2(0^+)\right),\\
	\end{cases}
\end{equation} 
where we changed the variables according to $ y=x-ct $.
	\subsection{Summary of previous results}
	This paper studies a problem arising from the combination of a classical Stefan problem with the radiative transfer equation. It is therefore worth revising the most important results for these two particular problems, which as far as we know were consider together rigorously firstly in our previous paper \cite{StefanRad1}.
	
	Starting from the seminal work \cite{Stefan3} of J. Stefan, who also discovered the well-known Stefan-Boltzmann law for the total emission of a black body (cf. \cite{Stefan5}), the Stefan problem for melting of ice has been comprehensively studied in both the one-phase and the two-phase formulations, in the case of classical (i.e. strong) and weak enthalpy solutions, i.e. the weak solutions of the enthalpy formulation of the problem.
	
	The well-posedness theory for classical solutions to the Stefan problem has been considered in many works, like for instance \cite{CannonHenryKotlow,CannonPrimicerio,Friedman5,Friedman1,Friedman2,FriedmanKinderlehrer,Meirmanov,Rubinstein}, using among other methods fixed-point equations for Volterra-type integrals and the maximum principle, the Baiocchi transform, a variational inequality. Concerning the long time behavior of the one-dimensional, one-phase Stefan problem, \cite{Friedman1,Meirmanov} prove that the temperature approaches a self-similar profile as $ t\to\infty $, which is given by an error function. The works \cite{Friedman3,Friedman4} deal with the well-posedness theory of weak (enthalpy) solutions for the one and two-phase free boundary problem.
	
	Another interesting question emerging for the higher dimensional local and non-local (c.f. fractional Laplacian) Stefan problem concerns the regularity of the free boundary, which can be studied through its formulation as a parabolic obstacle problem. This has been considered in \cite{Baiocchi,Caffarelli1,Caffarelli2,Duvaut1,Figalli}.

	Finally, an important class of results addresses of the formation of supercooled liquid (i.e. liquid regions where $ T<T_M $) or superheated solid (i.e. $ T>T_M $) for the classical solutions of the Stefan problem (c.f. \cite{LaceyShillor}) as well as the creation of mushy regions (i.e. where $ T=T_M $) of positive measure for the weak enthalpy formulation of the freee boundary problem, c.f. 
	\cite{Bertsch,FasanoPrimicerio1,FasanoPrimicerio2,LaceyTayler,PrimicerioUghi,Visintin1,Visintin2}.
	
	Besides the theory of free boundary problems, this paper deals also with the theory of radiative transfer, an issue extensively studied starting from the pioneer works of Compton \cite{compton} in 1922 and of Milne \cite{Milne} in 1926. The kinetic equation describing the interaction of photons with matter is the radiative transfer equation, whose derivation and main properties can be found in \cite{Chandrasekhar,mihalas,oxenius,Rutten,Zeldovic}. 
	
	In recent years, several different problems concerning the study of the distribution of temperature due to radiation have been considered, such as well-posedness results for the stationary radiative transfer equation as in \cite{dematt42024compactness,jang}, diffusion approximation (see \cite{Golse3,Golse6,dematte2024equilibrium,dematte2023diffusion} and the references therein), the interaction of radiation and fluids (for instance in \cite{Golse1,Golse2,mihalas,Zeldovic}) and in Boltzmann gases (c.f. \cite{dematte,paper} and the reference therein). Also the study of heat transfer due to conduction and radiation as well as homogenization problems have been studied, we refer to the literature of our previous article \cite{dematte2024equilibrium}.
	
	Finally, we want to mention that free boundary problems where heat is transported by conduction and radiation have been considered numerically in engineering applications in terms of melting problems (see for instance \cite{engineering2,engineering5,engineering1,engineering3,engineering4}) and in numerical applications in the context of vaporization problems for droplets (cf. \cite{droplets3,droplets2,droplets5,droplets7,droplets8,droplets1,droplets6})

	\subsection{Main results, plan of the paper and notation}
		In this paper we will study problem \eqref{syst.1} and we will see that the addition of the radiative operator to the one-dimensional two-phase Stefan problem yields interesting phenomena which differs from the well-known results for the classical Stefan problem. 
	Specifically, we will also show that there exist traveling wave solutions for the problem \eqref{syst.1}. This is very different from the classical two-phase Stefan problem, for which self-similar profiles exist while traveling wave solutions are impossible to obtain. We will show also that the interface moves towards the liquid region, i.e. in our case $ \dot{s}=-c<0 $, implying that the traveling wave solutions exist only when the solid expands.
	\begin{theorem}\label{thm.trav.wave.0}
		There exists $ c_{\max}>0 $ such that for any $ c\in(0,c_{\max}] $ there exists a solution to the system \eqref{syst.1} (without initial condition), such that the interface satisfies $ s(t)=-ct $ and the temperature is a traveling wave defined by $ T(y)=T(x-s(t)) $ with $ T>T_M $ for $ y<0 $ and $ T<T_M $ for $ y>0 $. Moreover, for $ c<0 $ there exists no solution with $ s(t) =-ct$ such that $ T $ is a traveling wave. Finally, for any $ c\in(0,c_{\max}] $ and for $ T_M $ small enough the traveling waves are unique.
	\end{theorem}
We will see that also for $ c=0 $ traveling waves exist. However, in this case the asymptotic behavior as $ y\to\infty $ is more involved and it has not been considered yet.
 
	 The results of Section \ref{Sec.exist.trav.wave} and of Section \ref{Sec.exist.limit} will imply Theorem \ref{thm.trav.wave.0}. Specifically, in Section \ref{Sec.exist.trav.wave} we will show the existence of traveling wave solutions in the case of negative speed of the interface, i.e. when the ice is expanding. While the traveling waves in the liquid are given by the well-known solution to the ODE $ \ddot{y}=\lambda \dot{y} $ for $ x<0 $, the existence of traveling wave solutions in the solid is more involved. By a variational argument we will prove the existence of such traveling waves (cf. subsection \ref{subs.existenec}), which will be shown to be monotonically increasing with respect to the melting temperature ((cf. subsection \ref{subs.monoton})). In Subsection \ref{subs.small} we will also show that for very small melting temperatures there exists a unique strictly positive traveling wave solution, which also converges with exponential rate to a positive constant as $ x\to\infty $. In Section \ref{Sec.exist.limit} the analysis of the traveling wave is carried on. In particular several applications of the maximum principle will be used together with blow-up limits, Liouville-type theorems, and Harnack-type arguments in order to show that the traveling waves have a limit as $ x\to\infty $.  Finally, in Section \ref{Sec.pict.} we will conclude this paper using asymptotic arguments with a formal picture of the long time asymptotic of the solutions to \eqref{syst.1} for arbitrary values of $ \lim\limits_{y\to-\infty}T(y)=T(-\infty) $ and $ \lim\limits_{y\to\infty}T(y)=T(\infty) $. \\
	 
	Throughout this article we will denote by $ C^{k,\beta} (U)$, where $ U\subseteq\R $ is possibly unbounded, the space of $ k $-times continuous differentiable functions $ f $ with
	 \[\Arrowvert f\Arrowvert_{k,\beta}=\max\limits_{0\leq j\leq k}\left(\sup\limits_{U}\left|\partial_x^j f\right| \right)+\sup\limits_{x,y\in U}\frac{\left|\partial^k_x f(x)-\partial_x^kf(y)\right|}{|x-y|^\beta}<\infty.\]
	 Notice that $ f\in C^{k,\beta}(U) $ has all $ k $ derivatives bounded.   
		\section{Existence of traveling wave solutions}\label{Sec.exist.trav.wave}
	 In the following section we construct traveling wave solutions solving \eqref{trav.wave.1} and we will prove some important properties satisfied by such functions. First of all, we will see that the traveling waves propagate necessarily with negative velocity. Hence, the interface is moving towards the liquid part and the ice is expanding. This behavior is intriguing at a first glance. However, it can be expected. Indeed, while in the liquid the heat is transferred only by conduction, in the solid the heat is also transferred by radiation. Since a radiative source is absent in this problem, the radiation escaping from the solid is helping the ice to cool faster. 
		
		Recall that the system \eqref{trav.wave.1} has been obtained considering solutions to the original problem \eqref{syst.3} of the form $ T(t,x)=T(t,x-s(t)):=T(y) $ and $ s(t)=-ct $ with $ c\in\R $. First of all we see that in order to obtain the existence of bounded solutions to the problem \eqref{trav.wave.1} $ c $ must be positive, thus since $ \dot{s}(t)=-c<0 $ the ice is expanding. Indeed, if $ c<0 $ then the temperature of the liquid should satisfy
		\begin{equation}\label{trav.wave.2}
			\begin{cases}
				c\partial_y T_1(y)=\kappa \partial_y^2 T_1(y) &y<0\\
				T_1(0)=T_M\\
			\end{cases}
		\end{equation} 
		and hence 
		\begin{equation}\label{unbnd.trav}
			T_1(y)=T_M+\left|\frac{A}{c}\right|\kappa \left(e^{-\frac{|c|}{\kappa}y}-1\right)\to \infty \quad\text{ as }y\to-\infty.
		\end{equation}
		Let thus $ c>0 $. In the next subsections we will prove the following theorem.
		\begin{theorem}\label{thm.trav.wave.all}
			For $ c<0 $ the problem \eqref{trav.wave.1} does not admit any bounded solution. However, there exists $ c_{\max}>0 $ such that for any $ c\in(0,c_{\max}] $ there exists traveling waves $ T_1,T_2 $ solving \eqref{trav.wave.1}. Moreover, for $ c\in(0,c_{\max}] $ the solutions satisfy $ T_1(y)>T_M $ for $ y<0 $ and $ 0<T_2(y)<T_M $ for $ y>0 $, and the limits $ \lim\limits_{y\to-\infty}T_1(y) $ and $ \lim\limits_{y\to\infty} T_2(y) $ exist.
			\begin{proof}
				As we have seen in \eqref{trav.wave.2} and in \eqref{unbnd.trav}, if $ c<0 $ the problem \eqref{trav.wave.1} does not have any bounded solution. Thus, we set $ c>0 $ and we see that for any $ c $ and any $ \alpha\in\R $ the solution to
				\[\begin{cases}
					c\partial_y T_1(y)&y<0\\
					T_1(0)=T_M\\
					\partial_yT_1(0)=-A
				\end{cases}\]
				is given by
				\[T_1(y)=T_M+\frac{A}{c}\kappa \left(1-e^{\frac{c}{\kappa}y}\right).\]
				Moreover, $ \lim\limits_{y\to-\infty}T_1(y)=T_M+\frac{A}{c}\kappa $ . Since $ T_1 $ describes the temperature in the liquid, we are interested only in $ A\geq0 $.\\
				
				In the following subsections we will study
				\begin{equation}\label{trav.wave.3}
					\begin{cases}
						\partial_y^2 f(y)-c\partial_y f(y)-f^4(y)=-\int_0^\infty \alpha \frac{E_1(\alpha(y-\eta))}{2}f^4(\eta)d\eta&y>0\\
						f(0)=T_M\\
						f\geq 0\\
					\end{cases}
				\end{equation} We will prove the existence of functions $ f\in  C^{2,1/2}(\R_+) $ solving the problem \eqref{trav.wave.3}. We will show also that there exists $ c_{\max}>0 $ such that $ \partial_y f(0^+)\leq -Lc$ for all $ 0<c<c_{\max} $. Then for $ c\in(0,c_{\max}) $ and $ A=-\frac{L c+\partial_y T_2(0^+)}{K} $ the functions $ T_1(y)=T_M+\frac{A}{c}\kappa \left(1-e^{\frac{c}{\kappa}y}\right) $ and $ T_2:=f $ are traveling waves solving \eqref{trav.wave.1}.
			\end{proof}
		\end{theorem} 
		Before moving to the existence theory for the solutions to \eqref{trav.wave.3} we do the following remark. It is enough to prove that the traveling wave solutions in the solid exists, that they are bounded from below and have a limit only for $ \alpha=1 $ and $ c>0 $. Indeed, let $ \alpha>0 $, $ c>0 $ and $ T_M>0 $ and let $ f $ solve \eqref{trav.wave.3}. Then the function $ \tilde{f}$ defined by \[f(y):=\alpha^{2/3}\tilde{f}(\alpha y)=\alpha^{2/3}\tilde{f}(\eta)\]
		satisfies the following equation
		\begin{equation}\label{trav.wave.3.0}
			\begin{cases}
				\partial_\eta^2 \tilde{f}(\eta)-c\alpha^{-1}\partial_\eta \tilde{f}(\eta)-\tilde{f}^4(\eta)=-\int_0^\infty  \frac{E_1(\eta-\xi)}{2}\tilde{f}^4(\xi)d\xi&\eta>0\\
				\tilde{f}(0)=T_M\alpha^{-2/3}\\
				\tilde{f}\geq 0\\
			\end{cases}
		\end{equation} 
		This is true since $ \partial_y f(y)=\alpha^{5/3}\partial_\eta\tilde{f}(\eta) $ as well as $ \partial_y^2 f(y)=\alpha^{8/3}\partial_\eta^2\tilde{f}(\eta) $. Notice also that $ f $ and $ 	\tilde{f} $ have the same regularity. Moreover, using that $ \eta=\alpha y $ and changing the variable $ \alpha\xi=z $ we have
		\[\int_0^\infty \alpha \frac{E_1(\alpha(y-\xi))}{2}f^4(\xi)d\xi=\int_0^\infty \alpha \frac{E_1(\eta-\alpha\xi)}{2}\alpha^{8/3}\tilde{f}^4(\alpha\xi)d\xi=\alpha^{8/3}\int_0^\infty \frac{E_1(\eta-z)}{2}\tilde{f}^4(z)dz.
		\]
		Hence, we see that defining $ E(x)=\frac{E_1(x)}{2} $ it is enough to consider the solutions to
		\begin{equation}\label{trav.wave.3.1}
			\begin{cases}
				\partial_y^2 f(y)-c\partial_y f(y)-f^4(y)=-\int_0^\infty E(y-\eta)f^4(\eta)d\eta&y>0\\
				f(0)=T_M\\
				f\geq 0\\
			\end{cases}
		\end{equation}
		
		\subsection{Existence of traveling wave solutions for $ y>0 $}\label{subs.existenec}
		Before proving the existence of traveling wave solutions for $ y>0 $ we prove the following technical proposition.
		\begin{prop}\label{prop.variational}
			Let $ c>0 $ and $ g\in C^{0,1/2}(\R_+) $ with $ -A^4<g\leq0 $ for some $ A>0 $. Let also \begin{equation*}\label{def.A.1}
				\mathcal{A}_{A,c}=\left\{f\geq 0\text{ measurable s.t. }f\in W^{1,2}\left(e^{-cy}dy,\R_+\right)\cap L^5\left(e^{-cy}dy,\R_+\right),\;f(0)=A>0\right\}.
			\end{equation*} Then the functional
			\[I_g[f]=\int_0^\infty e^{-cy}\left(\frac{(\partial_y f(y))^2}{2}+\frac{f(y)^5}{5}+g(y)f(y)\right)dy \]
			has a unique minimizer $ f\in \mathcal{A}_{A,c} $. Moreover, $ 0<f\leq A $ for $ y\in[0,\infty) $. Finally, $ f\in C^{2,1/2}(\R_+)$ solves the ODE
			\[\partial_y\left(e^{-cy}\partial_y f(y)\right)=\left(f^4(y)+g(y)\right)e^{-cy}\] and satisfies the bounds
			\[|f'(y|\leq \frac{A^4}{c},\;|f''(y)|\leq A^4\text{ and }[f'']_{1/2}\leq\max\left\{2A^4,2A^4c+\frac{4A^7}{c}+[g]_{1/2}\right\}.\]
			\begin{proof}
				Let us define the measure $ \mu $ given by the density $ d\mu(y)=e^{-cy}dy $. First of all we notice that if $ f\in W^{1,2}(\mu,\R_+)\cap L^5(\mu,\R_+) $, then $ fe^{-c/2 y}\in W^{1,2}(\R_+) $. Thus, by Morrey's embedding theorem $ fe^{-c/2y}\in C^{0,1/2}(\R_+) $. Hence, if $ f\in\mathcal{A} $, then $ f $ is continuous. This implies that the condition for $ f\in\mathcal{A}_{A,c} $ to be $ f\geq 0 $ holds everywhere in $ \R_+ $ as well as the boundary condition $ f(0)=A $, which for general functions in $ W^{1,2}(\mu,\R_+) $ is to be intended as trace condition, holds pointwise. These observations yield that $ \mathcal{A}_{A,c} $ is a closed and convex subset of $ W^{1,2}(\mu,\R_+)\cap L^5(\mu,\R_+) $. We also remark that the trace operator for $ \partial\R_+=\{0\}$ is a continuous operator with respect to the norm $ \Arrowvert \cdot \Arrowvert_{W^{1,2}(\mu\R_+)} $.
				
				Further, we notice that $ I_g $ is well-defined for $ f\in\mathcal{A}_{A,c}$ with 
				\[\left|I_g[f]\right|\leq \frac{1}{2}\Arrowvert f\Arrowvert^2_{W^{1,2}(\mu)}+\frac{1}{5}\Arrowvert f\Arrowvert_{L^5(\mu)}^5+\frac{1}{2c}\Arrowvert g\Arrowvert_\infty^2.\]
				Moreover, $ I_g[f] $ is bounded from below and coercive. Indeed, using both Young's inequality
				\[|g(y)|f(y)\leq \frac{8\cdot 2^{1/4}}{5c}|g(y)|^{5/4}+\frac{f^5(y)}{10}\] and H\"older's inequality
				\[\left(\int_0^\infty e^{-cy}|f(y)|^2dy\right)^{5/2}\leq \frac{1}{c^{3/2}}\Arrowvert f\Arrowvert_{L^5(\mu)}^5\] 
				we estimate
				\[I_0[f]\geq \min\left\{\frac{c^{3/2}}{10},\frac{1}{2}\right\}\left(\Arrowvert\partial_{y}f\Arrowvert_{L^2(\mu)}^2+\Arrowvert f\Arrowvert_{L^2(\mu)}^5\right)+\frac{1}{10}\Arrowvert f\Arrowvert_{L^5(\mu)}^5\to \infty \quad\text{ as }\Arrowvert f\Arrowvert_\mathcal{A}\to\infty\] if $ g\equiv 0 $ and 
				\[I_g[f]\geq \min\left\{\frac{c^{3/2}}{20},\frac{1}{2}\right\}\left(\Arrowvert\partial_{y}f\Arrowvert_{L^2(\mu)}^2+\Arrowvert f\Arrowvert_{L^2(\mu)}^5\right)+\frac{1}{20}\Arrowvert f\Arrowvert_{L^5(\mu)}^5-\frac{8\cdot 2^{1/4}}{5c}\Arrowvert g\Arrowvert_\infty^{4/5}\to \infty \quad\text{ as }\Arrowvert f\Arrowvert_\mathcal{A}\to\infty\] if $ g\not\equiv 0 $. 
				Moreover, $ I_0[f]\geq 0 $ as well as $ I_g[f]\geq -\frac{8\cdot 2^{1/4}}{5c}\Arrowvert g\Arrowvert_\infty^{4/5} $.\\
				
				Therefore, there exists a bounded minimizing sequence $ f_k\in\mathcal{A}_{A,c} $ such that $ I_g[f_k]\to \inf\limits_{f\in\mathcal{A}_{A,c}}I[f] $ as $ k\to\infty $. The boundedness of this sequence, the uniqueness of the weak and strong limit as well as the fact that $ L^2(\mu)\subset L^{5/4}(\mu)=\left(L^5(\mu)\right)^* $ imply the existence of a common subsequence $ f_{k_j} $ such that 
				\[f_{k_j}{\xrightarrow [\text{ptw. a.e.}]{L^2(\mu)}} f\in L^2(\mu) \quad \text{ and }\quad f_{k_j}{\xrightarrow [ \text{weak }L^5(\mu)]{\text{weak }W^{1,2}(\mu)}} f\in W^{1,2}(\mu)\cap L^5(\mu)\quad \text{ as }j\to\infty.\]
				The closedness and the convexity of $ \mathcal{A}_{A,c} $ imply also $ f\in\mathcal{A}_{A,c} $. Moreover, the pointwise convergence almost everywhere and the weak lower semicontinuity of the $ L^2 $ norm imply the weak lower semicontinuity of the functional $ I_g $. Hence, $ f $ is a minimizer of $ I_g $, i.e. $ I_g[f]=\inf\limits_{f\in\mathcal{A}_{A,c}}[f] $. In addition to that, since the functional $ I_g $ is strictly convex for non-negative functions, the minimizer is unique.
				
				We remark that $ f\in C^{0,1/2}_{\text{loc}}(\R_+) $ with $ f(y)\geq 0 $ for $ y\geq 0 $ and $ f(0)=A $. Next we prove that $ f\leq T_M $ if $ g\equiv 0 $ and that $ f\leq 5A $ if $ g\not\equiv 0 $. Both claims are a consequence of the uniqueness of the minimizer of $ I_g $ in $ \mathcal{A}_{A,c} $. If $ g\equiv 0 $ let us consider $ h_0=\min\left\{f,\;A\right\}\in\mathcal{A}_{A,c} $ , since the minimum of two Sobolev functions is a Sobolev function. Then the functional $ I_0 $ acting on $ h_0 $ gives
				\begin{equation*}
					I_0[h_0]=\int_0^\infty e^{-cy}\left(\frac{|\partial_y f|^2}{2}\mathds{1}_{\{f\leq A\}}+\frac{h_0^5}{5}\right)dy\leq \int_0^\infty e^{-cy}\left(\frac{|\partial_y f|^2}{2}+\frac{f^5}{5}\right)dy=I_0[f]=\inf\limits_{f\in\mathcal{A}_{A,c}}I_0[f] ,
				\end{equation*}
				where we used that $ 0\leq h_0\leq f $. By uniqueness we conclude $ 0\leq f\leq A $. In a similar way, if $ g\not\equiv 0 $ we consider  $ h_1=\min\{f,5A\} $. It is not difficult to see that 
				\[
				\mathds{1}_{\{f>5A\}}\left(\frac{h_1^5}{5}-|g|h_1\right)=\mathds{1}_{\{f>5A\}}\left(5^3 A^4-|g|\right)5A<\mathds{1}_{\{f>5A\}}\left(\frac{f^4}{5}-|g|\right)5A<\mathds{1}_{\{f>5A\}}\left(\frac{f^5}{5}-|g|f\right).
				\]
				For this chain of inequalities we used the definition of $ h_1 $ and the fact that $ |g|<A^4 $. Therefore $ 0<\left(5^3 A^4-|g|\right)<\left(\frac{f^4}{5}-|g|\right) $ in the set $ \{f>5A\} $.
				We conclude 
				\[I_g[h_1]=\int_0^\infty e^{-cy}\left(\frac{|\partial_y f|^2}{2}\mathds{1}_{\{f\leq 5A\}}+\frac{h_1^5}{5}+gh_1\right)dy\leq \int_0^\infty e^{-cy}\left(\frac{|\partial_y f|^2}{2}+\frac{f^5}{5}+gf\right)dy=I_g[f]=\inf\limits_{f\in\mathcal{A}_{A,c}}I_g[f].\]
				Hence, $ f=h\leq 5A $. These results show that $ f\in C_b(\R_+) $.\\
				
				We now study the Euler-Lagrange equations associated to the functional $ I_g $. It turns out that the minimizer $ f $ is the weak solution of the following inequality
				\begin{equation}\label{weak.ineq.all.0}
					-\partial_y\left(e^{-cy}\partial_y f\right)+e^{-cy}f^4(y)+e^{-cy}g\geq 0.
				\end{equation}
				Hence, $ f $ satisfies
				\begin{equation}\label{weak.ineq.all.1}
					0\leq \int_0^\infty e^{-cy}\left(\partial_y f\partial_y \psi+f^4(y)\psi(y)+g(y)\psi(y)\right)dy,
				\end{equation}
				for all $ \psi\geq 0 $, $ \psi\in C_c^\infty(\R_+) $ or also $ \psi\in W^{1,2}_0(\mu)\cap L^5(\mu)$. Moreover, on the open set $ \{f>0\} $ the minimizer $ f $ is a weak solution of the equation
				\begin{equation}\label{weak.eq.all.0}
					-\partial_y\left(e^{-cy}\partial_y f\right)+e^{-cy}f^4(y)+e^{-cy}g(y)= 0.
				\end{equation}
				Indeed, on the open set $ \{f>0\} $ for any $ \psi\in C_c^\infty(\{f>0\}) $ the function $ f+\eps\psi\in\mathcal{A} $ for $ \eps>0 $ small enough. Hence
				\begin{equation}\label{weak.eq.all.1}
					0=\partial_\eps I[f+\eps\psi]\big|_{\eps=0}= \int_0^\infty e^{-cy}\left(\partial_y f\partial_y \psi+f^4(y)\psi(y)+g(y)\psi(y)\right)dy,
				\end{equation}
				for all $ \psi\in C_c^\infty(\{f>0\}) $ or also $ \psi\in W^{1,2}_0(\mu,\{f>0\})\cap L^5(\mu,\{f>0\})$. We remark that equations \eqref{weak.ineq.all.0}-\eqref{weak.eq.all.1} hold for both $ g\equiv 0 $ and $ g\not\equiv 0 $.
				
				We aim to show that actually the minimizer $ f $ is a strong solution to \eqref{weak.eq.all.0} in the whole real line. To this end we will first show that $ \{f>0\}=\R_+ $, which implies that $ f $ is a weak solution of \eqref{weak.eq.all.0} in $ \R_+ $, and finally we will use elliptic regularity theory for \eqref{weak.eq.all.0}. \\
				
				Let us assume that $ \{f>0\}\subsetneq\R_+ $. Then there exists $ a \in\R_+ $ such that $ f(y)>0 $ for all $ y<a $ and $ f(a)=0 $. We have to consider two cases: first the case where $ f(y)\equiv0 $ in an interval $ (a,a+r) $ for some $ r>0 $ and second the case where $ f(y)\not\equiv0 $ on the interval $ (a,a+r) $ for any $ r>0 $. \\
				
				Let us assume first that there exists $ r>0 $ such that $ f(y)=0 $ for all $ y\in(a,a+r) $.  Since $ f $ is continuous there exists $ 0<\eps<\min\left\{r,\frac{c}{2}\right\} $ small enough such that $ f(a-\eps)=\delta\ll 1 $ as well as $ f(a+\eps)=0 $. Let us define for $ y\in[a-\eps,a+\eps] $ the following function
				\[\bar{f}(y)=\delta\left(1-\frac{y-(a-\eps)}{2\eps}\right).\]
				It is easy to see that $ 0<\bar{f}<\delta<1 $ for $ y\in(a-\eps,a+\eps) $, $ \bar{f}(a-\eps)=f(a-\eps) $ as well as $ \bar{f}(a+\eps)=f(a+\eps)=0 $. Moreover, $ \bar{f}(a)=\frac{\delta}{2}>0 $. Finally, since $ \bar{f}^4\leq \delta $ and $ 2\eps<c $, an easy computation shows 
				\begin{equation}\label{weak.ineq.10}
					-\bar{f}''(y)+c\bar{f}'(y)+\bar{f}'(y)=-\frac{c\delta}{2\eps}+\bar{f}^4(y)\leq \delta\left(1-\frac{c}{2\eps}\right)<0.
				\end{equation}
				Thus, $ -\left(e^{-cy}\bar{f}'(y)\right)'+e^{-cy}\bar{f}^4(y)<0 $. Since $ f(a)=0<\bar{f}(a) $, there exists an interval $ (y_0,y_1)\subseteq (a-\eps,a+\eps) $ such that $ f(y_0)=\bar{f}(y_0) $, $ f(y_1)=\bar{f}(y_1) $ and $ f(y)<\bar{f}(y) $ for $ y\in(y_0,y_1) $. Using the weak maximum principle we show now that this is not possible. Therefore, we test \eqref{weak.ineq.10} with a suitable test function $ \psi\geq 0 $. Let us consider the smooth solution to 
				\begin{equation}\label{positivity.small.interval3.all}
					\begin{cases}
						\partial_y^2\psi(y)-c\partial_y \psi(y)=-1&(y_0,y_1);\\
						\psi(y_0)=\psi(y_1)=0\\
					\end{cases}
				\end{equation} 
				The solution is given by the explicit formula $ \psi(y)=\frac{y-y_0}{c}-\frac{y_1-y_0}{c}\frac{e^{c(y-y_0)}-1}{e^{c(y_1-y_0)}-1} $. By a simple application of the maximum principle we see that $ \psi>0 $ in $ (y_0,y_1) $. Indeed, if $ \psi $ would have a minimum at $ y^*\in(y_0,y_1) $ on that point $ \psi   $ would not solve the equation, since $ \psi''(y^*)-c\psi'(y^*)\geq 0 $. Hence, let us consider $ \bar \psi $ as the extension by $ 0 $ of $ \psi $ in the whole positive real line, i.e. \begin{equation*}\label{bar.psi.all}
					\bar \psi(y)=\begin{cases}
						\psi(y)&y\in(y_0,y_1)\\ 0& \text{else}.
					\end{cases}
				\end{equation*}
				Clearly $ \bar \psi \in W_0^{1,2}(\R_+,\mu)\cap L^5(\R_+,\mu) $. Then, 
				\[-\left(e^{-cy}\bar{f}'(y)\right)'\bar \psi(y)+e^{-cy}\bar{f}^4(y)\bar \psi(y)\leq0,\]
				where we used that $ \bar{\psi}\equiv 0 $ on $ \R_+\setminus(a-\eps,a+\eps) $. Therefore, using also that the weak derivative of $ \bar{\psi} $ is supported also on $ [a-\eps,a+\eps] $ we obtain
				\begin{equation}\label{weak.bar.all.f}
					\int_0^\infty e^{-cy}\left(\bar{f}'(y)\bar{\psi}'(y)+\bar{f}^4(y)\bar{\psi}(y)\right)dy\leq 0.
				\end{equation}
				Hence, using \eqref{weak.ineq.all.1}, \eqref{weak.bar.all.f} and the definition of $ \bar \psi $ we have
				\begin{equation}\label{weak.max.all.1}
					\begin{split}
						0\leq&\int_0^\infty e^{-cy}\left(\partial_y (f-\bar f)\partial_y\bar \psi+\left(f^4-\bar f^4\right)\bar \psi+g\bar{\psi}\right)dy\\
						=& \int_{y_0}^{y_1}(f-\bar f) \partial_y\left( -e^{-cy}\partial_y \psi\right)+e^{-cy}\left(f^4-\bar f^4\right) \psi+e^{-cy}g{\psi} dy\\
						= &\int_{y_0}^{y_1} e^{-cy}\left((f-\bar f)\left(-\partial_y^2 \psi+c\partial_y \psi \right)+\left(f^4-\bar f^4\right) \psi+g{\psi}\right)dy<0
					\end{split}
				\end{equation}
				where we used also that $ \left(f-\bar f\right)\left.\right|_{\{y_0,y_1\}}=0 $, $ 0\leq f<\bar f $ on $ (y_0,y_1) $ as well as $ g\leq 0 $. This contradiction implies that $ f(y)\geq \bar{f}(y)>0 $ on $ (a-\eps,a+\eps) $. But since we assumed $ f(a)=0<\bar{f}(a) $ we conclude that there cannot exist any $ r>0 $ such that $ f(y)=0 $ for $ y\in(a,a+r) $.\\
				
				Hence, we assume that $ f(y)\not\equiv 0 $ for $ y\in(a,a+r) $ and $ r>0 $. Since $ f(a)=0 $ by continuity there exist $ 0<\eps_1,\eps_2<\min\{r,\frac{c}{4}\} $ small enough such that $ f(a-\eps_1)=\delta\ll1 $ and $ f(a+\eps_2)=\frac{\delta}{2} $. We then define for $ y\in[a-\eps_1,a+\eps_2] $ the function	\[\bar{f}(y)=\delta\left(1-\frac{y-(a-\eps_1)}{2(\eps_1+\eps_2)}\right).\]
				Also in this case $ \bar{f} $ satisfies $ 0<\frac{\delta}{2}<\bar{f}<\delta<1 $ for $ y\in(a-\eps_1,a+\eps_2) $, $ \bar{f}(a-\eps_1)=f(a-\eps_2) $, $ \bar{f}(a+\eps_2)=f(a+\eps_2)$, $ \bar{f}(a)\geq\frac{\delta}{2}>0 $, as well as
				\[-\bar{f}''(y)+c\bar{f}'(y)+\bar{f}'(y)=-\frac{c\delta}{2(\eps_1+\eps_2)}+\bar{f}^4(y)\leq \delta\left(1-\frac{c}{2(\eps_1+\eps_2)}\right)<0.\]
				We now argue as in the case $ f(a+\eps_2)=0 $. As we have seen before, since $ f(a)=0<\bar{f}(a) $, there exists an interval $ (y_0,y_1)\subseteq (a-\eps,a+\eps) $ such that $ f(y_0)=\bar{f}(y_0) $, $ f(y_1)=\bar{f}(y_1) $ and $ f(y)<\bar{f}(y) $ for $ y\in(y_0,y_1) $. Then, testing $ f-\bar{f} $ against the function $ \bar{\psi} $ defined as the zero extension of $ \psi $ in \eqref{positivity.small.interval3.all} we obtain the following contradiction as for \eqref{weak.max.all.1}
				\begin{equation*}\label{weak.max.all.2}
					\begin{split}
						0\leq&\int_0^\infty e^{-cy}\left(\partial_y (f-\bar f)\partial_y\bar \psi+\left(f^4-\bar f^4\right)\bar \psi+g\bar{\psi}\right)dy
						<0.
					\end{split}
				\end{equation*}
				This contradiction yields that $ \{f>0\}=\R_+ $. Thus, $ f $ is a weak solution to \eqref{weak.eq.all.0}.
				
				In the case where $ g\not \equiv 0 $, we proved that $ f\leq 5A $. We now prove that also $ f\leq A $ holds. To this end we consider for $ R>0 $ the function $ \phi_R(y) $ defined by $ \phi_R(y)=A+ 4A e^{c(y-R)}\geq A $. We see that $ \phi_R(0)> A=f(0) $ as well as $ \phi_R(R)=5A\geq f(R) $. By continuity we know that there exists some $ x_0\in[0,R] $ such that $ \min\limits_{[0,R]}\psi_R-f=\psi_R(x_0)-f(x_0) $. Hence, let us assume that $ \min\limits_{[0,R]}\psi_R-f=\psi_R(x_0)-f(x_0)<0 $. Since $ \phi_R-f\left.\right|_{\{0,R\}}\geq 0 $, there exists an interval $ x_0\in(a,b)\subset [0,R] $ in which $ \phi_R-f<0 $ and $ \phi_R-f(a)=\phi_R(b)-f(b)=0 $. We also see that $ \phi_R $ is a supersolution for the operator $ \mathcal{L}[\phi]=-\phi''+c\phi'+\phi^4+g $ on $ [0,R] $. Indeed
				\[
				\mathcal{L}[\phi_R]=\phi_R^4+g>A^4-A^4=0.
				\]
				Let us consider once again the zero extension $ \bar{\psi} $ of the function $ \psi>0 $ given by \eqref{positivity.small.interval3.all} on the interval $ (a,b) $. Then we see that 
				\[
				0\leq \int_0^\infty e^{-cy}\left(\partial_y \phi_R\partial_y  \bar{\psi}+\phi_R^4 \psi+g(y) \bar{\psi}(y)\right)dy.\]
				Therefore we obtain the following contradiction using once more that $ \left(f-\phi_R\right)\left.\right|_{\{a,b\}}=0 $, that $ 0< \phi_R<f $ on $ (a,b) $, and that $ f $ is a weak solution solving \eqref{weak.eq.all.0} 
				\begin{equation*}
					\begin{split}
						0\leq &
						\int_0^\infty e^{-cy}\left(\partial_y (\phi_R-f)\partial_y  \bar{\psi}+\left(\phi_R^4-f^4\right)  \bar{\psi}\right)dy\\
						=&\int_a^be^{-cy}\left((\phi_R-f)\left(-\partial_y^2 \psi+c\partial_y\psi\right)+\left(\phi_R^4-f^4\right) \psi \right)dy\\
						=&\int_a^b e^{-cy}\left((\phi_R-f)+\left(\phi_R^4-f^4\right) \psi \right)dy<0.
					\end{split}
				\end{equation*}
				Hence, for any $ y\in [0,R] $ we have $ f(y)\leq A +4Ae^{c(y-R)} $. Letting now $ R\to \infty $, we conclude that $ 0\leq f\leq A$.\\
				
				We finish the proof of Proposition \ref{prop.variational} showing that $ f $ is also a strong solution to \eqref{weak.eq.all.0}. This can be proved using the elliptic Schauder regularity. Indeed, since $ f\in\mathcal{A}_{A,c} $ is bounded and continuous, we have that $ f\in W^{1,2}(\mu)\cap L^\infty(\R_+) $. Hence, $ f\in W_{\text{loc}}^{1,2}(\R_+,dy)\cap L^\infty(\R_+)  $, so that also $ f^4e^{-cy}\in W_{\text{loc}}^{1,2}(\R_+,dy)\cap L^\infty(\R_+) $. Morrey's embedding theorem implies that $ f\in C_{\text{loc}}^{0,1/2}(\R_+) $, which yields also $ f^4e^{-cy}\in C_{\text{loc}}^{0,1/2}(\R_+) $. Applying now the elliptic regularity theory to the equation \eqref{weak.eq.all.0} we obtain that $ f\in C_{\text{loc}}^{2,1/2}(\R_+) $ since also $ ge^{-cy}\in C^{0,1/2}(\R_+) $. Thus, $ f\in C^2(\R_+) $ is a strong solution to \eqref{weak.eq.all.0}.\\
				
				We now show that $ f $ has also bounded first and second derivative. This is due to the fact that also $ f'\in W^{1,2}(\mu) $. Indeed, $$ \int_0^\infty e^{-cy}\left(|f''|^2+|f'|^2\right)dy\leq  \int_0^\infty e^{-cy}\left(|f^4+cf'+g|^2+|f'|^2\right)dy\leq C(A, c)\left(\Arrowvert f\Arrowvert_{W^{1,2}(\mu)}+\frac{A^8}{c}\right).$$ Hence, $ e^{-\frac{c}{2}y}f'\in W^{1,2}(\R_+,dy) $, which implies that $ e^{-cy}(f')^2 $ is bounded since its derivative $ 2e^{cy}f'f''-ce^{-cy}(f')^2 $ is integrable. Thus, the consequent boundedness of  $ e^{-\frac{c}{2}y}\left|f'\right| $ implies that \begin{equation}\label{limit.f_1.all}
					\lim\limits_{y\to\infty}e^{-cy}|f'|(y)=0.
				\end{equation}
				Since $ f $ solves \eqref{weak.eq.all.0}, using \eqref{limit.f_1.all} and integrating in $ (y,\infty) $ we obtain the desired estimate
				\begin{equation*}\label{derivative.f_1.all}
					|f'|(y)\leq e^{cy}\int_y^\infty e^{-c\xi}\left|f^4(\xi)+g(\xi)\right|d\xi\leq \frac{A^4}{c}.
				\end{equation*}
				Moreover, multiplying \eqref{weak.eq.all.0} by $ e^{cy} $ we conclude that $ f $ is a $ C^2 $-solution to
				\[f''-cf'=f^4+g\text{ on }\R_+.\]
				This yields the boundedness of the second derivative of $ f $ as 
				\begin{equation*}\label{s.derivative.f_1.all}
					|f''|(y)\leq c|f'|(y)+\left|f^4(y)+g(y)\right|\leq A^4,
				\end{equation*}
			where we used also $ 0\leq f^4\leq A^4 $ and $ -A^4\leq g\leq 0 $.
				These estimates imply that $ f\in C^{1,1}(\R_+) $ with bounded first and second derivatives. 
				Since $ cf'+f^4+g\in C^{0,1/2}(\R_+) $ we conclude that $ f\in C^{2,1/2}(\R_+) $ with H\"older seminorm bounded by
				\begin{equation*}\label{hoelder.derivative.f_1.all}
					[f'']_{1/2}\leq \max\left\{2\Arrowvert f''\Arrowvert_\infty, c\Arrowvert f''\Arrowvert_\infty+4\Arrowvert f\Arrowvert_\infty^3\Arrowvert f'\Arrowvert_\infty+[g]_{1/2}\right\}\leq\max\left\{2A^4, 2A^4c+\frac{4A^7}{c}+[g]_{1/2}\right\}.
				\end{equation*}
			\end{proof}
		\end{prop}
		Let us now consider the sequence $ f_n\in  C^{2}(\R_+)  $ with $ f_n\geq 0 $ such that
		\begin{equation}\label{trav.wave.4}
			\begin{cases}
				\partial_y^2 f_{n+1}(y)-c\partial_y f_{n+1}(y)-f_{n+1}^4(y)=-\int_0^\infty E(y-\eta)f_n^4(\eta)d\eta&y>0;\quad n\geq 1\\
				f_0=0&n=0\\
				f_{n+1}(0)=T_M\\
				f_{n+1}\geq 0\\
			\end{cases}
		\end{equation} 
		We prove the following theorem
		\begin{theorem}\label{thm.existence.all}
			Let $ T_M,c>0 $. Then there exists a solution $ f\in C^{2,1/2}(\R_+)  $ with $ f>0 $ at the interior of $ \R_+ $ solving \eqref{trav.wave.3.1}. Moreover, $ f $ is obtained as the limit of the monotone increasing bounded sequence \[0\leq f_1\leq f_2\leq ...\leq f_n\leq f_{n+1}\leq ...\leq T_M\] 
			with $ \left(f_n\right)_{n\in\mathbb{N}} \in C^{2,1/2}(\R_+) $ with $ \Arrowvert f_n\Arrowvert_{2,1/2} $ uniformly bounded and with $ f_n>0 $ in the interior of $ \R_+ $ solving the recursive system \eqref{trav.wave.4}.
			\begin{proof}
				We start considering the function $ f_1 $ solving the problem
				\begin{equation}\label{trav.wave.5}
					\begin{cases}
						\partial_y^2 f_1(y)-c\partial_y f_1(y)-f_1^4(y)=0&y>0;\\
						f_1(0)=T_M\\
						f_1\geq 0\\
					\end{cases}
				\end{equation} 
				The differential equation is equivalent to the elliptic ODE \[\left(e^{-cy}f'\right)'=e^{-cy}f^4.\]
				Hence, we consider the minimization problem of the functional 
				\[I_0[f]=\int_0^\infty e^{-cy}\left(\frac{(\partial_y f(y))^2}{2}+\frac{f(y)^5}{5}\right)dy\]
				on the set
				\begin{equation}\label{def.A}
					\mathcal{A}_{T_M,c}=\left\{f\geq 0\text{ measurable s.t. }f\in W^{1,2}\left(e^{-cy}dy,\R_+\right)\cap L^5\left(e^{-cy}dy,\R_+\right),\;f(0)=T_M\right\}.
				\end{equation}
				Proposition \ref{prop.variational} shows that there exists a unique  $ f_1\in \mathcal{A}_{T_M,c}$ minimizing the functional $ I_0 $. Moreover, $ f_1 $ solves \eqref{trav.wave.5} and satisfies $ 0<f_1(y)\leq T_M $ for $ y\geq 0 $. In addition to that, $ f_1\in C^{2,1/2}(\R_+) $ has bounded first and second derivative according to
				\[|f_1'(y)|\leq \frac{T_M^4}{c} \text{ and }|f_2''(y)|\leq T_M^4\]
				and H\"older seminorm bounded by
				\[[f_1'']_{1/2}\leq \max\{2T_M^4,2T_M^4c+\frac{4T_M^7}{c}\}.\]
				We now show the existence of the solutions $ f_n\in  C^{2,1/2}(\R_+)  $ of the equation \eqref{trav.wave.4} for $ n\geq 2 $. We do the proof only for $ n=2 $, since the very same arguments will work recursively for all $ n\geq 2 $. Let us define $ g=-\int_0^\infty E(y-\eta)f_1^4(\eta)d\eta $. We readily see that $ -T_M^4<g<0 $. Moreover, since $ f_1^4\in C^1(\R_+) $ with bounded derivative we conclude that $ g\in C^{0,1/2}(\R_+) $ with the seminorm $ [\cdot]_{1/2} $ bounded by
				\[[g]_{1/2}\leq \max\{2\Arrowvert g\Arrowvert_\infty,4\Arrowvert f_1\Arrowvert_\infty^3\Arrowvert f_1'\Arrowvert_\infty+\Arrowvert f_1\Arrowvert_\infty^4 \Arrowvert E\Arrowvert_{L^2}\}\leq \max\left\{2T_M^4,\frac{4T_M^7}{c}+T_M^4\Arrowvert E\Arrowvert_{L^2}\right\}.\]
				Indeed, the normalized exponential integral has the property that $ E\in L^q(\R) $for any $ q\in [1,2] $, since $ E\in L^1(\R)\cap L^2(\R) $. This yields together with the H\"older 's inequality that for $ b>a>0 $ and $ \delta\in[0,1/2] $
				\[\int_a^b E(\eta)d\eta \leq|a-b|^{\delta} \Arrowvert E\Arrowvert_{L^{\frac{1}{1-\delta}}}.\] Therefore, for $ v\in C^{0,\delta}(\R_+) $ and $ y>x>0 $ we estimate
				\begin{multline}\label{hoeld.of.I}
					\left|\int_0^\infty v^4(\eta)\left(E((y-\eta))-E((x-\eta))\right)d\eta\right|=\left|\int_{-y}^\infty v^4(\eta+y)E(\eta)d\eta-\int_{-x}^\infty v^4(\eta+x)E(\eta)d\eta\right|\\\leq\left| \int_{-x}^\infty E(\eta)\left(v^4(\eta+y)-v^4(\eta+x)\right)d\eta\right|+\left|\int_{-y}^{-x}E(\eta) v^4(\eta+y)d\eta\right|\\
					\leq[v^4]_\delta|x-y|^\delta+\Arrowvert v^4\Arrowvert_\infty\Arrowvert E\Arrowvert_{L^{\frac{1}{1-\delta}}}|x-y|^{\delta}.
				\end{multline} We remark that if $ v\in C^1(\R_+) $ with bounded derivative and if $ |x-y|<1 $, one can estimate
			\[	\left|\int_0^\infty v^4(\eta)\left(E((y-\eta))-E((x-\eta))\right)d\eta\right|\leq \Arrowvert (v^4)'\Arrowvert_\infty |x-y|^\delta+\Arrowvert v^4\Arrowvert_\infty\Arrowvert E\Arrowvert_{L^{\frac{1}{1-\delta}}}|x-y|^{\delta}\] 
			since also $ |(y+\eta)-(x+\eta)|<1 $.
			
			 Similarly as for the function $ f_1 $, we will consider a suitable minimization problem for which the unique minimizer will be $ f_2 $. Let us consider the minimization problem associated to the functional 
				\[I_g[f]=\int_0^\infty e^{-cy}\left(\frac{(\partial_y f(y))^2}{2}+\frac{f(y)^5}{5}+gf\right)dy\]
				on the set $ \mathcal{A}_{T_M,c} $ defined in \eqref{def.A}.
				Another application of Proposition \ref{prop.variational} shows that there exists $ f_2\in C^{2,1/2}(\R_+) $ solution to \eqref{trav.wave.4} for $ n=2 $ with
				\[|f_2'(y)|\leq \frac{T_M^4}{c},\;|f_2''(y)|\leq T_M^4 \text{ and }[f_2'']_{1/2}\leq \max\left\{2T_M^4,2T_M^4c+\frac{4T_M^7}{c}+\max\left\{2T_M^4,\frac{4T_M^7}{c}+T_M^4\Arrowvert E\Arrowvert_{L^2}\right\}\right\}.\] Moreover, $ 0<f_2(y)\leq T_M $ for $ y\geq 0 $. A recursive application of Proposition \ref{prop.variational} shows the existence of a sequence $ \left(f_n\right)_{n\in\mathbb{N}} \in   C^{2,1/2}(\R_+) $ with $ f_n>0 $ in the interior of $ \R_+ $ solving the recursive system \eqref{trav.wave.4}. Moreover, for all $ n\geq 1$ we have the uniform bounds
				\[f_n(y)\leq T_M,\;|f_n'(y)|\leq \frac{T_M^4}{c},|f_n''(y)|\leq T_M^4 \]and \[[f_n'']_{1/2}\leq \max\left\{2T_M^4,2T_M^4c+\frac{4T_M^7}{c}+\max\left\{2T_M^4,\frac{4T_M^7}{c}+T_M^4\Arrowvert E\Arrowvert_{L^2}\right\}\right\},\] 
				where the uniform bound of the H\"older seminorm is a consequence of the uniform bounds of $ f_{n-1},\;f'_{n-1}$ and $ \;f''_{n-1} $.
				We will now prove that the solutions form a monotonous sequence such that $ 0\leq f_1\leq f_2\leq ...\leq f_n\leq f_{n+1}\leq T_M $. We only need to show that $ f_n\leq f_{n+1} $ for all $ n\in\mathbb{N} $. We prove it by induction. Let us consider $ n=1 $. Then we define $ \varphi=f_2-f_1 $ and \[ a_1(y)=f_1^3(y)+f_2^3(y)+f_1^2(y)f_2(y)+f_1(y)f_2^2(y)> 0.\]
				The strict positivity is due to the fact that by construction $ f_n>0 $ in any open set of $ \R_+ $ and in $ y=0 $. Let $  R>0  $. Then $ \varphi(0)=0 $ as well as $ |\varphi(R)|\leq T_M $. Moreover,
				\[\varphi''-c\varphi'-a_1(y)\varphi(y)\leq 0.\]
				Let us consider now $ \psi_R(y)=-T_M e^{c(y-R)} $. Then we have on one hand that $ \varphi(0)-\psi_R(0)>0 $ as well as $ \varphi(R)-\psi_R(R)\geq 0 $ and on the other hand that \[
				\psi_R''-c\psi_R'-a_1(y)\psi_R=-a_1(y)\psi_R\geq 0.\]
				Hence, an application of the maximum principle to the function $ \varphi-\psi_R $ shows that there is no negative minimum on $ [0,R] $ since
				\[(\varphi-\psi_R)''-c(\varphi-\psi_R)'-a_1(\varphi-\psi_R)\leq 0.\]
				Therefore, $ f_2(y)-f_1(y)\geq -T_Me^{c(y-R)} $ for all $ y\leq R $. Hence, for $ R\to\infty $ we conclude $ f_2\geq f_1 $. 
				
				Let us assume now that for $ n\in\mathbb{N} $ it is true that $ f_{n-1}\leq f_{n} $. We shall now show that $ f_{n}\leq f_{n+1} $. We define $ \varphi_{n}=f_{n+1}-f_n $ and $ a_n(y)=f_n^3(y)+f_{n+1}^3(y)+f_n^2(y)f_{n+1}(y)+f_n(y)f_{n+1}^2(y)> 0 $. Moreover, since by induction $ 0<f_{n-1}\leq f_{n} $ we also have that \[\int_0^\infty E(y-\eta)\left(f_n^4(\eta)-f_{n-1}^4(\eta)\right)d\eta\geq 0. \]
				Hence, we have once more that $ \varphi_n(0)-\psi_R(0)>0 $ and $ \varphi_n(R)-\psi_R(R)\geq 0 $ as well as
				\[(\varphi_n-\psi_R)''-c(\varphi_n-\psi_R)'-a_n(\varphi_n-\psi_R)\leq 0\]
				on $ [0,R] $. We can conclude with the maximum principle that $ f_n-f_{n+1}\geq -T_Me^{c(y-R)} $ for all $ y\leq R $. This yields the claim $ f_n\geq f_{n+1} $.
				
				This concludes the proof of the existence $ f_n\in C^{2,1/2}(\R_+) $ with uniformly bounded $ C^{2,1/2} $-norm solving the recursive system \eqref{trav.wave.4} and satisfying
				\[0\leq f_1\leq f_2\leq ...\leq f_n\leq f_{n+1}\leq T_M.\]
				
				We now prove the existence of a solution to \eqref{trav.wave.3.1}. Let $ f(y)=\lim\limits_{n\to\infty}f_n(y) $. This function exists, since the sequence is monotone and bounded. Moreover, on any compact set $ [0,R] $ the sequence converges also uniformly in $ C^{2,1/4}([0,R]) $ to the function $ f $. Hence, Lebesgue dominated convergence theorem assures that \[\int_0^\infty E(y-\eta)f_n^4(\eta)d\eta\to\int_0^\infty E(y-\eta)f^4(\eta)d\eta\;\;\text{ as }n\to\infty\]
				and the $ C^2 $-uniform convergence in compact sets implies that $ f\in C^2(\R_+)\cap C^{1,1}(\R_+)\cap C^{2,1/2}_{\text{loc}}(\R_+) $ solves \eqref{trav.wave.3.1}, where the $ C^{2,1/2}- $regularity is once again a consequence of elliptic regularity theory. Finally, we prove that $ f\in C^{2,1/2}(\R_+) $ globally. Indeed, $ f\in C^{1,1}(\R_+) $ solves strongly \eqref{trav.wave.3.1}. Thus, $$ f''=cf'+f^4-\int_0^\infty E(\cdot-\eta)f^4 (\eta)d\eta\in C^{0,1/2}(\R_+), $$
				where we used that the convolution of a H\"older continuous function with the exponential integral $ E $ is H\"older continuous as we have proven in \eqref{hoeld.of.I}. This concludes the proof of the existence of traveling wave for \eqref{trav.wave.1} if $ y>0 $. Moreover, the monotonicity of the sequence $ f_n $ implies also $ f(y)>0 $ for any $ y>0 $.
			\end{proof}
		\end{theorem}
		In order to finish the proof of Theorem \eqref{thm.trav.wave.all} we have to show the existence of such $ c_{\max}>0 $. This will be done in the following Lemma and Corollary.
		\begin{lemma}\label{lem.neg.der.}
			Let $ T_M>0 $ and $ c>0 $. Let $ f\in C^2(\R_+)\cap C_b(\R_+) $ be a solution to \eqref{trav.wave.3.1} with $ |f|\leq T_M $. Then $ f(y) >T_M$ for all $ y>0 $ and $ \partial_y f(0^+)<0 $.
			\begin{proof}
				The proof is an adaptation of the proof of Hopf-Lemma. First of all, we notice that by the maximum principle $ f(y)<T_M $ for any $ y\in(0,R) $ with $ R>0 $. Indeed, by assumption we have $ \max\limits_{[0,R]}f=T_M $. If we assume that there exists $ y_0\in(0,R) $ such that $ f(y_0)=T_M $, we obtain the following contradiction
				\[0=f''(y_0)-cf'(y_0)-T_M^4+\int_0^\infty E(\eta-y_0)f^4(\eta)d\eta\leq T_M^4\left(-1+\int_{-R}^\infty E(\eta)d\eta\right)<0.\]
				Thus, since $ f^4(0)-\int_0^\infty E(\eta)f^4(\eta)\geq \frac{T_M^4}{2}>0 $ by continuity there exists $ \delta>0 $ such that $ f(\delta)<T_M $ and $ f^4(y)-\int_0^\infty E(\eta-y)f^4(\eta)d\eta>0 $ for all $ y\in(0,\delta) $.
				
				Let us now consider the operator $ \mathcal{L}=\partial_y^2-c\partial_y $. By construction we see $ \mathcal{L}(f)(y)>0 $ for all $ y\in(0,\delta) $. For $ \alpha>c $ and $ 0<\eps<\frac{T_M-f(\delta)}{e^{\alpha\delta}-1} $ we define the auxiliary function $ z(y)=e^{\alpha y}-1 $. Then a simple computation shows 
				\[\mathcal{L}(f+\eps z)(y)>0 \text{ for all }y\in(0,\delta) \text{ as well as }f(0)+\eps z(0)=T_M>f(\delta)+\eps (\delta) .\]
				Hence, the maximum principle for $ \mathcal{L} $ implies that $ f(y)+\eps z(y)\leq T_M $ for all $ y\in(0,\delta) $. This yields that
				\[f'(0^+)+\eps z'(0^+)=f'(0^+)+\eps\alpha\leq 0\]
				and therefore since $ \alpha>0 $ we conclude $ \partial_yf(0^+)<0 $.
			\end{proof}
		\end{lemma}
		A direct consequence of Lemma \ref{lem.neg.der.} is the following Corollary.
		\begin{corollary}
			There exists $ c_{\max}>0 $ such that for any $ c\in(0,c_{\max}) $ the solution $ f^c $ of \eqref{trav.wave.3} constructed as in Theorem \ref{thm.existence.all} satisfies $ \partial_y f^c(0^+)<-Lc $.
			\begin{proof}
				Let $ c>0 $ and let $ f^c\in C^{2,1/2}(\R_+)$ be the solution to \eqref{trav.wave.3} given by $ f=\alpha^{2/3}\tilde{f}^{\tilde{c}}(\alpha y) $, where $ \tilde{f}^{\tilde{c}} $ is the solution of \eqref{trav.wave.3.1} of Theorem \ref{thm.existence.all} for $ \tilde{c}=\frac{c}{\alpha} $ and melting temperature $ \tilde{T}_M=\frac{T_M}{\alpha^{2/3}} $. Using the bound of the first derivative obtained in Theorem \ref{thm.existence.all} and the definition of the rescaling, we conclude
				\[\Arrowvert \partial_y f^c\Arrowvert_\infty=\alpha^{5/3}\Arrowvert \partial_\eta\tilde{f}^{\tilde{c}}\Arrowvert_\infty\leq\alpha^{5/3}\frac{\tilde{T}_M}{\tilde{c}} =\frac{T_M^4}{c}.\]
				Lemma \ref{lem.neg.der.} implies $ \partial_yf^c(0^+)<0 $. Thus, the set $ \{c>0:\partial_yf^{c}(0^+)<-Lc\} $ is not empty. We hence define
				\[c_{\max}:=\sup\{c>0:\partial_yf^{c}(0^+)<-Lc\}.\]
			\end{proof}
		\end{corollary}
		In the next section we will prove that in the solid the traveling waves are bounded from below by a positive constant and they converge to a positive constant as $ y\to\infty $.
		\subsection{Monotonicity with respect to the melting temperature of the traveling wave solutions for $ y>0 $}\label{subs.monoton}
		In this section we will show that for $ y>0 $ to the traveling waves constructed in the previous section are monotone increasing with respect to the melting temperature, i.e. if $ f_1(0)=\theta_1 $ and $ f_2(0)=\theta_2 $ with $ \theta_1<\theta_2 $ and $ f_1,f_2 $ solve \eqref{trav.wave.3}, then $ f_1\leq f_2 $. We prove the following Lemma
		\begin{lemma}\label{lemma.monotonicity}
			Let $ 0<\theta_1<\theta_2 $ and let $ f_1,f_2\in C^{2,1/2}(\R_+)  $ be the two solutions of \eqref{trav.wave.3} constructed with the iterative scheme in Theorem \ref{thm.existence.all} for $ T_M=\theta_1 $ and $ T_M=\theta_2 $, respectively. Then $ f_1\leq f_2 $.
			\begin{proof}
				Let $ f_1,f_2 $ be given by the limit of the monotone bounded sequences $ f_i^n\in  C^{2,1/2}(\R_+)  $ solving the recursive problem
				\begin{equation*}\label{trav.wave.4.i}
					\begin{cases}
						\partial_y^2 f_i^{n+1}(y)-c\partial_y f_i^{n+1}(y)-\left(f_i^{n+1}(y)\right)^4=g_i^n(y)&y>0;\quad n\geq 1\\
						f_i^0=0&n=0\\
						f_i^{n+1}(0)=\theta_i\\
						f_i^{n+1}\geq 0\\
					\end{cases}
				\end{equation*} 
				where \[g_i^n(y)=-\int_0^\infty E(y-\eta)\left(f_i^n(\eta)\right)^4d\eta.\]
				We show by induction that $ f_1^n\leq f_2^n $ for all $ n\in\mathbb{N} $. This will imply the lemma, since $ f_i(y):=\lim\limits_{n\to\infty}f_i^n(y) $.\\
				
				Let us define $ \varphi_n= f_2^n-f_1^n $. Then $ \varphi_0=0 $ and for $ n\geq 1 $ it solves 
				\begin{equation*}\label{monotonicity.n}
					\begin{cases}
						\partial_y^2 \varphi_n(y)-c\partial_y \varphi_n(y)-a_n(y)\varphi_n(y)=h_{n-1}(y)&y>0;\quad n\geq 1\\
						\varphi_n(0)=\theta_2-\theta_1>0\\
						\varphi_n\in[-\theta_1,\theta_2]\\
					\end{cases}
				\end{equation*}
				where $ h_{n-1}(y)=g_2^{n-1}(y)-g_1^{n-1}(y)= \int_0^\infty E(y-\eta)\left(f_1^{n-1}(\eta)^4-f_2^{n-1}(\eta)^4\right)d\eta$ and \[a_n(y)=\frac{f_2^n(y)^4-f_1^n(y)^4}{f_2^n(y)-f_1^n(y)}=f_2^n(y)^3+f_1^n(y)^3+f_2^n(y)^2f_1^n(y)+f_2^n(y)f_1^n(y)^2>0. \]	
				The positivity of $ a_n(y) $ is given by the strict positivity of $ f_i^n $ in the interior of $ \R_+ $ and in $ y=0 $ as shown before. Moreover, $ \varphi_n\in[-\theta_1,\theta_2] $ since $ 0\leq f_i^{n}\leq \theta_i $ for $ i\in{1,2} $ by the construction in Theorem \ref{thm.existence.all}. We show inductively that $ \varphi_n\geq 0 $ for all $ n\geq 1 $. To this end we consider for $ R>0 $ the function $ \psi_R=-\theta_1e^{c(y-R)} $. It satisfies $ \psi_R(0)\geq -\theta_1 $ as well as $ \psi_R(R)=-\theta_1 $. Hence, on $ [0,R] $ we have
				\begin{equation*}\label{monotonicity.psi.R.n}
					\begin{cases}
						\partial_y^2 \left(\varphi_n(y)-\psi_R(y)\right)-c\partial_y \left(\varphi_n(y)-\psi_R(y)\right)-a_n(y)\left(\varphi_n(y)-\psi_R(y)\right)\\\;\;\;\;\;\;\;=h_{n-1}(y)+a_n(y)\psi_R(y)\leq h_{n-1}(y)&y\in[0,R];\quad n\geq 1\\
						\varphi_n(0)-\psi_R(0)>0\\
						\varphi_n(R)-\psi_R(R)\geq 0\\
					\end{cases}
				\end{equation*} 
				Let us now consider $ n=1 $. Since $ h_0=0 $ the supersolution $ \varphi_1-\psi_R $ solves
				\[\partial_y^2 \left(\varphi_1(y)-\psi_R(y)\right)-c\partial_y \left(\varphi_1(y)-\psi_R(y)\right)-a_1(y)\left(\varphi_1(y)-\psi_R(y)\right)\leq 0.\]
				An application of the maximum principle assuming the existence of a negative minimum,
				gives $ \varphi_1=f_2^1-f_1^1\geq -\theta_1 e^{c(y-R)} $ for $ y\in[0,R] $. Thus, letting $ R\to\infty $ we conclude $ f_2^1\geq f_1^1 $. \\
				
				Let us now assume that for $ n\in \mathbb{N} $ we know that $ f_2^{n-1}\geq f_1^{n-1} $. We show that $ f_2^n\geq f_1^n $. First of all we see that by the induction step we have $ h_{n-1}\leq 0 $, since $  \left(f_1^{n-1}\right)^4\leq \left(f_2^{n-1}\right)^4 $. Then the maximum principle applied to the supersolution $ \varphi_n-\psi_R $ solving
				\[\partial_y^2 \left(\varphi_n(y)-\psi_R(y)\right)-c\partial_y \left(\varphi_n(y)-\psi_R(y)\right)-a_n(y)\left(\varphi_n(y)-\psi_R(y)\right)	\leq 0\]
				implies as before $ f_2^n\leq f_1^n $. This concludes the proof of the lemma.
		\end{proof}\end{lemma}
		
		In the following we aim to show that for $ y>0 $ the constructed traveling wave solutions are bounded from below by a positive constant. This can be proved using the monotonicity property of the traveling wave solutions with respect to the melting temperature. We will indeed show that for very small melting temperature the traveling wave solutions are unique, strictly positive and with a positive limit. 
		\subsection{Traveling wave solutions for small melting temperatures for $ y>0 $}\label{subs.small}
		In this section we will show that for any $ T_M=\eps<\eps_0 $ with $ \eps_0>0 $ small enough there exists a unique solution $ f $ to \eqref{trav.wave.3} which converges to a positive constant with exponential rate $ y\to\infty $. Moreover, $ f $ is bounded from below by a positive constant. We will show it in several steps. We will first prove that any solution $ f $ obtained in Theorem \ref{thm.existence.all} for $ T_M=\eps $ small enough has a limit $ f_\infty $ as $ y\to\infty $ and converges to $ f_\infty $ with exponential rate. Afterwards, we will prove that both $ f $ and $ f_\infty $ are positive and bounded from below by a positive constant. Finally, we will prove that for $ T_M=\eps $ small enough there exists a unique solution to \eqref{trav.wave.3} converging with exponential rate to a constant.
		\begin{lemma}\label{lemma.limit.small.melting}
			Let $ f $ be a solution to \eqref{trav.wave.3.1} as in Theorem \ref{thm.existence.all}. Then for $ T_M=\eps>0 $ small enough there exists $ A>0 $, $ \alpha\in(0,1) $ and $ f_\infty \in[0,T_M] $ such that \[|f(y)-f_\infty|\leq \eps^4Ae^{-\alpha y}.\]
			\begin{proof}
				Let $ f $ be the function obtained in Theorem \ref{thm.existence.all}. First of all we notice that it is equivalent to consider $ f $ solving the equation 
				\begin{equation}\label{trav.wave.3.2}
					\begin{cases}
						\partial_y^2 f(y)-c\partial_y f(y)-\eps^3f^4(y)=-\eps^3\int_0^\infty E(y-\eta)f^4(\eta)d\eta&y>0\\
						f(0)=1\\
						f\geq 0\\
					\end{cases}
				\end{equation}
				Indeed, \eqref{trav.wave.3.2} is obtained considering $ \tilde{f} $ defined by $ \eps \tilde{f}(y)=f(y) $.  Clearly, if $ \tilde{f} $ converges with exponential rate to a constant $ \tilde{f}_\infty $ as $ y\to\infty $, then also $ f $ converges with same rate to $ f_\infty=\eps\tilde{f}_\infty $. Therefore, we will show the lemma for $ \tilde{f} $. In order to simplify the notation we will consider in this proof $ f=\tilde{f} $ solving \eqref{trav.wave.3.2}. 
				
				Since $ f $ is bounded and it solves strongly \eqref{trav.wave.3.2}, then it solves also 
				\[\left(e^{-cy}f'\right)'=\eps^3e^{-cy}\left(f^4-\int_0^\infty E(y-\eta)f^4(\eta)d\eta\right).\] Hence, using that by the boundedness of the first derivative we have $ \lim\limits_{y\to\infty}e^{-cy}f'(y)=0 $, we obtain integrating in $ (y,\infty) $
				\[f'(y)=-\eps^3e^{cy}\int_y^\infty e^{-c\eta}\left(f^4(\eta)-\int_0^\infty E(\eta-z)f^4(z)dz\right)d\eta.\]
				Integrating once more in $ (0,y) $, we conclude that $ f $ solves also the following fixed-point equation
				\begin{equation}\label{fix.eq.small.melting}
					f(y)=1+\eps^3\int_0^y e^{c\xi}\int_\xi^\infty e^{-c\eta}\left(\int_0^\infty E(\eta-z)f^4(z)dz-f^4(\eta)\right)d\eta d\xi.
				\end{equation}
				We define now \begin{equation*}\label{def.osc}
					\osc\limits_{(R,R+1)}f=\sup\limits_{y_1,y_2\in(R,R+1)} |f(y_1)-f(y_2)|.
				\end{equation*}
				Since $ f $ is non-negative and it is bounded by $ 1 $, we know that $ \osc\limits_{(R,R+1)} f\leq 1 $ for all $ R>0 $. For $ M>0 $ we also define
				\begin{equation*}\label{def.lamba}
					\lambda(M)=\sup\limits_{R\geq M} \osc\limits_{(R,R+1)} f.
				\end{equation*} 
				Notice that $ \lambda(M) $ is decreasing with $ \lambda(M)\leq \lambda(0)\leq 1 $. We will show that $ \lambda(M) $ decays like $ e^{-\frac{M}{2}} $. To this end we consider for $ M>0 $ and $ R\geq M $ the points $ y_1,y_2 \in[R,R+1]$ (w.l.o.g. $ y_1\leq y_2 $) and we compute
				\begin{equation*}\label{est.lambda.1}
					\begin{split}
						|f(y_1)-f(y_2)|\leq & \eps^3 \int_{y_1}^{y_2} e^{c\xi}\int_\xi^\infty e^{-c\eta}\left|\int_0^\infty E(\eta-z)f^4(z)dz-f^4(\eta)\right|d\eta d\xi\\
						\leq&  \eps^3 \int_{y_1}^{y_2} e^{c\xi}\int_{y_1}^\infty e^{-c\eta}\left|\int_0^\infty E(\eta-z)f^4(z)dz-f^4(\eta)\right|d\eta d\xi\\
						=&\eps^3 \frac{e^{cy_2}-e^{cy_1}}{c}\int_{y_1}^\infty e^{-c\eta}\left|\int_0^\infty E(\eta-z)f^4(z)dz-f^4(\eta)\right|d\eta,
					\end{split}
				\end{equation*}
				where in the first inequality we used the triangle inequality, in the second we used that $ \xi\geq y_1 $ and the last equality is given by  integrating with respect to $ \xi $.
				We use now that $ 0\leq y_2-y_1\leq 1 $, so that 
				\[\frac{e^{cy_2}-e^{cy_1}}{c}=e^{cy_2}\frac{1-e^{-c(y_2-y_1)}}{c}\leq e^{cy_2}|y_2-y_1|\leq e^{cy_2}\leq \exp(c)e^{cy_1}.\]
				Thus,we can further estimate
				\begin{equation}\label{est.lambda.2}
					\begin{split}
						|f(y_1)-f(y_2)|\leq & \eps^3\exp(c)\int_{y_1}^\infty e^{-c(\eta-y_1)}\left|\int_0^\infty E(\eta-z)f^4(z)dz-f^4(\eta)\right|d\eta\\
						=&\eps^3\exp(c)\int_{y_1}^\infty e^{-c(\eta-y_1)}\left|\int_{-\eta}^\infty E(z)f^4(z+\eta)dz-f^4(\eta)\right|d\eta\\
						=&\eps^3\exp(c)\int_{y_1}^\infty e^{-c(\eta-y_1)}\left|\int_{-\eta}^\infty E(z)\left(f^4(z+\eta)-f^4(\eta)\right)dz-\left(\int_\eta^\infty E(z)dz\right)f^4(\eta)\right|d\eta\\
						\leq & \eps^3\frac{\exp(c)}{2}\int_{y_1}^\infty e^{-c(\eta-y_1)}e^{-\eta}d\eta+\eps^3\exp(c)\int_{y_1}^\infty e^{-c(\eta-y_1)}\left|\int_{-\eta}^\infty E(z)\left(f^4(z+\eta)-f^4(\eta)\right)dz\right|d\eta\\
						\leq & \eps^3\frac{\exp(c)}{2}e^{-M}+\eps^3\exp(c)\int_{y_1}^\infty e^{-c(\eta-y_1)}\left|\int_{-\eta}^\infty E(z)\left(f^4(z+\eta)-f^4(\eta)\right)dz\right|d\eta
					\end{split}
				\end{equation}
				where the first equality follows by a change of coordinates $ z\to z-\eta $ using the symmetry of the kernel $ E $ and the second one is a consequence of the normalization of the kernel $ E $. Moreover, the last inequality uses the boundedness of $ f\leq 1 $ and the estimate 
				\begin{equation}\label{important.est.E}
					\int_a^\infty E(z)dz\leq \frac{e^{-a}}{2} 
				\end{equation} for any $ a>0 $. Finally, we considered $ \eta-y_1\geq 0 $ as well as $ y_1 \geq M $. 
				
				We now estimate the second term in the last line of \eqref{est.lambda.2}. First of all, using that $ |f^4(a)-f^4(b)|\leq 4 |f(a)-f(b)|\leq 4 $ we can rewrite it as the sum of three integrals
				\begin{equation}\label{est.lambda.3}
					\begin{split}
						\eps^3\exp(c)&\int_{y_1}^\infty e^{-c(\eta-y_1)}\left|\int_{-\eta}^\infty E(z)\left(f^4(z+\eta)-f^4(\eta)\right)dz\right|d\eta\\
						\leq & 4 \eps^3\exp(c) \int_{y_1}^\infty e^{-c(\eta-y_1)} \int_{-\eta}^{-M} E(z)dz d\eta\\&+4 \eps^3\exp(c) \int_{y_1}^\infty e^{-c(\eta-y_1)}\int_{0}^\infty E(z)\left|\left(f(z+\eta)-f(\eta)\right)\right|dz d\eta\\&+4 \eps^3\exp(c) \int_{y_1}^\infty e^{-c(\eta-y_1)}\int_{-M}^0 E(z)\left|\left(f(z+\eta)-f(\eta)\right)\right|dz d\eta\\
						\leq & A_1+A_2+A_3
					\end{split}
				\end{equation}
				The first integral term can be estimated easily by
				\begin{equation}\label{est.lambda.4}
					A_1\leq 2\eps^3\frac{\exp(c)}{c}e^{-M},
				\end{equation}
				where we used \eqref{important.est.E} and we solved $ \int_{y_1}^\infty e^{-c(\eta-y_1)}d\eta=\frac{1}{c} $. 
				For the terms $ A_2 $ and $ A_3 $ we will argue in a different way. We recall that $ \lambda(M) $ is decreasing. Hence, if $ z\in (0,1) $ for $ \eta\geq y_1\geq  M $ we have $ |f(\eta)-f(\eta+z)|\leq \lambda(\eta)\leq \lambda (M) $ as well as $ |f(\eta)-f(\eta-z)|\leq \lambda(\eta-1)\leq \lambda (M-1)  $. Thus, using a telescopic sum for $ \eta\geq y_1\geq M $ we compute
				\begin{equation}\label{est.lambda.5}
					\begin{split}
						\int_{0}^\infty E(z)\left|f(z+\eta)-f(\eta)\right|dz
						=&\sum_{n=0}^{\infty} \int_n^{n+1}E(z)\left|f(z+\eta)-f(\eta)\right|dz\\
						\leq&\sum_{n=0}^{\infty} \int_n^{n+1}E(z)\left(|f(\eta+z)-f(n+\eta)|+\sum_{k=1}^n\left|f(\eta+k)-f(\eta+k-1)\right|\right)dz\\
						\leq &\lambda(M)\sum_{n=0}^{\infty} \int_n^{n+1}E(z) (n+1)dz
						\leq  \lambda(M)\sum_{n=0}^{\infty} \int_n^{n+1}E(z) (z+1)dz\\
						=& \lambda(M) \int_0^\infty E(z) (z+1)dz\leq \lambda(M),
					\end{split}
				\end{equation}
				where at the end we used also $ E(a)a\leq \frac{e^{-a}}{2} $ for all $ a>0 $. Thus, \eqref{est.lambda.5} implies
				\begin{equation}\label{est.lambda.6}
					A_2\leq 4\eps^3\frac{\exp(c)}{c}\lambda(M).
				\end{equation}
				Similarly as we did in \eqref{est.lambda.5}, using again a telescopic sum and estimating $ \lambda(0)\leq 1 $, we estimate for $ \eta\geq y_1\geq M $
				\begin{multline}\label{est.lambda.7}
					\int_{-M}^0 E(z)\left|f(z+\eta)-f(\eta)\right|dz=\int_0^M
					E(z)\left|f(\eta)-f(\eta-z)\right|dz	
					=\sum_{n=1}^{M} \int_{n-1}^{n}E(z)\left|f(\eta)-f(\eta-z)\right|dz\\
					\leq \sum_{n=1}^{M} \int_{n-1}^{n}E(z)\left(|f(\eta-z)-f(\eta-(n-1))|+\sum_{k=1}^{n-1}\left|f(\eta-(k-1))-f(\eta-k)\right|\right)dz\\\leq \sum_{n=1}^{M} \int_{n-1}^{n}E(z)\left(\lambda(M-n)+\sum_{k=1}^{n-1} \lambda(M-k)\right)dz
					\leq\sum_{n=1}^{M}\lambda(M-n)\int_{n-1}^{n}E(z)ndz\\
					\leq \sum_{n=1}^{M}\lambda(M-n)\int_{n-1}^{n}E(z)(z+1)dz
					\leq  \sum_{n=1}^{M}\lambda(M-n)\int_{n-1}^{\infty}E(z)(z+1)dz		\end{multline}
					\[\quad\quad\quad\quad\quad\quad\quad\quad\quad\quad\quad\quad\quad\quad\quad\quad\leq  \sum_{n=1}^{M}\lambda(M-n) e^{-(n-1)}\leq e^{-(M-1)}+ e\sum_{n=1}^{M-1}e^{-n}\lambda(M-n).\]
		
				Hence, we have also the following estimate
				\begin{equation}\label{est.lambda.8}
					A_3\leq 4\eps^3\frac{\exp(c+1)}{c}\left[e^{-M}+\sum_{n=1}^{M-1}e^{-n}\lambda(M-n)\right].
				\end{equation}
				Finally, putting together \eqref{est.lambda.2}, \eqref{est.lambda.3}, \eqref{est.lambda.4}, \eqref{est.lambda.6} and \eqref{est.lambda.8} we obtain for $ M\leq R\leq y_1\leq y_2\leq R+1 $
				\begin{equation}\label{est.lambda.9}
					\osc\limits_{(R,R+1)}f\leq |f(y_1)-f(y_2)|\leq \eps^3\exp(c)\left(\frac{1}{2}+\frac{2}{c}+\frac{4e}{c}\right)e^{-M}+4\eps^3\frac{\exp(c)}{c}\lambda(M)+4\eps^3\frac{\exp(c+1)}{c}\sum_{n=1}^{M-1}e^{-n}\lambda(M-n).
				\end{equation}
				Let us take \begin{equation}\label{esp.1}
					\eps<\eps_1(c)=\sqrt[3]{\frac{c}{8\exp(c)}}
				\end{equation} and let us define $ B(c)=\exp(c)\left(1+\frac{4+8e}{c}\right) $. Then taking the supremum over all $ R\geq M $ we have
				\begin{equation}\label{est.lambda.10}
					\lambda(M)\leq B \eps^3 e^{-M}+B\eps^3\sum_{n=1}^{M-1}e^{-n}\lambda(M-n).
				\end{equation}
				We now show by induction that $ \lambda(M)\leq 2B\eps^3 e^{-M/2}$ for all $ \eps<\min\{\eps_1(c),\eps_2(c)\} $, where \begin{equation}\label{eps.2}
					\eps_2(c)=\sqrt[3]{\frac{1}{2B \gamma}}
				\end{equation} for $ \gamma=5\frac{1/2}{1-e^{-1/2}} =\frac{5}{2}\sum_{n=0}^\infty e^{-n/2}$.
				Moreover, since $ \frac{1/2}{1-e^{-1/2}}>\frac{1}{2}$ we have $ \gamma>2 $. This implies also that $ B\eps^3<\frac{1}{2\gamma}<\frac{1}{4} $.
			 First of all we see that if $ M=0 $ the estimates \eqref{est.lambda.7} and \eqref{est.lambda.8} reduce to $ A_3=0 $. Hence, using \eqref{est.lambda.2}, \eqref{est.lambda.3}, \eqref{est.lambda.4}, \eqref{est.lambda.6} we obtain
				\[\lambda(0)\leq \eps^3\left(\frac{\exp(c)}{2}+\frac{2\exp(c)}{c}\right)+4\eps^3\frac{\exp(c)}{c}\lambda(0).\]
				Thus, for $ \eps<\eps_1 $ we have
				\[\lambda(0)\leq B \eps^3\leq 2B \eps^3. \]
				Let us consider $ M=1 $. In this case \eqref{est.lambda.7} and \eqref{est.lambda.8} reduce to $ A_3=\frac{4\eps^3\exp(c)}{c}\lambda(0)e^{-0}\leq \frac{4\eps^3\exp(c+1)}{c}\lambda(0)e^{-1} $, where we used $ \lambda(0)\leq 1 $. Thus, we obtain once more for $ \eps<\min\{\eps_1,\eps_2\} $
				\[\lambda(1)\leq B \eps^3 e^{-1}\leq 2B \eps^3 e^{-1/2}. \]
				Let us now consider $ M=2 $. In this case the sum on the right hand side of \eqref{est.lambda.9} is non-zero. We compute using \eqref{est.lambda.10}
				\[\lambda(2)\leq B \eps^3 e^{-2}+B\eps^3\lambda(1)e^{-1}.\]
				Using now the estimate for $ \lambda(1) $ and that $ \eps<\eps_2 $ and so that $ B\eps^3<1/4 $ we have
				\[\lambda(2)\leq B \eps^3\left(e^{-2}+\frac{e^{-3/2}}{2}\right) \leq 2 B \eps^3 e^{-1}.\]
				Let us now assume that $ \lambda (k) $ satisfies \[\lambda(k)\leq 2B \eps^3 e^{-k/2}\] for $ k=2,..., M\in\mathbb{N} $. We show that also \[\lambda(M+1)\leq 2B \eps^3 e^{-(M+1)/2}.\]
				This is a consequence of the choice of $ \eps_2 $ depending on $ \gamma $. Indeed, by \eqref{est.lambda.10} we have
				\begin{multline*}
					\lambda(M+1)\leq B \eps^3 e^{-(M+1)}+B\eps^3\sum_{n=1}^{M}e^{-n}\lambda(M+1-n)\\
					\leq  B \eps^3 e^{-(M+1)}+B\eps^3e^{-(M+1)/2}\sum_{n=1}^{M}2B \eps^3e^{-n/2}\\
					\leq  B \eps^3 e^{-(M+1)}+B\eps^3e^{-(M+1)/2}\frac{2B\eps^3}{1-e^{-1/2}}
					< 2 B \eps^3 e^{-(M+1)/2},
				\end{multline*}
				where at the end we used the definition of $ \gamma $ as well as $ B\eps^3\gamma<1/2 $ for $ \eps<\eps_2 $.
				This concludes the proof of the exponential decay of $ \lambda(M) $. We will use this result in order to prove the convergence at exponential rate of $ f $.
				Let us consider $ x,y\in\R_+ $ with $ x<y $. Then there exists $ A>0 $ such that $ |f(x)-f(y)|\leq \eps^3A e^{-x/2} $. Indeed we have
				\begin{multline*}\label{est.fx-fy}
					|f(x)-f(y)|\leq \left|f(x)-f(\lfloor x\rfloor)\right|+\left|f(y)-f(\lfloor y\rfloor)\right|+\left|f(\lfloor x\rfloor)-f(\lfloor y\rfloor)\right|\\
					\leq  2B \eps^3 e^{1/2}\left(e^{-x/2}+e^{-y/2}\right)+\sum_{n=\lfloor x\rfloor}^{\lfloor y\rfloor-1} |f(n)-f(n+1)|\\
					\leq 4B \eps^3 e^{1/2}e^{-x/2}+2B \eps^3 \frac{e^{-\lfloor x\rfloor/2}-e^{-\lfloor y\rfloor/2}}{1-e^{-1/2}}\leq \eps^3 A e^{-x/2},
				\end{multline*}
				where $ A=4B e^{1/2}+2B \frac{e^{1/2}}{1-e^{-1/2}} $. Therefore, $ |f(x)-f(y)|\leq \eps^3 A e^{-x/2}\to 0 $ as $ x,y\to\infty $. This implies that for any increasing sequence $ \{x_n\}_{n\in\mathbb{N}}\subset \R_+ $ with $ \lim\limits_{n\to\infty}x_n=\infty $, the sequence $ f(x_n) $ is a Cauchy sequence and hence has a limit as $ n\to\infty $. Indeed, 
				\[|f(x_n)-f(x_m)|\leq \eps^3 A e^{-\frac{\min\{x_n,x_m\}}{2}}\to 0 \;\;\text{ as }n,m\to\infty.\]
				Let hence, $\{x_n\}_{n\in\mathbb{N}}\subset \R_+  $ and $ \{y_n\}_{n\in\mathbb{N}}\subset \R_+  $ be two increasing sequences with $ x_n,y_n\to\infty $ as $ n\to\infty $ and such that \[f_{\infty-}=\liminf\limits_{y\to\infty}f(y)=\lim\limits_{n\to\infty}f(y_n)\leq\lim\limits_{n\to\infty}f(x_n)=\limsup\limits_{y\to\infty}f(y)=f_{\infty+}.\]
				Let $ \delta>0 $. Then there exists some $ N_0\in\mathbb{N} $ such that 
				\[\eps^3 A e^{-\frac{\min\{x_n,y_n\}}{2}}<\frac{\delta}{3}\;\;\text{ for all } n\geq N_0\]
				and \[|f_{\infty+}-f(x_n)|<\frac{\delta}{3}\;\; \text{ as well as }\;\;|f_{\infty-}-f(y_n)|<\frac{\delta}{3}\;\;\text{ for all } n\geq N_0.\]
				Hence, for all $ n\geq N_0 $ we conclude
				\[|f_{\infty-}-f_{\infty+}|\leq |f_{\infty+}-f(x_n)|+|f_{\infty-}-f(y_n)|+|f(x_n)-f(y_n)| <\delta.\]
				This implies that $ f $ has a limit for $ y\to\infty $ which is denoted by \[\liminf\limits_{y\to\infty}f(y)=\limsup\limits_{y\to\infty}f(y)=\lim\limits_{y\to\infty}f(y)=f_{\infty}.\]
				A consequence of the existence of a limit is that any sequence $ \{f(x_n)\}_{n\in\mathbb{N}} $ defined by an increasing diverging sequence $ \{x_n\}_{n\in\mathbb{N}} $ has to converge to $ f_\infty $. Hence, also for $ y\in\R_+ $ we have $ \lim\limits_{n\to\infty}f(y+n)=f_\infty. $

				Finally, let $ y\in\R_+ $. We show that $ f $ converges to $ f_\infty $ with an exponential rate. 
				\begin{equation*}\label{exponential.rate}
					\begin{split}
						|f(y)-f_\infty|= &\sum_{n=0}^\infty \left|f(y+n)-f(y+n+1)\right|
						\leq A\eps^3e^{-\frac{y}{2}}\sum_{n=0}^\infty e^{-\frac{n}{2}}=\frac{A \sqrt{e}\eps^3}{\sqrt{e}-1}e^{-\frac{y}{2}}.
					\end{split}
				\end{equation*}
			\end{proof}
		\end{lemma}
		We continue the theory for small melting temperatures showing that the solution $ f $ of theorem \ref{thm.existence.all} is bounded from below by a positive constant. This will imply that also the limit $ f_\infty $ is strictly positive. We prove the following lemma.
		\begin{lemma}\label{lemma.positivity.small}
			Let $ f $ be a solution to \eqref{trav.wave.3.1} as in Theorem \ref{thm.existence.all}. Then for $ T_M=\eps>0 $ small enough there exists $ c_0>0 $ such that \[f(y)\geq c_0\eps\;\;\text{ for all }y\in\R_+.\]
			This implies also $ f_\infty\geq c_0\eps $.
			\begin{proof}
				As for Lemma \eqref{lemma.limit.small.melting} we consider $ f=\eps \tilde{f} $, where $ \tilde{f} $ solves \eqref{trav.wave.3.2}. We will show that $ \tilde{f}(y)\geq c_0 $ for all $ y\in\R_+ $. This implies clearly the claim of Lemma \ref{lemma.positivity.small}. In order to simplify the notation, we will denote in this proof $ \tilde{f} $ by $ f $. 
				By Lemma \ref{lemma.limit.small.melting} there exist $ f_\infty $ and $ A>0 $ such that $ |f(y)-f_\infty|\leq A\eps^3e^{-y/2} $ for $ \eps>0 $ small enough. As we have seen in Lemma \ref{lemma.limit.small.melting} the solution $ f $ to \eqref{trav.wave.3.2} solves the fixed-point equation \eqref{fix.eq.small.melting}. This can be rewritten as 
				\begin{equation*}\label{fix.eq.small.melting.1}
					f(y)=1+\eps^3\int_0^y e^{c\xi}\int_\xi^\infty e^{-c\eta}\left(\int_0^\infty E(\eta-z)\left[\left(f(z)-f_\infty\right)+f_\infty\right]^4dz-\left[\left(f(\eta)-f_\infty\right)+f_\infty\right]^4\right)d\eta d\xi.
				\end{equation*}
				We recall that \[\left[\left(f-f_\infty\right)+f_\infty\right]^4=(f-f_\infty)^4+4(f-f_\infty)^3 f_\infty+6(f-f_\infty)^2f_\infty^2+4(f-f_\infty)f_\infty^3+f_\infty^4.\]
				Hence, using on the one hand that $ 0\leq f_\infty\leq 1 $, $ |f-f_\infty|\leq 1 $ and that $ |f(y)-f_\infty|\leq A\eps^3e^{-y/2} $ we see easily that 
				\begin{equation}\label{estimate.positivity.1}
					\left[\left(f(y)-f_\infty\right)+f_\infty\right]^4\leq f_\infty^4+15 \eps^3 Ae^{-y/2}.
				\end{equation}
				On the other hand, using in addition that $ (f-f_\infty)^4\geq 0 $ as well as $ (f-f_\infty)^2f_\infty^2\geq 0 $ we have
				\begin{equation}\label{estimate.positivity.2}
					\left[\left(f(y)-f_\infty\right)+f_\infty\right]^4\geq f_\infty^4-8 \eps^3 Ae^{-y/2}.
				\end{equation}
				We can hence estimate from below $ f $ as
				\begin{equation}\label{estimate.positivity.3}
					\begin{split}
						f(y)\geq&1-\eps^3f_\infty^4\int_0^y e^{c\xi}\int_\xi^\infty e^{-c\eta}\int_\eta^\infty E(z)dzd\eta d\xi\\
						&-8\eps^6A\int_0^y e^{c\xi}\int_\xi^\infty e^{-\left(c+\frac{1}{2}\right)\eta}\int_{-\eta}^\infty E(z)e^{-\frac{z}{2}}dzd\eta d\xi-15\eps^6A\int_0^y e^{c\xi}\int_\xi^\infty e^{-\left(c+\frac{1}{2}\right)\eta}d\eta d\xi\\
						\geq& 1-\frac{\eps^3}{2}\int_0^y e^{c\xi}\int_\xi^\infty e^{-\left(c+1\right)\eta}d\eta d\xi-\eps^6A\left(16\artanh\left(\frac{1}{2}\right)+15\right)\int_0^y e^{c\xi}\int_\xi^\infty e^{-\left(c+\frac{1}{2}\right)\eta}d\eta d\xi\\
						=&1-\frac{\eps^3}{2(c+1)}\int_0^y e^{-\xi} d\xi-\frac{\eps^6A}{c+\frac{1}{2}}\left(16\artanh\left(\frac{1}{2}\right)+15\right)\int_0^y e^{-\frac{\xi}{2}}d\xi\\
						\geq& 1-\eps^3\left(\frac{1}{2(c+1)}+\frac{4\eps^3A}{2c+1}\left(16\artanh\left(\frac{1}{2}\right)+15\right)\right).
					\end{split}
				\end{equation}
				We used for the second inequality the fact that $ 0\leq f_\infty\leq 1 $, as well as \eqref{important.est.E}. Moreover, for the equality we used that for any $ a\in[0,1) $ \begin{equation*}\label{arctanh}
					\int_{-\infty}^\infty E(z)e^{-az}dz=\int_{0}^\infty E_1(z)\cosh(az) dz=\frac{\artanh(a)}{a}.
				\end{equation*}
				Equation \eqref{estimate.positivity.3} implies that for any $ c_0\in(0,1) $ defining 
				\begin{equation}\label{esp3}
					\eps_3=\min\left\{1,\sqrt[3]{(1-c_0)\left(\frac{1}{2(c+1)}+\frac{4A}{2c+1}\left(16\artanh\left(\frac{1}{2}\right)+15\right)\right)^{-1}}\right\}
				\end{equation}
				and choosing $ \eps<\min\{\eps_1,\eps_2,\eps_3\} $ according to \eqref{esp.1}, \eqref{eps.2} and \eqref{esp3}, we conclude that the function $ f $ satisfies 
				\[f(y)\geq c_0.\]
				This concludes the proof of the lemma.
			\end{proof}
		\end{lemma}
		Lemma \ref{lemma.monotonicity} and Lemma \ref{lemma.positivity.small} imply the following Corollary.
		\begin{corollary}\label{cor.positivity}
			Let $ T_M>0 $ and let $ f $ be a solution to \eqref{trav.wave.3.1} as in Theorem \ref{thm.existence.all}. Then there exists $ \lambda>0 $ such that $ f(y)\geq \lambda>0 $ for all $ y\geq 0 $.
		\end{corollary}
		Finally, we show that if $ T_M=\eps $ small enough the solution of \eqref{trav.wave.3.1} of Theorem \ref{thm.existence.all} is also unique. Indeed, we show that there is a unique solution of the fixed-point equation \eqref{fix.eq.small.melting} converging to a constant with exponential rate. This is stated in the following theorem.
		\begin{theorem}\label{thm.small.melting}
			Let $ T_M=\eps $. Then, for $ \eps<\eps_0$ small enough there exists a unique solution $f\in C^{2,1/2}(\R_+)    $ of \eqref{trav.wave.3.1} with $ \lim\limits_{y\to\infty}f(y)=f_\infty $ and $ |f(y)-f_\infty|\leq Ae^{-y/2} $. Moreover, $ f(y)\geq c_0\eps $ as well as $ f_\infty\geq c_0\eps $ for $ c_0\in(0,1) $.
			\begin{proof}
				First of all we remark that it is enough to prove the existence and uniqueness of the solution $ \tilde{f} $ to the equation \eqref{trav.wave.3.2}. Indeed, then $ f=\eps\tilde{f} $ is the desired unique solution of Theorem \ref{thm.small.melting}. We will indeed prove the theorem for $ \tilde{f} $, which is denoted in the rest of the proof by the sake of simplicity $ \tilde{f}=f $. \\
				
				Moreover, it is enough also to show the existence and uniqueness of the solution to the fixed-point equation \eqref{fix.eq.small.melting}. Indeed, any strong solution $ f $ to \eqref{trav.wave.3.2} satisfies $ f\in  C^{2,1/2}(\R_+) $ and it solves \eqref{fix.eq.small.melting}.
				
				Let us consider for $ B>1$ and for $ A>0 $ the following space
				\begin{equation*}\label{def.X}
					\mathcal{X}=\left\{f\in C_b(\R_+):|f(y)|\leq B,\; \exists f_\infty \text{ s.t. }|f(y)-f_\infty|\leq Ae^{-y/2} \right\}
				\end{equation*}
				equipped with the metric $ d_\mathcal{X} $ induced by the following norm
				\begin{equation*}\label{def.norm}
					\Arrowvert f\Arrowvert_{\mathcal{X}}=|f_\infty|+\sup\limits_{y\in\R_+}e^{y/2}|f(y)-f_\infty|.
				\end{equation*}
				We also define the following seminorm
				\begin{equation*}\label{def.seminorm}
					[f]_{\mathcal{X}}=\sup\limits_{y\in\R_+}e^{y/2}|f(y)-f_\infty|
				\end{equation*}
				so that $ 	\Arrowvert f\Arrowvert_{\mathcal{X}}=|f_\infty|+[f]_{\mathcal{X}} $. One can prove that $ \left(\mathcal{X},d_\mathcal{X} \right) $ is a complete metric space. We omit the elementary proof.\\
				
				We will now prove that the map 
				\[\mathcal{L}[f](y)=1+\eps^3\int_0^y e^{c\xi}\int_\xi^\infty e^{-c\eta}\left(\int_0^\infty E(\eta-z)f^4(z)dz-f^4(\eta)\right)d\eta d\xi\]
				is a selfmap $ \mathcal{L}:\mathcal{X}\to\mathcal{X} $ and that it is a contraction for $ \eps<\eps_4 $ small enough. The Banach fixed-point theorem will imply the existence of a unique fixed-point $ f $ solving \eqref{trav.wave.3.2}.\\
				
				Let now $ f\in\mathcal{X} $.
				We observe that if $ f\in C_b(\R_+) $ then $ \mathcal{L}[f] $ is continuous.
				We move on proving that for $ f\in\mathcal{X} $ also $ \mathcal{L}[f] $ is bounded. Indeed, using that $ |f|\leq B $ as well as $ |f-f_\infty|\leq 2B $, we obtain similarly as for \eqref{estimate.positivity.1} and for  \eqref{estimate.positivity.2} that
				\begin{equation*}\label{estimate.bnd.1}
					\left[\left(f(y)-f_\infty\right)+f_\infty\right]^4\leq f_\infty^4+40B^3   Ae^{-y/2}
				\end{equation*}
				and
				\begin{equation*}\label{estimate.bnd.2}
					\left[\left(f(y)-f_\infty\right)+f_\infty\right]^4\geq f_\infty^4-20B^3  Ae^{-y/2}.
				\end{equation*}
				Thus, we estimate similarly as in \eqref{estimate.positivity.3}
				\begin{equation}\label{estimate.bnd.3}
					\begin{split}
						\left|\mathcal{L}[f](y)\right|\leq 1+\eps^3B^3\left(\frac{B}{2(c+1)}+\frac{4AB}{2c+1}\left(40\artanh\left(\frac{1}{2}\right)+20\right)\right).
					\end{split}
				\end{equation}
				Hence, defining by $ c_1(A,B)=\left(\frac{B}{2(c+1)}+\frac{4AB}{2c+1}\left(40\artanh\left(\frac{1}{2}\right)+20\right)\right) $ and taking \begin{equation}\label{eps.5}
					\eps_5=\frac{1}{B}\sqrt[3]{\frac{B-1}{c_1(A,B)}} 
				\end{equation}we see that for $ \eps<\eps_5 $ we have 
				\[	\left|\mathcal{L}[f](y)\right|\leq B.\]
				We have now to show that $ \mathcal{L}[f] $ has also a limit as $ y\to\infty $, which we will call $ \mathcal{L}_\infty[f] $. Moreover, we shall show that $ \left|\mathcal{L}[f](y)-\mathcal{L}_\infty[f]\right|\leq Ae^{-y/2} $. This is the consequence of the convergence of $ f $ to $ f_\infty $ with exponential rate. Let us define
				\begin{equation*}\label{limit.L.f}
					\mathcal{L}_\infty[f]=1+\eps^3\int_0^\infty e^{c\xi}\int_\xi^\infty e^{-c\eta}\left(\int_0^\infty E(\eta-z)f^4(z)dz-f^4(\eta)\right)d\eta d\xi.
				\end{equation*}
				By \eqref{estimate.bnd.3} we know that $ \mathcal{L}_\infty[f] $ is bounded. Moreover, using that
				\[\left|f^4(y)-f^4_\infty\right|\leq 4B^3 \left|f(y)-f_\infty\right|\leq 4AB^3 e^{-y/2}\]
				we can estimate
				\begin{equation*}\label{limit.L.f.1}
					\begin{split}
						\left|\mathcal{L}[f](y)\right.&\left.-\mathcal{L}_\infty[f]\right|\leq \eps^3\int_y^\infty e^{c\xi}\int_\xi^\infty e^{-c\eta}\left|\left(\int_0^\infty E(\eta-z)f^4(z)dz-f^4(\eta)\right)\right|d\eta d\xi\\
						= &\eps^3\int_y^\infty e^{c\xi}\int_\xi^\infty e^{-c\eta}\left|\left(\int_0^\infty E(\eta-z)\left(f^4(z)-f_\infty^4\right)dz-\left(f^4(\eta)-f_\infty^4\right)\right)-\int_\eta^\infty E(z)dz f_\infty^4\right|d\eta d\xi\\
						\leq &\eps^3 4AB^3\int_y^\infty e^{c\xi}\int_\xi^\infty e^{-(c+1/2)\eta}\left(\int_{-\infty}^\infty E(z)e^{-z/2}dz+1\right)d\eta d\xi+\eps^3 \frac{B^4}{2}\int_y^\infty e^{c\xi}\int_\xi^\infty e^{-(c+1)\eta}d\eta d\xi\\
						\leq & \eps^3B^3\left[4A\frac{2\artanh\left(\frac{1}{2}\right)+1}{2c+1}\right]e^{-y/2}+\eps^3 \frac{B^4}{2(c+1)}e^{-y}.
					\end{split}
				\end{equation*}
				Hence, defining 
				\begin{equation}\label{eps.6}
					\eps_6=\left[\frac{B^3}{A}\left(4A\frac{2\artanh\left(\frac{1}{2}\right)+1}{2c+1}+\frac{B}{2(c+1)}\right)\right]^{-\frac{1}{3}}
				\end{equation}
				we can conclude that there exists a limit 
				\[\lim\limits_{y\to\infty}\mathcal{L}[f](y)=\mathcal{L}_\infty[f]\]
				such that 
				\[\left|\mathcal{L}[f](y)-\mathcal{L}_\infty[f]\right|\leq Ae^{-y/2}\]
				for all $ \eps<\min\left\{\eps_5,\eps_6\right\} $, defined in \eqref{eps.5} and \eqref{eps.6}. This concludes the proof of $ \mathcal{L} $ being a self-map. We now finish the proof of the theorem showing that $ \mathcal{L} $ is also a contraction map.
				
				We first prove that there exists a constant $ c_2(A,B) $ such that if $ f,g\in\mathcal{X} $, then \[\Arrowvert f^4-g^4\Arrowvert_{\mathcal{X}}\leq c_2(A,B)\Arrowvert f-g\Arrowvert_{\mathcal{X}}.\]
				We recall that if $ g\in\mathcal{X} $, then $ [g]_{\mathcal{X}}\leq A $ and $ |g|\leq B $. Hence, we have the estimate
				\begin{equation*}\label{equivalence.norm}
					\begin{split}
						e^{y/2}\left|\right.&\left.\left(f(y)^4-g(y)^4\right)-\left(f_\infty^4-g_\infty^4\right)\right|=e^{y/2}\left|\left(f(y)-f_\infty+f_\infty\right)^4-\left(g(y)-g_\infty+g_\infty\right)^4-\left(f_\infty^4-g_\infty^4\right)\right|\\
						\leq & e^{y/2}\left|\left(f(y)-f_\infty\right)^4-\left(g(y)-g_\infty\right)^4\right|+e^{y/2}\left|4f_\infty\left(f(y)-f_\infty\right)^3-4g_\infty\left(g(y)-g_\infty\right)^3\right|\\
						&+e^{y/2}\left|6f_\infty^2\left(f(y)-f_\infty\right)^2-6g_\infty^2\left(g(y)-g_\infty\right)^2\right|+e^{y/2}\left|4f_\infty^3\left(f(y)-f_\infty\right)-4g_\infty^3\left(g(y)-g_\infty\right)\right|\\
						\leq & \left(16+24+12\right)AB^2|f_\infty-g_\infty|+\left(32+48+24+4\right)B^3[f-g]_{\mathcal{X}}\\=&52 AB^2|f_\infty-g_\infty|+108B^3[f-g]_{\mathcal{X}}.
					\end{split}
				\end{equation*}
				Thus, using that $ |f_\infty^4-g_\infty^4|\leq 4B^3 |f_\infty-g_\infty| $ and defining $ c_2(A,B)=\max\{52AB^2+4B^3,108B^3\} $ we conclude that
				\[\Arrowvert f^4-g^4\Arrowvert_{\mathcal{X}}\leq c_2(A,B)\Arrowvert f-g\Arrowvert_{\mathcal{X}}.\]
				Moreover, we see that $ [f^4-g^4]_{\mathcal{X}}\leq c_2(A,B)\Arrowvert f-g\Arrowvert_{\mathcal{X}} $, which implies \[ \left|\left(f(y)^4-g(y)^4\right)-\left(f_\infty^4-g_\infty^4\right)\right|\leq c_2(A,B))\Arrowvert f-g\Arrowvert_{\mathcal{X}} e^{-y/2}.\]
				Hence, we estimate
				\begin{equation*}\label{contraction.1}
					\begin{split}
						\left|\mathcal{L}_\infty[f]\right.&\left.-\mathcal{L}_\infty[g]\right|\leq \eps^3\int_0^\infty e^{c\xi}\int_\xi^\infty e^{-c\eta}\left|\left(\int_0^\infty E(\eta-z)\left(f^4(z)-g^4(z)\right)dz-\left(f^4(\eta)-g^4(z)\right)\right)\right|d\eta d\xi\\
						=& \eps^3\int_0^\infty e^{c\xi}\int_\xi^\infty e^{-c\eta}\left|\left(\int_0^\infty E(\eta-z)\left(f^4(z)-f^4_\infty-\left(g^4(z)-g^4_\infty\right)\right)dz-\left(f^4(\eta)-f_\infty^4-\left(g^4(z)-g_\infty^4\right)\right)\right)\right.\\&\phantom{spacespacespacespacespace}\left.-\int_\eta^\infty E(z)dz\left(f_\infty^4-g_\infty^4\right)\right|d\eta d\xi\\
						\leq& \eps^3 [f^4-g^4]_{\mathcal{X}}\left(2\artanh\left(\frac{1}{2}\right)+1\right)\int_0^\infty e^{c\xi}\int_\xi^\infty e^{-{(c+1/2)}\eta}d\eta d\xi+\eps^3 \frac{\left|f_\infty^4-g_\infty^4\right|}{2}\int_0^\infty e^{c\xi}\int_\xi^\infty e^{-{(c+1)}\eta}d\eta d\xi\\
						\leq & \eps^3 c_2(A,B)\left(4\frac{2\artanh\left(\frac{1}{2}\right)+1}{2c+1}+\frac{1}{c+1}\right)\Arrowvert f-g\Arrowvert_{\mathcal{X}}.
					\end{split}
				\end{equation*}
				In a similar way we can estimate
				\begin{equation*}\label{contraction.2}
					\begin{split}
						\left|\mathcal{L}[f](y)-\right.&\left.\mathcal{L}[g](y)-\left(\mathcal{L}_\infty[f]-\mathcal{L}_\infty[g]\right)\right|\\\leq& \eps^3\int_y^\infty e^{c\xi}\int_\xi^\infty e^{-c\eta}\left|\left(\int_0^\infty E(\eta-z)\left(f^4(z)-g^4(z)\right)dz-\left(f^4(\eta)-g^4(z)\right)\right)\right|d\eta d\xi\\
						\leq &	\eps^3 [f^4-g^4]_{\mathcal{X}}\left(2\artanh\left(\frac{1}{2}\right)+1\right)\int_y^\infty e^{c\xi}\int_\xi^\infty e^{-{(c+1/2)}\eta}d\eta d\xi+\eps^3 \frac{\left|f_\infty^4-g_\infty^4\right|}{2}\int_y^\infty e^{c\xi}\int_\xi^\infty e^{-{(c+1)}\eta}d\eta d\xi\\	
						\leq& \eps^3 c_2(A,B)\left(4\frac{2\artanh\left(\frac{1}{2}\right)+1}{2c+1}\right)\Arrowvert f-g\Arrowvert_{\mathcal{X}}e^{-y/2}+\eps^3\frac{ c_2(A,B)}{c+1}|f_\infty-g_\infty|e^{-y}.
					\end{split}
				\end{equation*}
				This estimate implies easily 
				\begin{equation*}\label{contraction.3}
					\begin{split}
						\left[\mathcal{L}[f](y)-\mathcal{L}[g](y)\right]_\mathcal{X}\leq\eps^3 c_2(A,B)\left(4\frac{2\artanh\left(\frac{1}{2}\right)+1}{2c+1}+\frac{1}{c+1}\right)\Arrowvert f-g\Arrowvert_{\mathcal{X}}.
					\end{split}
				\end{equation*}
				Therefore, taking $ \theta\in(0,1) $ and 
				\begin{equation}\label{eps.7}
					\eps_7=\left[\frac{2c_2(A,B)}{\theta}\left(4\frac{2\artanh\left(\frac{1}{2}\right)+1}{2c+1}+\frac{1}{c+1}\right)\right]^{-\frac{1}{3}}
				\end{equation}
				we conclude that the map $ \mathcal{L}:\mathcal{X}\to \mathcal{X} $ is a contraction self-map for $ \eps<\min\{\eps_5,\eps_6,\eps_7\}=\eps_4 $, given in \eqref{eps.5}, \eqref{eps.6} and \eqref{eps.7}. Hence, there exists a unique fixed-point $ \tilde{f} $ of the equation \eqref{fix.eq.small.melting}, which solves also \eqref{trav.wave.3.2}. Finally, $ f=\eps \tilde{f}\in C^{2,1/2}(\R_+)   $ solves \eqref{trav.wave.3.1}. Taking now $ \eps_0=\min\{\eps_1,\eps_2,\eps_3,\eps_4\} $, for $ \eps_i $ defined in \eqref{esp.1}, \eqref{eps.2}, \eqref{esp3} and above, Lemma \ref{lemma.limit.small.melting} and Lemma \ref{lemma.positivity.small} imply Theorem \ref{thm.small.melting}.
			\end{proof}
		\end{theorem}
		\section{Existence of the limit of the traveling wave solutions as $ y\to\infty $}\label{Sec.exist.limit}
		We have proved in Theorem \ref{thm.existence.all} the existence for any $ c>0 $ of a traveling wave $ f $ in $ \R_+ $ solving \eqref{trav.wave.3.1} and with the property that that $  f\in C^{2,1/2}(\R_+)  $. Moreover, as we have seen in Corollary \ref{cor.positivity}, $ f $ is bounded from below by a positive constant as long as $ T_M>0 $. In this section we will show that $ f $ has a limit as $ y\to\infty $.
		
		We will proceed as follows. We will show that for any sequence $ \{y_n\}_{n\in\mathbb{N}} $ increasing and diverging, the sequence $ f_n(y)=f(y+y_n) $ has a subsequence converging to a function, which will be denoted by an abuse of notation as $ \omega $-limit. This definition relies on the similarity with the notion of $ \omega $-limit point for dynamical systems. Analogously, the $ \omega $-limit set is given in this setting by all the existing limit functions $ \lim\limits_{k\to\infty}f(y+y_k) $, i.e \[\omega(f):=\left\{\bar f:\R\to\R \text{ such that }\exists\{y_k\}_k,\;y_k<y_{k+1},\; y_k\to\infty \text{ as }k\to\infty\text{, satisfying }\lim\limits_{k\to\infty}f(y+y_k)=\bar f(y)\right\}.\] We will prove that any $ \omega $-limit is a constant function. This will be used in the end in order to show that $ f $ has a limit.\\
		\subsection{Elementary properties of the $ \omega $-limits of the traveling waves}
		Let us consider $ \{y_n\}_{n\in\mathbb{N}}\subset \R_+ $ any increasing sequence with $ \lim\limits_{n\to\infty} y_n=\infty$ and let us consider $ f_n(y):=f(y+y_n) $. Then $ f_n:[-y_n,\infty)\to \R_+ $ solves for $ \lambda>0 $ small enough
		\begin{equation*}\label{trav.wave.3.n}
			\begin{cases}
				\partial_y^2 f_n(y)-c\partial_y f_n(y)-f_n^4(y)=-\int_{-y_n}^\infty E(y-\eta)f_n^4(\eta)d\eta&\text{ for }y>-y_n\\
				f(-y_n)=T_M\\
				f\geq \lambda>0.\\
			\end{cases}
		\end{equation*}
		Since $ f_n\in C^{2,1/2}[-y_n,\infty) $, by compactness a diagonal argument shows that there exists a subsequence $ f_{n_k} $ such that $ f_{n_k}\to \bar f $ in $ C^{2,\alpha}([-R,R]) $ for $ \alpha\in \left(0,\frac{1}{2}\right) $ and for any $ R>0 $. Therefore $ \bar f\in C^2(\R) $ and by the uniform boundedness of $ f_n' $ and $ f_n'' $ we also have $ \Arrowvert \bar f'\Arrowvert_\infty\leq \frac{2T_M^4}{c} $ and $ \Arrowvert \bar f''\Arrowvert_\infty\leq 4T_M^4 $. Moreover, an application of the dominated convergence theorem yields that $ \bar f $ solves 
		\begin{equation}\label{trav.wave.3.n.R}
			\begin{cases}
				\partial_y^2 \bar f(y)-c\partial_y \bar f(y)-\bar{f}^4(y)=-\int_{-\infty}^\infty E(y-\eta)\bar{f}^4(\eta)d\eta&y\in\R\\
				0<\lambda\leq \bar f\leq T_M.\\
			\end{cases}
		\end{equation}
		Hence, regularity theory implies $ \bar f  \in C^{2,1/2}(\R)$, since the convolution $ E*\bar{f}^4\in C^{0,1/2}(\R) $.
		
		\begin{lemma}\label{lemma.no.sup.interior}
			Let $ f $ solve \eqref{trav.wave.3.n.R}. Then $ f $ does not attain its supremum and infimum at the interior, unless $ f $ is constant. 
			\begin{proof}
				The proof is a direct consequence of the maximum principle. Let us assume that $ f $ is not constant and that there exists $ y_m\in \R $ or $ y_M\in\R $ such that $ \sup\limits_{\R} f=f(y_M)$ or $ \inf\limits_\R f=f(y_m) $. Then by the positivity of $ f $ we see that $ f^4(y)-f^4(y_M)\leq 0 $ as well as $ f^4(y)-f^4(y_m)\geq 0 $. Moreover, $ f $ differs from its maximum and minimum in sets of positive measures, since $ f $ is continuous and non-constant. Hence, we obtain the following contradictions
				\[0= f''(y_M)-cf'(y_M)+\int_\R E(\eta-y)\left[f^4(\eta)-f^4(y_M)\right]d\eta <0\]
				if the supremum is attained at the interior or
				\[0= f''(y_m)-cf'(y_m)+\int_\R E(\eta-y)\left[f^4(\eta)-f^4(y_m)\right]d\eta >0\]
				if the infimum is attained at the interior. This concludes the proof of the lemma.
			\end{proof}
		\end{lemma}
		This result implies that, if $ \bar f $ is not constant, it have to attain its supremum and infimum at $ +\infty $ or $ -\infty $. We will prove that this is not possible. We start showing that $ \bar{f} $ does not attain its supremum and infimum at $ +\infty$.
		\begin{lemma}\label{lemma.no.sup.infty}
			Let $ f $ solve \eqref{trav.wave.3.n.R}. Then $ f $ does not attain its supremum at $ +\infty $, i.e. $ \limsup\limits_{y\to\infty}f(y)<\sup\limits_{\R}f $, unless $ f $ is constant. 
			\begin{proof}
				The proof is based once again on the maximum principle. Let us assume that $ f $ is not constant and that $  \limsup\limits_{y\to\infty}f(y)=\sup\limits_{\R}f=:A  $. We consider the function $ \omega=A-f\geq 0 $. Moreover, since $ f $ is not constant also $ \omega>0 $ at the interior by Lemma \ref{lemma.no.sup.interior}. Hence, $ \omega $ solves
				\begin{equation}\label{eq.omega}
					-\omega''(y)+c\omega'(y)-(A-\omega(y))^4+\int_\R E(y-\eta)(A-\omega(\eta))^4d\eta=0.
				\end{equation}
				We will show that $ \omega(y)>0 $ as $ y\to\infty $, which is a contradiction with the assumption of $ f $ attaining its supremum at $ +\infty $. To this end we construct a suitable family of subsolutions $ \psi_\delta(y) $ with the property $ f\geq \psi_\delta $ and such that $ \psi_\delta>0 $ for $ y\in[0,R_\delta) $ for a suitable $ R_\delta\to\infty $ as $ \delta\to 0 $. \\
				
				We define the following constants. First of all we take $ \theta=\frac{1}{5} $ and $ R >0$ fixed so that
				\begin{equation}\label{def.R}
					\int_y^{R+y}E(\eta)d\eta>\int_{R+y}^\infty E(\eta)d\eta \;\;\;\;\;\text{ for all }y>0.
				\end{equation}
				Moreover, we define $ c_0=\min\{1,c\} $ and we take $ \beta\in\left(0,\frac{c_0}{4}\right) $ fixed so that \\
				\begin{equation}\label{def.beta}
					\frac{\artanh(4\beta)}{4\beta} <\frac{3}{2}\;\;\;\;\;\text{ and }\;\;\;\;\;\beta^2-\frac{c_0}{4}\beta+4A^3\left(\frac{\artanh(\beta)}{\beta}-1\right)\leq 0.
				\end{equation}
				For a suitable constant $ C(\beta,A,\theta)>0 $, which will be computed later, we also fix 
				\begin{equation}\label{def.eps}
					\eps<\min\left\{\min\limits_{\left[-R,\frac{\ln(2)}{\beta}\right]}\{A-f(y)\},\; C(\beta,A,\theta)\right\}.
				\end{equation}
				Finally, for $ \delta_0=\frac{\eps\theta}{2} $ we consider the following family of subsolutions
				\begin{equation}\label{def.psi.delta}
					\psi_\delta(y)=\begin{cases}
						0&\quad y<-R\\
						\eps-\delta e^{\beta y}& \quad y\in[-R,0)\\
						\eps\theta-\delta e^{\beta y}& \quad y\in[0,R_\delta]\\
						0& \quad y>R_\delta,
					\end{cases}
				\end{equation}
				where $ R_\delta=\frac{1}{\beta}\ln\left(\frac{\eps\theta}{\delta}\right)\to \infty $ as $ \delta\to 0 $ as well as $ \eps\theta-\delta e^{\beta R_\delta}=0 $.
				By construction, $ \psi_\delta \leq \omega $ for $ y\in \R\setminus(0,R_\delta) $. We will show that on $ (0,R_\delta) $ the family $ \psi_\delta $ consists of subsolutions to \eqref{eq.omega}.
				However, before  moving to the proof of this claim we show that equations \eqref{def.R} and \eqref{def.beta} are well-defined. We first show the function \[h(y)=\int_y^{R+y}E(\eta)d\eta-\int_{R+y}^\infty E(\eta)d\eta\]
				is a decreasing function. 
				Using the definition of the kernel $ E $, we notice
				\[h(0)=\frac{1}{2}-e^{-R}+2RE(R)>0\quad\quad\quad \text{ for }R>0\text{ large enough}.\]
				Moreover, $ \lim\limits_{y\to\infty}h(y)= 0$. We compute also for $ R>\max\{1,\ln(2)\}=1 $ 
				\[h'(y)=2E(R+y)-E(y) \quad \text{ and }\quad h''(y)=\frac{e^{-y}}{2y}-\frac{e^{-(y+R)}}{y+R}>\frac{e^{-y}}{y}\left(\frac{1}{2}-e^{-R}\right)>0.\]
				Since $ \lim\limits_{y\to 0}h'(y)=-\infty  $ and $  \lim\limits_{y\to \infty}h'(y)=0 $, we conclude $ h'(y)<0 $. This implies that $ h $ is monotonically decreasing for $ R>1 $. Therefore, there exists an $ R>0 $ such that \eqref{def.R} holds.\\
				
				We move to the existence of $ \beta\in\left(0,\frac{c_0}{4}\right) $ solving \eqref{def.beta}. First of all, let us define  $ g(y)=\frac{\artanh(y)}{y} $. Then $ g:(0,1)\to\R_+ $. Hence, $ \beta\in\left(0,\frac{c_0}{4}\right)  $ is well-defined. Moreover, elementary calculus implies 
				\[\lim\limits_{y\to 0}g(y)=1,\quad \lim\limits_{y\to 1}g(y)=\infty,\quad g'(y)\geq 0\text{ with }\quad g'(0)=0, \quad \text{ and }\quad g''(y)\geq 0.\]
				
				Therefore,
				$ g $ is a convex monotone non-decreasing function with $ g(0)=1 $. Hence, there exists $ \beta_0\in\left(0,\frac{c_0}{4}\right) $ such that $ g(4\beta)<\frac{3}{2} $ for all $ \beta<\beta_0 $. Moreover, the function
				$ k(\beta)=\beta^2-\frac{c_0}{4}\beta+4A^3(g(\beta)-1) $ is convex as sum of convex functions. Since $ k(0)=0 $, $ k\left(\frac{c_0}{4}\right)>0 $ as well as \[k'(\beta)=-\frac{c_0}{4}+[2\beta+4A^3g'(\beta)]\underset{\beta\to0}{\longrightarrow} -\frac{c_0}{4}<0\]
				we conclude the existence of a $ \beta $ satisfying \eqref{def.beta}.\\
				
				We prove now that $ \psi_\delta $ are subsolutions to \eqref{eq.omega} for $ y\in(0,R_\delta) $, where the functions are smooth.
				We compute for $ y\in(0,R_\delta) $
				\begin{equation}\label{subsol.1}
					\begin{split}
						-\psi_\delta''(y)&+c\psi_\delta'(y)-\left(A-\psi_\delta(y)\right)^4+\int_\R E(y-\eta)\left(A-\psi_\delta(\eta)\right)^4d\eta\\=&-\left(c-\frac{c_0}{4}\right)\beta\delta e^{\beta y}
						+\delta e^{\beta y}\left(\beta^2-\frac{c_0}{4}\beta\right)-\left(A-\eps \theta+\delta e^{\beta y}\right)^4+\int_{-\infty}^{-R}E(y-\eta)A^4 d\eta\\&+\int_{-R}^0E(y-\eta)\left(A-\eps+\delta e^{\beta \eta}\right)^4d\eta+\int_{0}^{R_\delta}E(y-\eta)\left(A-\eps\theta+\delta e^{\beta \eta}\right)^4d\eta\\
						&+\int_{R_\delta}^\infty E(y-\eta)\left(A-\eps\theta+\delta e^{\beta R_\delta}\right)^4d\eta\\
						\leq&-\frac{3c_0}{4}\beta\delta e^{\beta y}
						+\delta e^{\beta y}\left(\beta^2-\frac{c_0}{4}\beta\right)-\left(A-\eps \theta+\delta e^{\beta y}\right)^4+\int_{-\infty}^{-R}E(y-\eta)\left(A+\delta e^{\beta\eta}\right)^4 d\eta\\&+\int_{-R}^0E(y-\eta)\left(A-\eps+\delta e^{\beta \eta}\right)^4d\eta+\int_{0}^{R_\delta}E(y-\eta)\left(A-\eps\theta+\delta e^{\beta \eta}\right)^4d\eta\\&+\int_{R_\delta}^\infty E(y-\eta)\left(A-\eps\theta+\delta e^{\beta \eta}\right)^4d\eta,\\
					\end{split}
				\end{equation}
				where we used the definition of $ c_0 $, the fact that $ A^4\leq\left(A+\delta e^{\beta\eta}\right)^4 $ as well as that $ e^{\beta R_\delta}\leq e^{\beta \eta} $ for $ \eta>R_\delta $. Expanding the power-law, ordering terms together and using that \begin{equation}\label{artanh}
					\int_\R E(\eta-y)e^{\alpha \eta}d\eta=\frac{\artanh(\alpha)}{\alpha}e^{\alpha y}
				\end{equation} for all $ |\alpha|<1 $, we compute
				\begin{equation}\label{subsol.2}
					-\psi_\delta''(y)+c\psi_\delta'(y)-\left(A-\psi_\delta(y)\right)^4+\int_\R E(y-\eta)\left(A-\psi_\delta(\eta)\right)^4d\eta\phantom{spacespacespacespacespacespacespace}
				\end{equation}
				\vspace{-0.4cm}\begin{alignat*}{3}
					\leq&-\frac{3c_0}{4}\beta\delta e^{\beta y}+\delta e^{\beta y}\left(\beta^2-\frac{c_0}{4}\beta+4A^3\left(\frac{\artanh(\beta)}{\beta}-1\right)\right)&(I^1_1)\\
					&+4A^3\eps\left[\theta-\int_{-R}^0 E(y-\eta) d\eta-\theta \int_0^\infty E(y-\eta)d\eta\right]&\quad\quad\quad\quad\quad\quad\quad\quad\quad\quad\quad\quad\quad(I^1_2)\\
					&+4A\eps^3\left[\theta^3-\int_{-R}^0 E(y-\eta) d\eta-\theta^3 \int_0^\infty E(y-\eta)d\eta\right]&(I^1_3)\\ 
					&-6A^2\eps^2\left[\theta^2-\int_{-R}^0 E(y-\eta) d\eta-\theta^2 \int_0^\infty E(y-\eta)d\eta\right]&(I^1_4)\\
					&-\eps^4\left[\theta^4-\int_{-R}^0 E(y-\eta) d\eta-\theta^4 \int_0^\infty E(y-\eta)d\eta\right]&(I^1_5)\\
					&+4\delta^3e^{3\beta y}\eps\left[\theta-\int_{-R}^0 E(y-\eta)e^{3\beta(\eta-y)} d\eta-\theta \int_0^\infty E(y-\eta)e^{3\beta(\eta-y)}d\eta\right]&(I^1_6)\\
					&+4\delta e^{\beta y}\eps^3\left[\theta^3-\int_{-R}^0 E(y-\eta)e^{\beta(\eta-y)} d\eta-\theta^3 \int_0^\infty E(y-\eta)e^{\beta(\eta-y)}d\eta\right]&(I^1_7)\\
					&-6\delta^2 e^{2\beta y}\eps^2\left[\theta^2-\int_{-R}^0 E(y-\eta)e^{2\beta(\eta-y)} d\eta-\theta^2 \int_0^\infty E(y-\eta)e^{2\beta(\eta-y)}d\eta\right]&\quad\quad(I^1_8)\\
					&-6A^2\delta^2e^{2\beta y}\left[1-\int_\R E(\eta-y)e^{2\beta(\eta-y)}d\eta\right]&(I^1_9)\\
					&-4A\delta^3e^{3\beta y}\left[1-\int_\R E(\eta-y)e^{3\beta(\eta-y)}d\eta\right]&(I^1_{10})\\
					&-\delta^4e^{4\beta y}\left[1-\int_\R E(\eta-y)e^{4\beta(\eta-y)}d\eta\right]&(I^1_{11})\\
					&+12A^2\eps\delta e^{\beta y}\left[\theta-\int_{-R}^0E(y-\eta)e^{\beta(\eta-y)}d\eta-\theta\int_0^\infty E(y-\eta)e^{\beta(\eta-y)}d\eta\right]&(I^1_{12})\\
					&+12A\eps\delta^2 e^{2\beta y}\left[\theta-\int_{-R}^0E(y-\eta)e^{2\beta(\eta-y)}d\eta-\theta\int_0^\infty E(y-\eta)e^{\beta(2\eta-y)}d\eta\right]&\phantom{spacespacespacesp}(I^1_{13})\\
					&-12A\eps^2\delta e^{\beta y}\left[\theta^2-\int_{-R}^0E(y-\eta)e^{\beta(\eta-y)}d\eta-\theta^2\int_0^\infty E(y-\eta)e^{\beta(\eta-y)}d\eta\right].&(I^1_{14})\\
				\end{alignat*}
				We proceed now estimating all different terms in \eqref{subsol.2}. By the choice of $ \beta $ in \eqref{def.beta} we have
				\begin{equation}\label{subsol.3}
					(I^1_1)\leq -\frac{3c_0}{4}\beta \delta e^{\beta y}.
				\end{equation}
				We now proceed estimating the terms $ (I^1_2) $-$ (I^1_5) $. Using the symmetry of $ E $, the definition of $ R $ and the choice of $ \theta=\frac{1}{5} $, we compute
				\begin{equation}\label{subsol.4}
					\begin{split}
						(I^1_2)=&4A^3\eps\left[\theta\int_y^\infty E(\eta)d\eta-\int_{y}^{R+y} E(\eta) d\eta\right]
						=4A^3\eps\left[\theta\int_{R+y}^\infty E(\eta)d\eta-(1-\theta)\int_{y}^{R+y} E(\eta) d\eta\right]\\
						=&4A^3\eps\left[-\theta\left(\int_{y}^{R+y} E(\eta) d\eta-\int_{R+y}^\infty E(\eta)d\eta\right)-3\theta\int_{y}^{R+y} E(\eta) d\eta\right]
						\leq -12 A^3\eps\theta\int_{y}^{R+y} E(\eta) d\eta.
					\end{split}
				\end{equation}
				Similarly, since $ 1-\theta^3= 124 \theta^3 $ we have
				\begin{equation}\label{subsol.5}
					\begin{split}
						(I^1_3)		\leq & -492 A\eps^3\theta^3\int_{y}^{R+y} E(\eta) d\eta.
					\end{split}
				\end{equation}
				Choosing $ \eps<2A\theta $, which by the choice of $ \theta $ implies that $ \eps<492 A\theta^3 $, we obtain
				\begin{equation}\label{subsol.6}
					(I^1_4)		= -6A^2\eps^2\left[\theta^2\int_y^\infty E(\eta)d\eta-\int_{y}^{R+y} E(\eta) d\eta\right]
					\leq 6A^2\eps^2\int_{y}^{R+y} E(\eta) d\eta
					\leq 12 A^3	\eps\theta \int_{y}^{R+y} E(\eta) d\eta
				\end{equation}
				and 
				\begin{equation}\label{subsol.7}
					\begin{split}
						(I^1_5)	=-\eps^4\left[\theta^4\int_y^\infty E(\eta)d\eta-\int_{y}^{R+y} E(\eta) d\eta\right]
						\leq \eps^4\int_{y}^{R+y} E(\eta) d\eta
						\leq 492 A \eps^3\theta^3 \int_{y}^{R+y} E(\eta) d\eta.
					\end{split}
				\end{equation}
				Hence, \eqref{subsol.4}-\eqref{subsol.7} imply
				\begin{equation}\label{subsol.8}
					(I^1_2)+(I^1_3)+(I^1_4)+(I^1_5)\leq 0.
				\end{equation}
				Besides the choice of $ \beta $ as in \eqref{def.beta} we use in the remaining estimates the fact that for $ y\in(0,R_\delta) $ the following holds true
				\[\delta e^{\beta y}\leq \delta e^{yR_\delta}=\eps\theta<\eps.\]
				Hence, we see that
				\begin{equation}\label{subsol.9}
					(I^1_6)\leq 4\delta^3 e^{3\beta y}\eps\theta\leq 4 \eps^3\delta e^{\beta y},\quad \quad
					(I^1_7)\leq 4\delta e^{\beta y}\eps^3\theta^3\leq 4 \eps^3\delta e^{\beta y},
				\end{equation}
				and
				\begin{multline}\label{subsol.11}
					(I^1_8)\leq 6\delta^2 e^{2\beta y}\eps^2\left(\int_{-(R+y)}^{-y} E(\eta)e^{2\beta\eta} d\eta+\theta^2 \int_{-y}^\infty E(\eta)e^{2\beta\eta}d\eta\right)\\\leq 6\delta^2 e^{2\beta y}\eps^2\int_{-\infty}^{\infty} E(\eta)e^{2\beta\eta} d\eta 6\delta^2 e^{2\beta y}\eps^2 \frac{\artanh(2\beta)}{2\beta}\leq 9\eps^3\delta e^{\beta y},
				\end{multline}
				where we used also that $ \beta \mapsto \frac{\artanh(\beta)}{\beta} $ is monotonically increasing. Finally, for the last six terms we estimate
				\begin{equation}\label{subsol.12}
					(I^1_9)= 6 A^2\delta^2 e^{2\beta y}\left(\frac{\artanh(2\beta)}{2\beta}-1\right)\leq 3A^2\delta^2e^{2\beta y}\leq3 \eps A^2 \delta e^{\beta y},
				\end{equation}
				\begin{equation}\label{subsol.13}
					(I^1_{10})= 4 A\delta^3 e^{3\beta y}\left(\frac{\artanh(3\beta)}{3\beta}-1\right)\leq 2 \eps^2 A \delta e^{\beta y},
					\quad\quad
					(I^1_{11})= \delta^4 e^{4\beta y}\left(\frac{\artanh(4\beta)}{4\beta}-1\right)\leq \frac{1}{2} \eps^3 \delta e^{\beta y},
				\end{equation}
				\begin{equation}\label{subsol.14.1}
					(I^1_{12})\leq 12 A^2 \eps \delta e^{\beta y},\quad\quad
					(I^1_{13})\leq 12 A \eps^2 \delta e^{\beta y},\quad \text{ and }\quad
					(I^1_{14})\leq 18 A \eps^2 \delta e^{\beta y}.
				\end{equation}
				Therefore, defining the constant in equation \eqref{def.eps}
				\begin{equation*}\label{constant.eps}
					C(\beta,A,\theta)=\min\left\{1,2A\theta,\frac{c_0\beta}{2(18+15A^2+32A)}\right\}
				\end{equation*}
				and combining the estimates \eqref{subsol.3},\eqref{subsol.8}-\eqref{subsol.14.1} we conclude that
				\begin{equation*}\label{subsol.15}
					-\psi_\delta''(y)+c\psi_\delta'(y)-\left(A-\psi_\delta(y)\right)^4+\int_\R E(y-\eta)\left(A-\psi_\delta(\eta)\right)^4d\eta\leq -\frac{c_0}{4}\beta\delta e^{\beta y}<0
				\end{equation*}
				for all $ y\in(0,R_\delta) $.\\
				
				We now notice that by the choice of $ \delta_0 $ we have $ \psi_{\delta_0}\leq \omega $ on $ \R $. In particular, since $ R_{\delta_0}=\frac{1}{\beta}\ln(2) $ the definition of $ \eps $ in \eqref{def.eps} implies that $ \psi_{\delta_0}(y)\leq \psi_{\delta_0} (0)=\frac{\eps\theta}{2}<\omega$ on $ [0,R_{\delta_0}] $ as well as $ \inf\limits_{y>0} \left(\omega-\psi_{\delta_0}\right)\geq\frac{\eps\theta}{2}> 0 $.  We remark that on $ \{y>0\} $ the functions $ \psi_\delta $ are continuous.
				
				We aim to show that $ \psi_\delta\leq \omega $ on $ [0,R_\delta] $ for all $ \delta\leq \delta_0 $. To this end we assume the contrary, i.e. we assume that there exists some $ 0<\delta<\delta_0 $ such that 
				\begin{equation}\label{subsol.17}
					\inf\limits_{y>0}\left( \omega-\psi_{\delta}\right)< 0.
				\end{equation}
				By construction this yields that $ \inf\limits_{y>0} \left(\omega-\psi_{\delta}\right)=\min\limits_{[0,R_\delta]}\left(\omega-\psi_{\delta}\right)< 0 $.
				The uniform continuity of $ [\delta,\delta_0]\ni\bar \delta \mapsto \psi_{\bar \delta} $ as functions on $ [0,R_\delta] $ and their monotonicity ($ \delta\mapsto\psi_{\delta} $ is increasing) imply that there exists
				\begin{equation}\label{subsol.18}
					\delta^*:=\sup\left\{\delta<\delta^*<\delta_0:\; \min\limits_{[0,R_\delta]} \left(\omega-\psi_{\delta^*}\right) <0\right\}
				\end{equation}
				such that 
				\[\min\limits_{[0,R_\delta]}\left( \omega-\psi_{\delta^*}\right) =\omega (y_0)-\psi_{\delta^*}=0\]
				for some $ y_0\in(0,R_{\delta^*}) $. Indeed, by construction $ \psi_{\delta^*}<\omega $ on $ y\geq R_{\delta^*} $ as well as $ \psi_{\delta^*}(0)=\eps\theta-\delta^*<\omega(0) $. Hence, we can apply the maximum principle for \eqref{eq.omega} at the point $ y_0 $ since on $ (0,R_{\delta^*}) $ the function $ \psi_{\delta^*} $ is smooth. We obtain the following contradiction
				\begin{equation}\label{subsol.16}
					\begin{split}
						0<&-\left(\omega-\psi_{\delta^*}\right)''(y_0)+c\left(\omega-\psi_{\delta^*}\right)'(y_0)-\left(A-\omega(y_0)\right)^4+\left(A-\psi_{\delta^*}\right)^4\\&+\int_\R E(y-\eta)\left(A-\omega\right)^4d\eta-\int_\R E(y-\eta)\left(A-\psi_{\delta^*}\right)^4d\eta\\
						\leq &\int_\R E(y-\eta)\left(A-\omega\right)^4d\eta-\int_\R E(y-\eta)\left(A-\psi_{\delta^*}\right)^4d\eta\leq 0,
					\end{split}
				\end{equation}
				since by construction $ 0\leq \psi_{\delta^*}\leq \omega $ for all $ y\in \R\setminus (0,R_{\delta^*}) $. Moreover, $  0\leq \psi_{\delta^*}\leq \omega $ for $ y\in (0,R_{\delta^*}) $. Thus, $ \left(A-\psi_{\delta^*}\right)\geq \left(A-\omega\right)\geq 0 $ on $ \R $.
				This contradiction implies that such $ \delta^* $ as in \eqref{subsol.18} and consequently such $ \delta $ satisfying \eqref{subsol.17} do not exist. Therefore we conclude that 
				\[\inf\limits_{y>0} \left(\omega-\psi_{\delta}\right)\geq 0\]
				for all $ \delta<\delta_0 $. 
				
				This implies that for all $ y\in[0,R_\delta] $ we can estimate $ w(y)\geq \eps\theta-\delta e^{\beta y} $ for all $ \delta<\delta_0 $. Thus, taking the pointwise limit as $ \delta\to 0 $ we conclude
				\[A-f(y)=w(y)\geq \eps\theta>0.\]
				This is clearly a contradiction to the assumption that $ \limsup\limits_{y\to\infty} f(y)=A $. Hence, $ f $ does not attain its supremum at $ +\infty $.
			\end{proof}
		\end{lemma}
		A similar argument shows that $ \bar f $, solution to \eqref{trav.wave.3.n.R}, does not attain its infimum at $ +\infty $, unless it is constant.
		\begin{lemma}\label{lemma.no.inf.infty}
			Let $ f $ solve \eqref{trav.wave.3.n.R}. Then $ f $ does not attain its infimum at $ +\infty $, i.e. $ \liminf\limits_{y\to\infty}f(y)<\inf\limits_{\R}f $, unless $ f $ is constant. 
			\begin{proof}
				We assume again that $ f $ is not constant and that $  \liminf\limits_{y\to\infty}f(y)=\inf\limits_{\R}f=:B>0  $. We consider the function $ \omega=f-B\geq 0 $. Moreover, since $ f $ is not constant also $ \omega>0 $ at the interior by Lemma \ref{lemma.no.sup.interior}. Hence, $ \omega $ solves
				\begin{equation}\label{eq.omega.2}
					-\omega''(y)+c\omega'(y)+(B+\omega(y))^4-\int_\R E(y-\eta)(B+\omega(\eta))^4d\eta=0.
				\end{equation}
				As we did in Lemma \ref{lemma.no.sup.infty} we will show that $ \omega(y)>0 $ as $ y\to\infty $, which is a contradiction to the assumption of $\liminf\limits_{y\to\infty}f(y)=\inf\limits_{\R}f$. We will consider the family of functions $ \psi_\delta $ defined as in \eqref{def.psi.delta} for $ \theta=\frac{1}{5} $, $ \beta $ as in \eqref{def.beta} and $ R $ defined in \eqref{def.R}. Moreover, we take $ \eps>0 $ satisfying
				\begin{equation}\label{def.eps.1}
					\eps<\min\left\{\min\limits_{\left[-R,\frac{\ln(2)}{\beta}\right]}\{f(y)-B\},\; C(\beta,B,\theta)\right\},
				\end{equation}
				where $ C(\beta,B,\theta)>0 $ is a constant that will be computed later. Finally, we consider $ \delta<\delta_0=\frac{\eps\theta}{2} $. \\
				
				By construction we see that $ \psi_\delta\leq \omega $ on $ \R\setminus (0,R_\delta) $. Moreover, it is important to remark that for $ y\in(0,R_\delta) $ the functions $ \psi_\delta $ are smooth, as well as $ \psi_\delta $ are continuous on $ y\geq 0 $. Hence, we can compute for $ y\in(0,R_\delta) $ the following 
				\begin{equation}\label{inf.subsol.1}
					\begin{split}
						-\psi_\delta''(y)&+c\psi_\delta'(y)+\left(B+\psi_\delta(y)\right)^4-\int_\R E(y-\eta)\left(B+\psi_\delta(\eta)\right)^4d\eta\\=&-\left(c-\frac{c_0}{4}\right)\beta\delta e^{\beta y}
						+\delta e^{\beta y}\left(\beta^2-\frac{c_0}{4}\beta\right)+\left(B+\eps \theta-\delta e^{\beta y}\right)^4-\int_{-\infty}^{-R}E(y-\eta)B^4 d\eta\\&-\int_{-R}^0E(y-\eta)\left(B+\eps-\delta e^{\beta \eta}\right)^4d\eta-\int_{0}^{R_\delta}E(y-\eta)\left(B+\eps\theta-\delta e^{\beta \eta}\right)^4d\eta\\&-\int_{R_\delta}^\infty E(y-\eta)\left(B+\eps\theta-\delta e^{\beta R_\delta}\right)^4d\eta\\
						\leq&-\frac{3c_0}{4}\beta\delta e^{\beta y}
						+\delta e^{\beta y}\left(\beta^2-\frac{c_0}{4}\beta\right)+\left(B+\eps \theta-\delta e^{\beta y}\right)^4-\int_{-\infty}^{-R}E(y-\eta)\left(B-\delta e^{\beta\eta}\right)^4 d\eta\\
						&-\int_{-R}^0E(y-\eta)\left(B+\eps-\delta e^{\beta \eta}\right)^4d\eta-\int_{0}^{R_\delta}E(y-\eta)\left(B+\eps\theta-\delta e^{\beta \eta}\right)^4d\eta\\&-\int_{R_\delta}^\infty E(y-\eta)\left(B+\eps\theta-\delta e^{\beta R_\delta}\right)^4d\eta.\\
					\end{split}
				\end{equation}
				As for \eqref{subsol.1} we used here the definition of $ c_0 $ as well as the fact that $ B^4\geq\left(B-\delta e^{\beta\eta}\right)^4 $ for $ \eta<-R $. Expanding the power-law, ordering terms together and using \eqref{artanh}, we compute
				\begin{equation}\label{inf.subsol.2}
					-\psi_\delta''(y)+c\psi_\delta'(y)+\left(B+\psi_\delta(y)\right)^4-\int_\R E(y-\eta)\left(B+\psi_\delta(\eta)\right)^4d\eta\phantom{spacespacespacespacespacespacespace}
				\end{equation}
				\vspace{-0.4cm}\begin{alignat*}{3}
					\leq&-\frac{3c_0}{4}\beta\delta e^{\beta y}+\delta e^{\beta y}\left(\beta^2-\frac{c_0}{4}\beta+4B^3\left(\frac{\artanh(\beta)}{\beta}-1\right)\right)&(I^2_1)\\
					&+4B^3\eps\left[\theta-\int_{-R}^0 E(y-\eta) d\eta-\theta \int_0^\infty E(y-\eta)d\eta\right]&(I^2_2)\\
					&+4B\eps^3\left[\theta^3-\int_{-R}^0 E(y-\eta) d\eta-\theta^3 \int_0^\infty E(y-\eta)d\eta\right]&(I^2_3)\\
					&+6B^2\eps^2\left[\theta^2-\int_{-R}^0 E(y-\eta) d\eta-\theta^2 \int_0^\infty E(y-\eta)d\eta\right]&(I^2_4)\\
					&+\eps^4\left[\theta^4-\int_{-R}^0 E(y-\eta) d\eta-\theta^4 \int_0^\infty E(y-\eta)d\eta\right]&(I^2_5)\\
					&-4\delta^3e^{3\beta y}\eps\left[\theta-\int_{-R}^0 E(y-\eta)e^{3\beta(\eta-y)} d\eta-\theta \int_0^\infty E(y-\eta)e^{3\beta(\eta-y)}d\eta\right]&(I^2_6)\\
					&-4\delta e^{\beta y}\eps^3\left[\theta^3-\int_{-R}^0 E(y-\eta)e^{\beta(\eta-y)} d\eta-\theta^3 \int_0^\infty E(y-\eta)e^{\beta(\eta-y)}d\eta\right]&(I^2_7)\\
					&+6\delta^2 e^{2\beta y}\eps^2\left[\theta^2-\int_{-R}^0 E(y-\eta)e^{2\beta(\eta-y)} d\eta-\theta^2 \int_0^{R_\delta} E(y-\eta)e^{2\beta(\eta-y)}d\eta-\theta^2 \int_{R_\delta}^\infty E(y-\eta)e^{2\beta\left(R_\delta-y\right)}d\eta\right]&(I^2_8)\\
					&+6B^2\delta^2e^{2\beta y}\left[1-\int_{-\infty}^{R_\delta} E(\eta-y)e^{2\beta(\eta-y)}d\eta-\int_{R_\delta}^\infty E(\eta-y)e^{2\beta\left(R_\delta-y\right)}d\eta\right]&(I^2_9)\\
					&-4B\delta^3e^{3\beta y}\left[1-\int_\R E(\eta-y)e^{3\beta(\eta-y)}d\eta\right]&(I^2_{10})\\
					&+\delta^4e^{4\beta y}\left[1-\int_{-\infty}^{R_\delta} E(\eta-y)e^{4\beta(\eta-y)}d\eta-\int_{R_\delta}^\infty E(\eta-y)e^{4\beta\left(R_\delta-y\right)}d\eta\right]&(I^2_{11})\\	&-12 B^2\eps \delta e^{\beta y}\left[\theta-\int_{-R}^{0} E(\eta-y)e^{\beta(\eta-y)}d\eta-\theta\int_{0}^\infty E(\eta-y)e^{\beta\left(\eta-y\right)}d\eta\right]&(I^2_{12})\\	&-12 B\eps^2 \delta e^{\beta y}\left[\theta^2-\int_{-R}^{0} E(\eta-y)e^{\beta(\eta-y)}d\eta-\theta^2\int_{0}^\infty E(\eta-y)e^{\beta\left(\eta-y\right)}d\eta\right]&(I^2_{13})	\\&+12 B\eps \delta^2 e^{2\beta y}\left[\theta-\int_{-R}^{0} E(\eta-y)e^{2\beta(\eta-y)}d\eta-\theta\int_{0}^{R_\delta} E(\eta-y)e^{2\beta\left(\eta-y\right)}d\eta-\theta\int_{R_\delta}^\infty E(\eta-y)e^{2\beta\left(R_\delta-y\right)}d\eta\right],&(I^2_{14})
				\end{alignat*}
				where we used that $ e^{n\beta R_\delta}\leq e^{n\beta \eta} $ for any $ n=1,2,3,4 $ and $ \eta\geq R_\delta $.\\
				
				Arguing as in the proof of \eqref{subsol.3}, \eqref{subsol.4} and \eqref{subsol.5} we see that also
				\begin{alignat*}{3}
					(I^2_1)\leq -\frac{3c_0}{4}\beta\delta e^{\beta y},\qquad\qquad(I^2_2)\leq 0\qquad\text{ and }\qquad (I^2_3)\leq 0.
				\end{alignat*}
				As we argued for $ \eqref{subsol.4} $ and \eqref{subsol.5} using that $ 1-\theta^2=24\theta^2 $ and $ 1-\theta^4=624\theta^4 $ we estimate 
				\begin{alignat}{3}\label{inf.subsol.4}
					(I^2_4)\leq -138B^2\eps^2\theta^2\int_{y}^{R+y} E(\eta) d\eta<0\qquad\text{ and }\qquad (I^2_5)\leq -623 \eps^4\theta^4\int_{y}^{R+y} E(\eta) d\eta<0.
				\end{alignat}
				Finally, estimating only the positive terms, using that $ \delta e^{\beta y}\leq \eps\theta $ for $ y<R_\delta $ and using the definition of $ \beta $ in \eqref{def.beta}, we compute
				\begin{alignat*}{5}\label{inf.subsol.5}
					(I^2_6)\leq& 4\eps\frac{\artanh(3\beta)}{3\beta}\delta^3 e^{3\beta y}\leq 6\eps^3 \delta e^{\beta y}, &\quad\quad&
					(I^2_7)\leq4\eps^3\frac{\artanh(\beta)}{\beta}\delta e^{\beta y}\leq 6\eps^3 \delta e^{\beta y},&\quad\quad&
					(I^2_8)\leq 6\eps^3 \delta e^{\beta y}\\
					(I^2_9)\leq & 6\eps B^2 \delta e^{\beta y},&\quad\quad&
					(I^2_{10})\leq 6\eps^2 B \delta e^{\beta y},&\quad\quad&
					(I^2_{11})\leq \eps^3 \delta e^{\beta y},\\
					(I^2_{12})\leq&  18B^2\eps \delta e^{\beta y},&\quad\quad&
					(I^2_{13})\leq  18 B\eps^2 \delta e^{\beta y},&\text{ and }\quad&
					(I^2_{14})\leq  12 B\eps^2 \delta e^{\beta y}.
				\end{alignat*}
				Hence, choosing in the definition \eqref{def.eps.1} of $ \eps $ the constant $ C(\beta,B,\theta)>0 $ as 
				\begin{equation*}\label{def.constant}
					C(\beta,B,\theta)=\min\left\{1,\frac{c_0\beta}{2(18+36B+24B^2)}\right\},
				\end{equation*}
				we conclude that
				\[-\psi_\delta''(y)+c\psi_\delta'(y)+\left(B+\psi_\delta(y)\right)^4-\int_\R E(y-\eta)\left(B+\psi_\delta(\eta)\right)^4d\eta=-\frac{c_0}{4}\beta\delta e^{\beta y}<0.\]
				We see once more that by the choice of all the parameters we have $ \psi_{\delta_0}\leq \omega $ on $ \R $ as well as $ \psi_{\delta_0}<\omega $ on $ [0,R_{\delta_0}] $. Moreover, for all $ \delta<\delta_0 $ it is true that $ \psi_\delta\leq \omega $ on $ \R\setminus (0,R_\delta) $ as well as $ \psi_R(0)<\omega $ and $ \psi_\delta(R_\delta)<\omega $. Hence, arguing as in the proof of Lemma \ref{lemma.no.sup.infty} we see that assuming the existence of some $ \delta<\delta_0 $ with \[\inf\limits_{y>0} \left(\omega-\psi_{\delta}\right)=\min\limits_{[0,R_\delta]}\left(\omega-\psi_{\delta}\right)< 0\]
				there exists also some $ \delta<\delta^*<\delta_0 $ defined by 
				$ \delta^*:=\sup\left\{\delta<\delta^*<\delta_0:\; \min\limits_{[0,R_\delta]} \left(\omega-\psi_{\delta^*} \right)<0\right\} $
				such that 
				\[\min\limits_{[0,R_\delta]}\left( \omega-\psi_{\delta^*} \right)=\omega (y_0)-\psi_{\delta^*}=0\]
				for some $ y_0\in(0,R_{\delta^*}) $. However, the application of the maximum principle for the equation \eqref{eq.omega.2} to the functions $ \omega$ and $\psi_{\delta^*} $ yields as in \eqref{subsol.16} the contradiction $ 0<-\int_\R E(y-\eta)\left[(B+\omega(\eta))^4-(B+\psi_{\delta^*}(\eta))^4\right] d\eta<0 $. Therefore, we conclude that  
				\[\inf\limits_{y>0} \left(\omega-\psi_{\delta}\right)\geq 0\]
				for all $ \delta<\delta_0 $, so that  $ w(y)\geq \eps\theta-\delta e^{\beta y} $ for all $ \delta<\delta_0 $ and all $ y\in[0,R_\delta] $.
				Thus, taking the pointwise limit as $ \delta\to 0 $ we establish
				\[f(y)-B=w(y)\geq \eps\theta>0,\]
				which contradicts the assumption that $ \liminf\limits_{y\to\infty} f(y)=B$. Hence, $ f $ does not attain its infimum at $ +\infty $.
			\end{proof}
		\end{lemma}
		\subsection{The $ \omega $-limits of the traveling waves are constant}
		Lemma \ref{lemma.no.sup.interior}, Lemma \ref{lemma.no.sup.infty} and Lemma \ref{lemma.no.inf.infty} imply that the limit function $ \bar f $ solving \eqref{trav.wave.3.n.R} is either constant or it takes the supremum and infimum at $ -\infty $, i.e. 
		\[\inf\limits_\R \bar f=\liminf\limits_{y\to-\infty} \bar f(y)<\limsup\limits_{y\to-\infty} \bar f(y)=\sup\limits_\R \bar f.\]
		We will show that $ \bar f $ is constant, showing that $ \liminf\limits_{y\to-\infty} \bar f(y)=\limsup\limits_{y\to-\infty} \bar f(y) .$
		We start proving the following Theorem, which is a fundamental stability result.
		\begin{theorem}\label{thm.stability}
			Let $ f $ solve \eqref{trav.wave.3.n.R} for $ 0<\lambda<T_M $. Then there exists an $ \eps_0=\eps_0(T_M,\lambda,c)>0 $ such that for all $ \eps<\eps_0 $ there exists $ L_0(\eps,T_m,\lambda,c)>0 $ with the property that if \[\osc\limits_{[-L,L]} f<\eps\] then also \[\osc\limits_{[L,\infty)} f<3\eps\] for all $ L>L_0 $.
			\begin{proof}
				Let us assume that $ f $ satisfy $ \osc\limits_{[-L,L]} f<\eps $ for some $ L>0 $ and some $ \eps>0 $. We show that for $ \eps>0 $ small enough and for $ L>0 $ large enough this assumption implies $ \osc\limits_{[L,\infty)} f<3\eps $. In the course of the proof we will also define $ \eps_0 $ and $ L_0(\eps) $.\\
				
				If $ \osc\limits_{[-L,L]} f<\eps $, then the maximum $ \max\limits_{[-L,L]} f=: M_L<T_M $ and the minimum $ \min\limits_{[-L,L]} f=:m_L>\lambda $ satisfy \[M_L-m_L<\eps.\]
				We now construct two suitable families of subsolutions and supersolutions, for which the maximum principle will show the claim in a similar way as in the proofs of Lemma \ref{lemma.no.sup.infty} and of Lemma \ref{lemma.no.inf.infty}. Let us consider the following functions
				\begin{equation}\label{def.psi.delta.1}
					\psi^L_\delta(y)=m_L-\eps+\begin{cases}
						-(m_L-\eps)&\quad y<-L\\
						\eps-\delta e^{\beta(y-L)}&\quad -L\leq y<L\\
						\eps\theta-\delta e^{\beta(y-L)}&\quad y\in[L,R_\delta]\\
						-(m_L-\eps)&\quad y>R_\delta,
					\end{cases}
				\end{equation}
				where $ R_\delta=\frac{1}{\beta}\ln\left(\frac{m_L-(1-\theta)\eps}{\delta}\right)+L $ is so defined that $ \psi^L_\delta $ is continuous on $ [L,\infty) $. Moreover, we notice that $ \psi^L_\delta $ is smooth in $ (L,R_\delta) $. We consider $ \eps<\lambda $ as well as $ \delta<m_L-\lambda $, so that $ m_L-\eps>0 $ and $ R_\delta>L $, and we study the family of functions $ \psi_\delta^L $ for $ \delta<\delta_0 $, where $ \delta_0(\eps,L)>0 $ will be specified later. We also fix $ \theta=\frac{1}{5} $ and $ c_0=\min\{c,1\} $. In addition we choose \begin{equation}\label{def.beta.2}
					\beta<\min\left\{\frac{1}{24},\frac{c_0}{2}, \frac{c_0}{2T_M^3} , \left(\frac{10}{8}\frac{c_0}{77 T_M^3}\right)^{24}\right\}
				\end{equation} satisfying also
				\begin{equation}\label{def.beta.1}
					\beta^2-\frac{c_0}{2}\beta+4T_M^3\left(\frac{\artanh(\beta)}{\beta}-1\right)+4T_M^3\left(\frac{\artanh(3\beta)}{3\beta}-1\right)\leq 0.
				\end{equation}
				For $ c_1=6\artanh\left(\frac{1}{6}\right) $ we take $ \eps<\eps_1 $ defined by
				\begin{equation}\label{def.eps1}
					\eps_1=\min\left\{1,\lambda,\frac{c_0\beta}{8(4c_1+16c_1T_M^2+12T_M^2+12T_Mc_1+6T_M)}\right\}.
				\end{equation}
				We also consider a fixed $ L>L_1(\eps) $ satisfying
				\begin{equation}\label{def.L.1}
					\left(1+\frac{T_M}{\eps\theta}+\frac{T_M^3}{(\eps\theta)^3}\right)\int_{L_1+y}^\infty E(z)dz<\int_{y-L_1}^{L_1+y} E(z)dz \text{ for all }y>L_1,\;e^{-L_1}<\beta^2 \text{ and } L_1>\frac{1}{\sqrt{\beta}}.
				\end{equation}
				We remark that $ \beta $ given by \eqref{def.beta.1} and $ L_1 $ defined by \eqref{def.L.1} are well-defined. For $ \beta $ one argues similarly as for \eqref{def.beta}, while for $ L_1 $ we need to adapt the proof for \eqref{def.R}. This adaptation however is easy. As we proved in \eqref{def.R}, one can readily see that for any $ A>0 $
				\[A\int_{N+y}^\infty E(z)dz<\int_0^{N+y} E(z)dz \text{ for }  N=N(A)>0 \text{ large enough and for all }  y>0 .\]Taking $ A=1+\frac{T_M}{\eps\theta}+\frac{T_M^3}{(\eps\theta)^3} $ and $ L_1=\frac{N(A)}{2} $, we conclude \eqref{def.L.1}. Moreover, we remark that the function
				\[N\mapsto A\int_{N+y}^\infty E(z)dz-\int_0^{N+y} E(z)dz \]
				is monotonically decreasing for $ N>N(A) $. Hence, \eqref{def.L.1} holds also for all $ L>L_1(\eps) $.
				
				Finally, we set 
				\begin{equation}\label{def.delta.0}
					\delta_0=(m_L-(1-\theta)\eps)e^{-\beta R_\eps}>0,
				\end{equation}
				where $ R_\eps $ is the distance such that \begin{equation}\label{def.R.eps}
					f(y)-m_L\geq -\frac{1-\theta}{2}\eps \quad\text{ for all }\quad y\in [L-R_\eps,L+R_\eps].
				\end{equation}It is important to notice that $ R_\eps $ is independent of $ L $. This can be proved using the uniform continuity of $ f $, according to which there exists $ R_\eps $ such that $ f(y)-f(x)\geq -\frac{1-\theta}{2}\eps $ for all $ |x-y|<R_\eps $. Finally, $ x=L $ and $ f(L)\geq m_L $ implies \eqref{def.R.eps}. Thus, with $ \delta_0 $ defined in \eqref{def.delta.0} we see that 
				\[R_{\delta_0}= \frac{1}{\beta}\ln\left(e^{\beta R_\eps}\right)=R_\eps+L.\]
				Hence, for all $ y\in[L,R_{\delta_0}] $ we have by construction $ \psi_{\delta_0}^L\leq m_L -(1-\theta)\eps $ as well as $ f(y)\geq m_L-\frac{1-\theta}{2}\eps. $ This implies \[f(y)-\psi_{\delta_0}^L\geq\frac{1-\theta}{2}\eps>0. \]for $ y\in[L,R_{\delta_0}] $. \\
				
				Moreover, by definition we know that $ \psi_\delta^L<f(y) $ for $y\in\R\setminus(L,R_\delta) $ and for all $ \delta\leq \delta_0 $. We remark also that $ \psi_\delta^L(L)=m_L-(1-\theta)\eps<m_L\leq f(L) $ as well as $ \psi_\delta^L(R_\delta)=0<f(R_\delta) $. Hence, $ \psi_{\delta_0}^L<f $ in $ \R $.\\
				
				We will now show that $ \psi_\delta^L $ is a subsolution to the equation \eqref{trav.wave.3.n.R} for $ y\in(L,R_\delta) $, where the function is also smooth.\\
				
				Let us first of all assume that $ y\in\left[R_\delta-\frac{1}{\sqrt{\beta}},R_\delta\right)\cap (L,R_\delta) $. We compute \[\delta e^{\beta(y-L)}\geq \delta e^{\beta(R_\delta-L)}e^{-\sqrt{\beta}}=(m_L-(1-\theta)\eps)e^{-\sqrt{\beta}}.\]
				Hence, 
				\[0\leq\psi_\delta^L(y)\leq (m_L-(1-\theta)\eps)\left(1-e^{-\sqrt{\beta}}\right).\]
				This implies that 
				\begin{equation}\label{subsol.1case.1}
					\begin{split}
						-\left(\right.&\left.\psi_\delta^L\right)''(y)+c\left(\psi_\delta^L\right)'(y)+\left(\psi_\delta^L(y)\right)^4-\int_\R E(\eta-y)\left(\psi_\delta^L(\eta)\right)^4d\eta\\
						<&(m_L-(1-\theta)\eps) (\beta^2-c\beta)+(m_L-(1-\theta)\eps)^4\left(1-e^{-\sqrt{\beta}}\right)^4
						\leq  (m_L-(1-\theta)\eps)^4\left[\frac{\beta^2-c_0\beta}{(m_L-(1-\theta)\eps)^3}+\beta^2\right]\\
						\leq&(m_L-(1-\theta)\eps)^4\left[-\frac{c_0\beta}{2(m_L-(1-\theta)\eps)^3}+\beta^2\right]
						\leq(m_L-(1-\theta)\eps)^4\left[-\frac{c_0\beta}{2T_M^3}+\beta^2\right]<0.\\
					\end{split}
				\end{equation}
				We used besides the definition of $ \beta $ in \eqref{def.beta.2} also that $ (m_L-(1-\theta)\eps)\leq T_M $, $ c\leq c_0 $ as well as $ 1-e^{-|x|}\leq |x| $.
				
				It remains to show that $ \psi_\delta^L $ is a subsolution also for $ y\in\left(L,R_\delta-\frac{1}{\sqrt{\beta}}\right) $. Without loss of generality we assume $ \left(L,R_\delta-\frac{1}{\sqrt{\beta}}\right) \neq\emptyset $, since this is true for $ \delta $ small enough. Moreover, for all $ \delta\leq\delta_0 $ with  $\left[R_\delta-\frac{1}{\sqrt{\beta}},R_\delta\right)\cap (L,R_\delta)= (L,R_\delta) $ estimate \eqref{subsol.1case.1} gives the result about $ \psi_\delta^L $ being a subsolution. We collect many estimates similar to the ones made for \eqref{subsol.2} and \eqref{inf.subsol.2}.
				For the following computation we use $ c\geq c_0 $ and that $ e^{\beta R_\delta}<e^{\beta \eta} $ for $ \eta>R_\delta $, we expand the power law, and we order similar terms together.
				
				\begin{equation}\label{subsol.2case.1}
					-\left(\psi_\delta^L\right)''(y)+c\left(\psi_\delta^L\right)'(y)+\left(\psi_\delta^L(y)\right)^4-\int_\R E(\eta-y)\left(\psi_\delta^L(\eta)\right)^4d\eta\\\phantom{spacespacespacespacespacespacespace}
				\end{equation}
				\vspace{-0.65cm}\begin{alignat*}{3}
					\leq&-\frac{c_0}{2}\beta\delta e^{\beta (y-L)}+\delta e^{\beta (y-L)}\left(\beta^2-\frac{c_0}{4}\beta-4(m_L-\eps)^3+4(m_L-\eps)^3\int_{-L}^\infty E(\eta-y) e^{\beta(\eta-y)}d\eta\right)&\quad\quad\quad\quad\quad(I^3_1)\\
					&-4(m_L-\eps)\delta^3 e^{3\beta (y-L)}\left(1-\int_{-L}^\infty E(\eta-y) e^{3\beta(\eta-y)}d\eta\right)&(I^3_2)\\\end{alignat*}\begin{alignat*}{3}
					&+6(m_L-\eps)^2\delta^2 e^{2\beta (y-L)}\left(1-\int_{-L}^{R_\delta} E(\eta-y) e^{2\beta(\eta-y)}d\eta\right)-6(m_L-\eps)^2\delta^2\int_{R_\delta}^\infty E(\eta-y) e^{2\beta(R_\delta-L)}&(I^3_3)\\
					&+\delta^4 e^{4\beta (y-L)}\left(1-\int_{-L}^{R_\delta} E(\eta-y) e^{4\beta(\eta-y)}d\eta\right)-\delta^4\int_{R_\delta}^\infty E(\eta-y) e^{4\beta(R_\delta-L)}&(I^3_4)\\
					&+ 4(m_L-\eps)^3\left[\eps\theta+\int_{-\infty}^{-L}E(\eta-y)(m_L-\eps)d\eta-\eps\int_{-L}^LE(\eta-y)d\eta-\eps\theta\int_{L}^\infty E(\eta-y)d\eta\right]&(I^3_5)\\
					&+ 4(m_L-\eps)\left[(\eps\theta)^3+\int_{-\infty}^{-L}E(\eta-y)(m_L-\eps)^3d\eta-\eps^3\int_{-L}^LE(\eta-y)d\eta-(\eps\theta)^3\int_{L}^\infty E(\eta-y)d\eta\right]&(I^3_6)\\
					&+ 6(m_L-\eps)^2\left[(\eps\theta)^2-\int_{-\infty}^{-L}E(\eta-y)(m_L-\eps)^2d\eta-\eps^2\int_{-L}^LE(\eta-y)d\eta-(\eps\theta)^2\int_{L}^\infty E(\eta-y)d\eta\right]&(I^3_7)\\
					&+(\eps\theta)^4-\int_{-\infty}^{-L}E(\eta-y)(m_L-\eps)^4d\eta-\eps^4\int_{-L}^LE(\eta-y)d\eta-(\eps\theta)^4\int_{L}^\infty E(\eta-y)d\eta&(I^3_8)\\
					&-4\eps^3\delta e^{\beta (y-L)}\left(\theta^3-\int_{-L}^L E(\eta-y) e^{\beta(\eta-y)}d\eta-\theta^3\int_{L}^\infty E(\eta-y) e^{\beta(\eta-y)}d\eta\right)&(I^3_9)\\
					&-4\eps\delta^3 e^{3\beta (y-L)}\left(\theta^3-\int_{-L}^L E(\eta-y) e^{3\beta(\eta-y)}d\eta-\theta\int_{L}^\infty E(\eta-y) e^{3\beta(\eta-y)}d\eta\right)&(I^3_{10})\\
					&+6\eps^2\delta^2 e^{2\beta (y-L)}\left(\theta^2-\int_{-L}^L E(\eta-y) e^{2\beta(\eta-y)}d\eta-\theta^2\int_{L}^{R_\delta} E(\eta-y) e^{2\beta(\eta-y)}d\eta\right)&(I^3_{11})\\&-6(\eps\theta\delta)^2e^{2\beta(R_\delta-L)}\int_{R_\delta}^\infty E(\eta-y)d\eta&\\
					&-12(m_L-\eps)^2\eps\delta e^{\beta (y-L)}\left(\theta-\int_{-L}^L E(\eta-y) e^{\beta(\eta-y)}d\eta-\theta\int_{L}^{\infty} E(\eta-y) e^{\beta(\eta-y)}d\eta\right)&(I^3_{12})\\
					&-12(m_L-\eps)\eps^2\delta e^{\beta (y-L)}\left(\theta^2-\int_{-L}^L E(\eta-y) e^{\beta(\eta-y)}d\eta-\theta^2\int_{L}^{\infty} E(\eta-y) e^{\beta(\eta-y)}d\eta\right)&(I^3_{13})\\
					&+12(m_L-\eps)\eps\delta^2 e^{2\beta (y-L)}\left(\theta-\int_{-L}^L E(\eta-y) e^{2\beta(\eta-y)}d\eta-\theta\int_{L}^{R_\delta} E(\eta-y) e^{2\beta(\eta-y)}d\eta\right)&\phantom{justab}(I^3_{14})\\
					&-12 (m_L-\eps)\eps\delta^2 e^{2\beta (R_\delta-L)}\int_{R_\delta}^\infty E(\eta-y)d\eta.&
				\end{alignat*}
				Next, using \eqref{artanh} and estimating $ (m-\eps)\leq T_M $ as well as $ \delta e^{\beta (y-L)}\leq m-(1-\theta)\eps\leq T_M $ we can compute
				\begin{equation}\label{subsol.2case.2}
					\begin{split}
						(I^3_1)+(I^3_2)\leq &-\frac{c_0}{2}\beta\delta e^{\beta (y-L)}+\delta e^{\beta (y-L)}\left[\beta^2-\frac{c_0}{4}\beta+4T_M^3\left(\frac{\artanh(\beta)}{\beta}-1\right)+4T_M^3\left(\frac{\artanh(3\beta)}{3\beta}-1\right)\right]\\\leq&-\frac{c_0}{2}\beta\delta e^{\beta (y-L)}
					\end{split}
				\end{equation}
				by the choice of $ \beta $ as in \eqref{def.beta.1}.
				Next we consider $ (I^3_7) $ and $ (I^3_8) $. Here we use \eqref{subsol.4}, \eqref{subsol.5} and \eqref{inf.subsol.4}, obtaining
				\begin{equation}\label{subsol.2case.3}
					(I^3_7)\leq 6(m_L-\eps)^2\left[(\eps\theta)^2-\eps^2\int_{-L}^LE(\eta-y)d\eta-(\eps\theta)^2\int_{L}^\infty E(\eta-y)d\eta\right]\leq 0
				\end{equation}
				as well as
				\begin{equation}\label{subsol.2case.4}
					(I^3_8)\leq(\eps\theta)^4-\eps^4\int_{-L}^LE(\eta-y)d\eta-(\eps\theta)^4\int_{L}^\infty E(\eta-y)d\eta\leq 0
				\end{equation}
				Using $ \theta=\frac{1}{5}$ and $ L\geq L_1 $ as defined in \eqref{def.L.1}, we also estimate
				\begin{equation}\label{subsol.2case.5}
					\begin{split}
						(I^3_5)\leq & 4(m_L-\eps)^3\left[\eps\theta\int_{-\infty}^{L}E(\eta-y)d\eta+T_M\int_{-\infty}^{-L}E(\eta-y)d\eta-\eps\int_{-L}^LE(\eta-y)d\eta\right]\\
						\leq &  4(m_L-\eps)^3\left[\eps\theta\int_{-\infty}^{-L}E(\eta-y)d\eta+T_M\int_{-\infty}^{-L}E(\eta-y)d\eta-4\eps\theta\int_{-L}^LE(\eta-y)d\eta\right]\\
						=&4(m_L-\eps)^3\eps\theta\left[\left(1+\frac{T_M}{\eps\theta}\right)\int_{-\infty}^{-L}E(\eta-y)d\eta-4\int_{-L}^LE(\eta-y)d\eta\right]\leq 0
					\end{split}
				\end{equation}
				and
				\begin{equation}\label{subsol.2case.6}
					\begin{split}
						(I^3_6)\leq & 4(m_L-\eps)\left[(\eps\theta)^3\int_{-\infty}^{L}E(\eta-y)d\eta+T_M^3\int_{-\infty}^{-L}E(\eta-y)d\eta-\eps^3\int_{-L}^LE(\eta-y)d\eta\right]\\
						\leq &4(m_L-\eps)(\eps\theta)^3\left[\left(1+\frac{T_M^3}{(\eps\theta)^3}\right)\int_{-\infty}^{-L}E(\eta-y)d\eta-124\int_{-L}^LE(\eta-y)d\eta\right]\leq 0.
					\end{split}
				\end{equation}
				Using $ \frac{\artanh(n\beta)}{n\beta}\leq c_1 $ for $ n\leq 4 $, $ \theta=\frac{1}{5} $ as well as the estimate $ \delta e^{\beta(y-L)}\leq T_M $ for $ y\leq R_\delta $ we compute furthermore
				\begin{equation}\label{subsol.2case.7}
					(I^3_9)\leq 4 \eps^3 \delta e^{\beta(y-L)} \frac{\artanh(\beta)}{\beta}\leq 4\eps^3 c_1 \delta e^{\beta(y-L)},
					\quad
					(I^3_{10})\leq 4 \eps \delta^3 e^{3\beta(y-L)} \frac{\artanh(3\beta)}{3\beta}\leq 4\eps T_M^2 c_1 \delta e^{\beta(y-L)},
				\end{equation} 
				\begin{equation}\label{subsol.2case.9}
					(I^3_{11})\leq  6\eps^2 T_M \delta e^{\beta(y-L)},
					\quad \text{ and }\quad
					(I^3_{12})+(I^3_{13})+(I^3_{14})\leq  12 \eps T_M \delta e^{\beta(y-L)}\left(T_M(c_1+1)+c_1\right).
				\end{equation}
				Finally, we have to estimate the remaining terms $ (I^3_3) $ and $ (I^3_4) $. To this end we recall that we are considering the case for which $ R_\delta-y\geq \frac{1}{\sqrt{\beta}} $ and that we have chosen $ L>L_1 $ such that $ e^{-L_1}\leq \beta^2 $. Additionally, we also use that \begin{equation}\label{est.beta}
					e^{-\frac{1}{\sqrt{\beta}}}<\beta^{\frac{5}{4}} \text{ for all } \beta>0.
				\end{equation}
				We hence estimate
				\begin{equation}\label{subsol.2case.10}
					\begin{split}
						(I_3)\leq & 6(m_L-\eps)^2 \delta^2 e^{2\beta(y-L)}\left(1-\int_{-L}^{R_\delta}E(\eta-y)e^{2\beta(\eta-y)}d\eta\right)\\
						=&6(m_L-\eps)^2 \delta^2 e^{2\beta(y-L)}\left(1-\frac{\artanh(2\beta)}{2\beta}+\int_{-\infty}^{-L}E(\eta-y)e^{2\beta(\eta-y)}d\eta+\int_{R_\delta}^{\infty}E(\eta-y)e^{2\beta(\eta-y)}d\eta\right)\\
						\leq&6 T_M^3 \delta e^{\beta(y-L)}\left(\int_{-\infty}^{-(L+y)}E(\eta)e^{2\beta\eta}d\eta+\int_{R_\delta-y}^{\infty}E(\eta)e^{2\beta\eta}d\eta\right)\\\end{split}\end{equation}\begin{equation*}\begin{split}
						\leq &6T_M^3 \delta e^{\beta(y-L)}\left(\frac{e^{-(L+y)}}{2}+\int_{\frac{1}{\sqrt{\beta}}}^{\infty}\frac{e^{-(1-2\beta)\eta}}{2}d\eta\right)
						\leq 6T_M^3 \delta e^{\beta(y-L)}\left(\frac{\beta^2}{2}+\frac{3}{5}e^{-\frac{5}{6\sqrt{\beta}}}\right)\quad\quad\quad\quad\quad\quad\quad\quad\\
						\leq &6T_M^3 \delta e^{\beta(y-L)}\left(\frac{\beta^2}{2}+\frac{3}{5}\beta \beta^{\frac{1}{24}}\right).
					\end{split}
				\end{equation*}
				We also used in the second inequality that $ \frac{\artanh(a)}{a}\geq 1 $, as well as in the third inequality the estimate $ e^{2\beta(\eta)}\leq 1 $ for $ \eta\leq -(L+y)\leq 0 $, the estimate \eqref{important.est.E} and the inequality $ E(z)\leq \frac{e^{-|z|}}{2} $ for $ |z|>1 $  since $ \beta^{-\frac{1}{2}}>1 $. For the fourth inequality we used $ (1-2\beta)\geq \frac{5}{6} $ since $ \beta<\frac{1}{24}<\frac{1}{12} $ and we concluded the fifth estimate with \eqref{est.beta}.
				In a very similar way, using again that $ (1-4\beta)\geq \frac{5}{6} $ since $ \beta<\frac{1}{24} $, we also have the estimate
				\begin{equation}\label{subsol.2case.11}
					\begin{split}
						(I_4)\leq &  \delta^4 e^{4\beta(y-L)}\left(1-\int_{-L}^{R_\delta}E(\eta-y)e^{4\beta(\eta-y)}d\eta\right)\\
						=& \delta^4 e^{4\beta(y-L)}\left(1-\frac{\artanh(4\beta)}{4\beta}+\int_{-\infty}^{-L}E(\eta-y)e^{4\beta(\eta-y)}d\eta+\int_{R_\delta}^{\infty}E(\eta-y)e^{4\beta(\eta-y)}d\eta\right)\\
						\leq& T_M^3 \delta e^{\beta(y-L)}\left(\int_{-\infty}^{-(L+y)}E(\eta)e^{4\beta\eta}d\eta+\int_{R_\delta-y}^{\infty}E(\eta)e^{4\beta\eta}d\eta\right)
						\leq T_M^3 \delta e^{\beta(y-L)}\left(\frac{\beta^2}{2}+\frac{3}{5}\beta \beta^{\frac{1}{24}}\right).\\
				\end{split}\end{equation}
				Putting now together all estimates \eqref{subsol.2case.2}-\eqref{subsol.2case.9} and \eqref{subsol.2case.10}-\eqref{subsol.2case.11}, and using that $ \beta<\left(\frac{10}{8}\frac{c_0}{77 T_M^3}\right)^{24} $ and $ \eps<1 $ we conclude for $ y\in\left(L,R_\delta-\frac{1}{\sqrt{\beta}}\right) $
				\begin{equation}\label{subsol.2case.12}
					\begin{split}
						-\left(\psi_\delta^L\right)''(y)&+c\left(\psi_\delta^L\right)'(y)+\left(\psi_\delta^L(y)\right)^4-\int_\R E(\eta-y)\left(\psi_\delta^L(\eta)\right)^4d\eta\\
						\leq &\delta e^{\beta (y-L)}\left(-\frac{c_0}{2}\beta+\eps\left(4c_1+16c_1T_M^2+12T_M^2+12c_1T_M+6T_M\right)+\frac{77T_M^3}{10}\beta\beta^{\frac{1}{24}}\right)\\
						\leq &\delta e^{\beta (y-L)}\left(-\frac{c_0}{2}\beta+\frac{c_0}{8}\beta+\frac{c_0}{8}\beta\right)=-\frac{c_0}{4}\beta\delta e^{\beta (y-L)}<0,\\
					\end{split}
				\end{equation}
				where at the end we used the choice of $ \eps<\eps_1 $ and of $ \beta $ according to \eqref{def.eps1} and \eqref{def.beta.2}, respectively.\\
				
				Estimates \eqref{subsol.1case.1} and \eqref{subsol.2case.12} show that for all $ \delta<\delta_0 $ and for all $ y\in(L,R_\delta) $ the functions $ \psi_\delta^L $ are subsolutions, i.e.
				\[-\left(\psi_\delta^L\right)''(y)+c\left(\psi_\delta^L\right)'(y)+\left(\psi_\delta^L(y)\right)^4\leq 0.\]
				Since by construction $ \psi_{\delta_0}^L<f $ in $ \R $ with $  \psi_{\delta_0}^L-f\leq -\frac{1-\theta}{2}\eps<0 $ for all $ y\geq L $, as well as $ \psi_{\delta}^L\leq f $ for $ y\in\R\setminus(L,R_\delta) $ with $ \psi_\delta^L\left.\right|_{\{L,R_\delta\}}<f\left.\right|_{\{L,R_\delta\}} $, applying the maximum principle in the same way as we did in the proof of Lemma \ref{lemma.no.sup.infty} and Lemma \ref{lemma.no.inf.infty} and using the uniform continuity and the increasing monotonicity of $ \delta\mapsto \psi_\delta^L $ on compact sets as well as the fact that  $ \psi_\delta^L $ are subsolutions on $ (L,R_\delta) $ we conclude that 
				\[\psi_\delta^L(y)\leq f(y) \quad \text{ for all }y\in\R \text{ and }\delta<\delta_0.\]
				Hence, for any $ y>L $ we have for $ \delta<\delta_0 $ small enough
				\[f(y)\geq m_L-(1-\theta)\eps -\delta e^{\beta(y-L)}.\]
				Taking $ \delta\to 0 $ and thus $ R_\delta\to\infty $ yields
				\begin{equation}\label{minimum}
					f(y)\geq m_L-(1-\theta)\eps\quad\text{ for all }y>L.
				\end{equation}
				In a similar way we show now that $ f(y)\leq M_L+(1-\theta)\eps $ for $ y>L $. We consider a similar family of functions called $ \{\psi_\gamma^L\} $ which we will prove to be supersolutions. In this case we define
				\begin{equation}\label{def.psi.gamma}
					\psi_\gamma^L(y)=M_L+\eps+\begin{cases}
						2T_M-(M_L+\eps)&\quad y<-L\\
						\gamma e^{\zeta(y-L)}-\eps&\quad -L\leq y<L\\
						\gamma e^{\zeta (y-L)}-\eps\theta&\quad L\leq y\leq R_\gamma\\
						2T_M-(M_L+\eps)&\quad y>R_\gamma,
					\end{cases}
				\end{equation}
				where $ R_\gamma=\frac{1}{\zeta}\ln\left(\frac{2T_M-(M_L+(1-\theta)\eps)}{\gamma}\right)+L $. We consider also $ \eps<T_M $ and $ \gamma<T_M-M_L $, so that $ 2T_M-M_L-\eps>0 $ as well as $ R_\gamma>L $. We remark that since $ f $ does not take supremum and infimum at the interior, $ T_M-M_L>0 $. Moreover, we notice that this family of functions is continuous on $( L,\infty) $ as well as smooth on $ (L,R_\gamma) $. For a $ \gamma_0(\eps,L)>0 $ defined later we study the family of functions $ \left\{\psi_\gamma^L\right\} $ for $ \gamma<\gamma_0 $. We also fix as usual $ \theta=\frac{1}{5} $ and $ c_0=\min\{1,c\} $. Additionally, we choose
				\begin{equation*}\label{def.zeta}
					\zeta<\min\left\{\frac{1}{4},\frac{c_0}{2}\right\}
				\end{equation*} such that 
				\begin{equation}\label{def.zeta.1}
					\frac{\artanh(4\zeta)}{4\zeta}<\frac{3}{2}\quad \text{ and }\quad\frac{c_0}{2}\zeta-\zeta^2-15(2T_M)^3
					\left(\frac{\artanh(4\zeta)}{4\zeta}-1\right)>0.
				\end{equation}
				We also consider $ \eps<\eps_2 $ defined by
				\begin{equation}\label{def.eps.2}
					\eps_2=\min\left\{1,T_M, \frac{\zeta c_0}{4}\frac{1}{4+27(2T_M)+28(2T_M)^2},\frac{2}{5}\lambda\right\}.
				\end{equation}
				We notice that $ \zeta $ depends only on $ c $, so that $ \eps_2=\eps_2(c,\lambda,T_M) $.
				Moreover, we study the family of functions for $ L>L_2(\eps,T_M) $ satisfying
				\begin{equation}\label{def.L.2}
					\left(1+\frac{2T_M}{\eps\theta}+\frac{(2T_M)^2}{(\eps\theta)^2}+\frac{(2T_M)^3}{(\eps\theta)^3}+\frac{(2T_M)^4}{(\eps\theta)^4}\right)\int_{L_2+y}^\infty E(z)dz<\int_{y-L_2}^{y+L_2}E(z)dz\text{ for all }y<L_2 .
				\end{equation}
				We remark that such $ \zeta $ as in \eqref{def.zeta.1} and such $ L_2 $ as in \eqref{def.L.2} exist, as we have seen already several times. Moreover, \eqref{def.L.2} holds true for all $ L>L_2 $.\\
				
				We will also consider $ \gamma_0=\left(2T_M-(M_L+(1-\theta)\eps)\right)e^{-\zeta R_\eps} $,
				where $ R_\eps $ is once again the distance such that $ M_L-f(y)\geq -\frac{1-\theta}{2} \eps $ for all $ y\in[L-R_\eps,L+R_\eps] $. By the uniform continuity of $ f $ and since $ f(L)\leq M_L $ we know that such $ R_\eps $ exists and it is independent of $ L $. Moreover, by definition we obtain
				\[R_{\gamma_0}=R_\eps,\]
				which implies that for all $ y\in[L,R_{\gamma_0}] $ we have
				\[\psi_{\gamma_0}^L(y)-f(y)>\frac{1-\theta}{2}\eps>0,\]
				since $ \psi_{\gamma_0}^L\geq M_L+(1-\theta)\eps $ and $ f(y)\leq M_L+\frac{1-\theta}{2}\eps  $.\\
				
				We also remark that by construction we have that for all $ \gamma\leq \gamma_0 $
				\[\psi_\gamma^L(y)>f(y)\quad \text{ for all }y\in\R\setminus(L,R_\gamma)\]
				and also $ \psi_\gamma^L(L)>M_L\geq f(y) $ as well as $ \psi_\gamma^L(R_\delta)=0<f(R_\gamma) $. Thus, $ \psi_{\gamma_0}^L>f $ in $ \R $.\\
				
				We now show that the functions $ \psi_\gamma^L $ are supersolutions to the equation \eqref{trav.wave.3.n.R} for $ y\in(L,R_\gamma) $, the interval where the functions are smooth. This will be done in the spirit of \eqref{subsol.2}, \eqref{inf.subsol.1} and \eqref{subsol.2case.1}. We use that $ -e^{\zeta (\eta-L)}\leq -e^{\zeta(R_\gamma-L)} $ for all $ \eta>R_\gamma $, we expand the power law, and we rearrange the terms. Moreover, using also $ c\geq c_0 $ we obtain
				\begin{equation}\label{supsol.1}
					-\left(\psi_\gamma^L\right)''(y)+c\left(\psi_\gamma^L\right)'(y)+\left(\psi_\gamma^L(y)\right)^4-\int_\R E(\eta-y)\left(\psi_\gamma^L(\eta)\right)^4d\eta\\\phantom{spacespacespacespacespacespacespace}
				\end{equation}
				\vspace{-0.65cm}\begin{alignat*}{3}
					\geq&\frac{c_0}{2}\zeta\gamma e^{\zeta (y-L)}+\gamma e^{\zeta (y-L)}\left(\frac{c_0}{2}\zeta-\zeta^2+4(M_L+\eps)^3-4(M_L+\eps)^3\int_{-L}^\infty E(\eta-y) e^{\zeta(\eta-y)}d\eta\right)&(I^4_1)\\
					&+4(M_L+\eps)\gamma^3 e^{3\zeta (y-L)}\left(1-\int_{-L}^\infty E(\eta-y) e^{3\zeta(\eta-y)}d\eta\right)&(I^4_2)\\
					&+6(M_L+\eps)^2\gamma^2 e^{2\zeta (y-L)}\left(1-\int_{-L}^{\infty} E(\eta-y) e^{2\zeta(\eta-y)}d\eta\right)&(I^4_3)\\
					&+\gamma^4 e^{4\zeta (y-L)}\left(1-\int_{-L}^{\infty} E(\eta-y) e^{4\zeta(\eta-y)}d\eta\right)&(I^4_4)\\
					&- 4(M_L+\eps)^3\left[\eps\theta+\int_{-\infty}^{-L}E(\eta-y)\left(2T_M-(M_L+\eps)\right)d\eta-\eps\int_{-L}^LE(\eta-y)d\eta-\eps\theta\int_{L}^\infty E(\eta-y)d\eta\right]&\quad(I^4_5)\\
					&- 4(M_L+\eps)\left[(\eps\theta)^3+\int_{-\infty}^{-L}E(\eta-y)\left(2T_M-(M_L+\eps)\right)^3d\eta-\eps^3\int_{-L}^LE(\eta-y)d\eta-(\eps\theta)^3\int_{L}^\infty E(\eta-y)d\eta\right]&(I^4_6)\\
					&+ 6(M_L+\eps)^2\left[(\eps\theta)^2-\int_{-\infty}^{-L}E(\eta-y)\left(2T_M-(M_L+\eps)\right)^2d\eta-\eps^2\int_{-L}^LE(\eta-y)d\eta-(\eps\theta)^2\int_{L}^\infty E(\eta-y)d\eta\right]&(I^4_7)\\
					&+(\eps\theta)^4-\int_{-\infty}^{-L}E(\eta-y)\left(2T_M-(M_L+\eps)\right)^4d\eta-\eps^4\int_{-L}^LE(\eta-y)d\eta-(\eps\theta)^4\int_{L}^\infty E(\eta-y)d\eta&(I^4_8)\\
					&-4\eps^3\gamma e^{\zeta (y-L)}\left(\theta^3-\int_{-L}^L E(\eta-y) e^{\zeta(\eta-y)}d\eta-\theta^3\int_{L}^\infty E(\eta-y) e^{\zeta(\eta-y)}d\eta\right)&(I^4_9)\\	
					&+4\eps \theta \gamma^3 e^{3\zeta(R_\gamma-L)}\int_{R_\gamma}^\infty E(\eta-y)d\eta&\\
					&-4\eps\gamma^3 e^{3\zeta (y-L)}\left(\theta-\int_{-L}^L E(\eta-y) e^{3\zeta(\eta-y)}d\eta-\theta\int_{L}^\infty E(\eta-y) e^{3\zeta(\eta-y)}d\eta\right)&(I^4_{10})\\
					&+4^3 \gamma^3 e^{\zeta(R_\gamma-L)}\int_{R_\gamma}^\infty E(\eta-y)d\eta&\\\end{alignat*}\begin{alignat*}{3}
					&+6\eps^2\gamma^2 e^{2\zeta (y-L)}\left(\theta^2-\int_{-L}^L E(\eta-y) e^{2\zeta(\eta-y)}d\eta-\theta^2\int_{L}^{R_\gamma} E(\eta-y) e^{2\zeta(\eta-y)}d\eta\right)\quad\quad\quad\quad\quad\quad\quad\quad\quad\quad\quad\quad&(I^4_{11})\\&-6(\eps\theta\gamma)^2e^{2\zeta(R_\gamma-L)}\int_{R_\gamma}^\infty E(\eta-y)d\eta&\\
					&-12(M_L+\eps)^2\eps\gamma e^{\zeta (y-L)}\left(\theta-\int_{-L}^L E(\eta-y) e^{\zeta(\eta-y)}d\eta-\theta\int_{L}^{R_\gamma} E(\eta-y) e^{2\zeta(\eta-y)}d\eta\right)&(I^4_{12})\\&+12(M_L+\eps)^2\eps\theta\gamma e^{\zeta(R_\gamma-L)}\int_{R_\gamma}^\infty E(\eta-y)d\eta&\\
					&-12(M_L+\eps)\eps\gamma^2 e^{2\zeta (y-L)}\left(\theta-\int_{-L}^L E(\eta-y) e^{2\zeta(\eta-y)}d\eta-\theta\int_{L}^{R_\gamma} E(\eta-y) e^{2\zeta(\eta-y)}d\eta\right)&(I^4_{13})\\&+12(M_L+\eps)\eps\theta\gamma^2e^{2\zeta(R_\gamma-L)}\int_{R_\gamma}^\infty E(\eta-y)d\eta&\\
					&+12(M_L+\eps)\eps^2\gamma e^{\zeta (y-L)}\left(\theta^2-\int_{-L}^L E(\eta-y) e^{\zeta(\eta-y)}d\eta-\theta^2\int_{L}^{R_\gamma} E(\eta-y) e^{\zeta(\eta-y)}d\eta\right)\quad\quad\quad\quad\quad\quad\quad\quad\quad&(I^4_{14})\\&-12(M_L+\eps)\eps^2\theta^2\gamma2e^{\zeta(R_\gamma-L)}\int_{R_\gamma}^\infty E(\eta-y)d\eta.&
				\end{alignat*}
				We now proceed to estimate all the terms. First of all, using the identity \eqref{artanh}, the estimate $ (M_L+\eps)\leq 2T_M $ as well as the definition of $ \zeta $ in \eqref{def.zeta.1} we compute
				\begin{equation}\label{supsol.2}
					\begin{split}
						(I^4_1)+(I^4_2)+(I^4_3)+(I^4_4)\geq &\frac{c_0}{2}\zeta\gamma e^{\zeta (y-L)}+\gamma e^{\zeta (y-L)}\left(\frac{c_0}{2}\zeta-\zeta^2-4(2T_M)^3\left(\frac{\artanh(\zeta)}{\zeta}-1\right)\right)\\
						&-4(2T_M)\gamma^3 e^{3\zeta (y-L)}\left(\frac{\artanh(3\zeta)}{3\zeta}-1\right)\\
						&-6(2T_M)^2\gamma^2 e^{2\zeta (y-L)}\left(\frac{\artanh(2\zeta)}{2\zeta}-1\right)+\gamma^4 e^{4\zeta (y-L)}\left(\frac{\artanh(4\zeta)}{4\zeta}-1\right)\\
						\geq \frac{c_0}{2}&\zeta\gamma e^{\zeta (y-L)}+\gamma e^{\zeta (y-L)}\left(\frac{c_0}{2}\zeta-\zeta^2-4(2T_M)^3\left(\frac{\artanh(4\zeta)}{4\zeta}-1\right)\right)\\
						-&4(2T_M)^3\gamma e^{\zeta (y-L)}\left(\frac{\artanh(4\zeta)}{4\zeta}-1\right)\\
						-&6(2T_M)^3\gamma e^{\zeta (y-L)}\left(\frac{\artanh(4\zeta)}{4\zeta}-1\right)+(2T_M)^3\gamma e^{\zeta (y-L)}\left(\frac{\artanh(4\zeta)}{4\zeta}-1\right)\\
						=\frac{c_0}{2}\zeta\gamma e^{\zeta (y-L)}&+\gamma e^{\zeta (y-L)}\left(\frac{c_0}{2}\zeta-\zeta^2-15(2T_M)^3\left(\frac{\artanh(4\zeta)}{4\zeta}-1\right)\right)\geq\frac{c_0}{2}\zeta\gamma e^{\zeta (y-L)},
					\end{split}
				\end{equation}
				where we used also that $ \zeta \mapsto \frac{\artanh(\zeta)}{\zeta}-1 $ is a monotonically increasing non-negative function and that $ \gamma e^{\zeta(y-L)}\leq 2T_M $ for $ y\leq R_\gamma $.
				
				We estimate the terms $ (I^4_5),\;(I^4_6),\;(I^4_7) $ and $ (I^4_8) $.
				We compute using $ \theta=\frac{1}{5} $ and the choice of $ L>L_2 $ as in \eqref{def.L.2}
				\begin{equation}\label{supsol.3}
					\begin{split}
						(I^4_5)+(I^4_7)\geq &- 4(M_L+\eps)^3\left[\eps\theta\int_{-\infty}^{-L} E(\eta-y)d\eta+\int_{-\infty}^{-L}E(\eta-y)\left(2T_M-(M_L+\eps)\right)d\eta-4\theta\eps\int_{-L}^LE(\eta-y)d\eta\right]\\
						+& 6(M_L+\eps)^2\left[(\eps\theta)^2\int_{-\infty}^{L} E(\eta-y)d\eta-\int_{-\infty}^{-L}E(\eta-y)\left(2T_M-(M_L+\eps)\right)^2d\eta-\eps^2\int_{-L}^LE(\eta-y)d\eta\right]\\
						\geq & -4(M_L+\eps)^3\eps\theta\left[\left(1+\frac{2T_M}{\eps\theta}\right)\int_{-\infty}^{-L} E(\eta-y)d\eta-4\int_{-L}^LE(\eta-y)d\eta\right]\\
						&+ 6(M_L+\eps)^2\left[(\eps\theta)^2\int_{-\infty}^{L} E(\eta-y)d\eta-\left(2T_M\right)^2\int_{-\infty}^{-L}E(\eta-y)d\eta-\eps^2\int_{-L}^LE(\eta-y)d\eta\right]\\\end{split}\end{equation}\begin{equation*}\begin{split}
						\geq & 4(M_L+\eps)^3\left[3\theta\eps\int_{-L}^LE(\eta-y)d\eta\right]- 6(M_L+\eps)^2\left[\eps^2\int_{-L}^LE(\eta-y)d\eta\right]\\	
						=&6(M_L+\eps)^2\eps\int_{-L}^LE(\eta-y)d\eta\left(\frac{2}{5}(M_L+\eps)-\eps\right)>0
					\end{split}
				\end{equation*}
				for $ \eps<\eps_2 $ as in \eqref{def.eps.2} since $ M_L+\eps>\lambda $.
				Since $ \frac{492}{125}>\frac{2}{5} $, in a similar way we can estimate
				\begin{equation*}
					\begin{split}
						(I^4_6)+(I^4_8)\geq &- 4(M_L+\eps)\left[(\eps\theta)^3\int_{-\infty}^{-L} E(\eta-y)d\eta+\int_{-\infty}^{-L}E(\eta-y)\left(2T_M\right)^3d\eta-124(\theta\eps)^3\int_{-L}^LE(\eta-y)d\eta\right]\\
						&+ (\eps\theta)^4\int_{-\infty}^{L} E(\eta-y)d\eta-\int_{-\infty}^{-L}E(\eta-y)\left(2T_M\right)^4d\eta-\eps^4\int_{-L}^LE(\eta-y)d\eta\\	
						\geq & 4(M_L+\eps)(\eps\theta)^3\left[123(\theta\eps)^3\int_{-L}^LE(\eta-y)d\eta\right]- \eps^4\int_{-L}^LE(\eta-y)d\eta,\\\end{split}
				\end{equation*} so that\begin{equation}\label{supsol.4}
					\begin{split}
						(I^4_6)+(I^4_8)\geq&\eps^3\int_{-L}^LE(\eta-y)d\eta\left(\frac{492}{125}(M_L+\eps)-\eps\right)>0.
					\end{split}
				\end{equation}
				
				We estimate the last terms using $ \frac{\artanh(4\zeta)}{4\zeta}<\frac{3}{2} $ as well as $ \gamma e^{\zeta(y-L)}\leq 2T_M $.
				\begin{equation}\label{supsol.5}
					(I^4_9)\geq -4\eps^3\theta^3\gamma e^{\zeta(y-L)}\geq -4\eps\gamma e^{\zeta(y-L)},\quad\quad 	(I^4_{10})\geq -4\eps\theta (2T_M)^2\gamma e^{\zeta(y-L)},
				\end{equation}
				%
				\begin{equation}\label{supsol.7}
					(I^4_{11})\geq -6\eps^2\gamma^2 e^{2\zeta(y-L)}\frac{\artanh(2\zeta)}{2\zeta}\geq -9\eps (2T_M)\gamma e^{\zeta(y-L)}
				\end{equation}
				and similarly
				\begin{equation}\label{supsol.7.1}
					(I^4_{12})+	(I^4_{13})+	(I^4_{14})\geq -12(2T_M)\eps\gamma e^{\zeta(y-L)}\left(2(2T_M)+18\right)
				\end{equation}
				Finally, using $ \eps<\eps_2 $ as given in \eqref{def.eps.2} and combining the equations \eqref{supsol.1}-\eqref{supsol.7.1} we conclude that $ \psi_\gamma^L $ are supersolutions in $ (L,R_\gamma) $, i.e.
				\begin{equation*}\label{supsol.8}
					-\left(\psi_\gamma^L\right)''(y)+c\left(\psi_\gamma^L\right)'(y)+\left(\psi_\gamma^L(y)\right)^4-\int_\R E(\eta-y)\left(\psi_\gamma^L(\eta)\right)^4d\eta\\\geq \frac{c_0}{4}\zeta\gamma e^{\zeta (y-L)}>0
				\end{equation*}
				for all $ y\in(L,R_\gamma) $.\\
				
				We recall that by construction $ \psi_{\gamma_0}^L>f $ in $ \R $ with $ \psi_{\gamma_0}^L-f\geq\frac{1-\theta}{2}\eps>0 $ for all $ y\geq L $ and $ \psi_{\gamma}^L\geq f $ for $ y\in\R\setminus(L,R_\gamma) $ with $ \psi_\gamma^L\left.\right|_{\{L,R_\gamma\}}>f\left.\right|_{\{L,R_\gamma\}} $. Hence, once again arguing with the maximum principle as we did in the proof of Lemma \ref{lemma.no.sup.infty} and of Lemma \ref{lemma.no.inf.infty} we conclude, by the uniform continuity of $ \gamma\mapsto \psi_\gamma^L $ on compact sets and their decreasing monotonicity, that 
				\[\psi_\gamma^L(y)\geq f(y) \quad \text{ for all }y\in\R \text{ and }\gamma<\gamma_0,\] since $ \psi_\gamma^L $ are supersolutions on $ (L,R_\gamma) $. 
				Hence, for any $ y>L $ we have for $ \gamma<\gamma_0 $ small enough
				\[f(y)\leq M_L+(1-\theta)\eps +\gamma e^{\zeta(y-L)}.\]
				Finally, taking $ \gamma\to 0 $ and thus $ R_\gamma\to\infty $ we conclude
				\begin{equation}\label{maximum}
					f(y)\leq M_L+(1-\theta)\eps\quad\text{ for all }y>L.
				\end{equation}
				Let us now define $ \eps_0(T_M,\lambda,c)=\min\{\eps_1,\eps_2\} $ for $ \eps_1 $ and $ \eps_2 $ defined in \eqref{def.eps1} and \eqref{def.eps.2}, respectively. For any given $ \eps<\eps_0 $ we define also $ L_0(\eps,T_M,\lambda,c)=\max\{L_1,L_2\} $, where $ L_1,\;L_2 $ are given in  \eqref{def.L.1} and \eqref{def.L.2}, respectively, and $ \theta=\frac{1}{5} $.\\
				
				The estimates \eqref{minimum} and \eqref{maximum} yield the proof of Theorem \ref{thm.stability}. Indeed, we have just proved that, if $ f $ solves \eqref{trav.wave.3.n.R}, there exists some $ \eps_0(T_M,\lambda,c)>0 $ such that, if 
				\[\osc\limits_{[-L,L]}f<\eps\]
				for $ \eps<\eps_0 $ and for $ L>L_0(\eps,T_M,\lambda,c) $, then 
				\begin{equation*}\label{osc.final}
					\osc\limits_{[L,\infty)}f\leq  M_L+(1-\theta)\eps-m_L+(1-\theta)\eps=(3-2\theta)\eps< 3\eps,
				\end{equation*}
				where we also use that $ m_L\leq f(L)\leq M_L $.
			\end{proof}
		\end{theorem}
		In order to use Theorem \ref{thm.stability} we need to have functions satisfying the oscillation assumption. The next lemma shows that there exist sequences of functions satisfying both \eqref{trav.wave.3.n.R} and the oscillation condition.
		\begin{lemma}\label{lemma.osc}
			Let $ f $ solve \eqref{trav.wave.3.n.R} for $ 0<\lambda<T_M $. Let us assume that $ f $ is not constant. Then there exist $ \{x_n\}_{n\in\mathbb{N}} $ and $ \{\xi_k\}_{k\in\mathbb{N}} $ monotonically decreasing sequences with $ \lim\limits_{n\to\infty}x_n= -\infty $ as well as $  \lim\limits_{k\to\infty}\xi_k= -\infty $ satisfying \[ \lim\limits_{n\to\infty}f(x_n)= \sup\limits_\R f\quad\text{ and }\quad \lim\limits_{k\to\infty}f(\xi_k)= \inf\limits_\R f.\]
			Moreover, they satisfy
			\[\lim\limits_{n\to\infty}\osc\limits_{[-L,L]} f(x_n+\cdot)=0\quad\text{ and }\quad \lim\limits_{k\to\infty}\osc\limits_{[-L,L]} f(\xi_k+\cdot)=0\]
			for all $ L>0 $.
			\begin{proof}
				Since $ f $ is not constant, according to Lemma \ref{lemma.no.sup.interior}, Lemma \ref{lemma.no.sup.infty} and Lemma \ref{lemma.no.inf.infty} it has to attain its supremum and infimum at $ -\infty $. Hence, there exist monotonically decreasing sequences $ \{x_n\}_{n\in\mathbb{N}} $ and $ \{\xi_k\}_{k\in\mathbb{N}} $ satisfying
				\[\lim\limits_{n\to\infty}x_n= -\infty ,\quad\lim\limits_{k\to\infty}\xi_k= -\infty ,\quad\lim\limits_{n\to\infty}f(x_n)= \sup\limits_\R f\quad\text{ and }\quad \lim\limits_{k\to\infty}f(\xi_k)= \inf\limits_\R f.\] We now prove the statement for the supremum. We define $ f_n=f(x_n+\cdot) $. Then $ f_n $ solves the same equation \eqref{trav.wave.3.n.R} by the translation invariance of this equation. Moreover, $ f_n\in C^{2,1/2}(\R)  $. Thus, by compactness for $ \alpha\in\left(0,\frac{1}{2}\right) $ there exists a subsequence $ f_{n_j}=f\left(x_{n_j}+\cdot\right)\to g $ in $ C^{2,\alpha}(\R) $ in every compact set and hence uniformly everywhere.
				Moreover, $ g$ solves also \eqref{trav.wave.3.n.R}. By regularity theory we see that  $ g\in  C^{2,1/2}(\R)  $. 
				
				It is important to notice also that 
				\[g(0)=\lim\limits_{j\to\infty}f\left(x_{n_j}\right)=\sup\limits_{\R} f\geq \sup\limits_{\R} g\geq g(0).\]Since $ g $ attains its supremum at the interior, it is constant according to Lemma \ref{lemma.no.sup.interior}. Thus, $ f_{n_j}\to g=\sup\limits_{\R} f  $ uniformly in every compact set. 
				
				Let $ \eps>0 $ and $ L>0 $. By the uniform convergence in $ [-L,L] $ there exists $ N_0(\eps,L)>0 $ such that for all $ j\geq N_0 $ we have
				\[\Arrowvert (f_{n_j}-g)\left.\right|_{[-L,L]}\Arrowvert_\infty<\frac{\eps}{2}.\]
				We thus conclude that
				\[\osc\limits_{[-L,L]}f_{n_j}=\max\limits_{[-L,L]}f_{n_j}-\min\limits_{[-L,L]}f_{n_j}<\frac{\eps}{2}+g-g+\frac{\eps}{2}=\eps.\]
				This proves $ \lim\limits_{j\to\infty}\osc\limits_{[-L,L]} f(x_{n_j}+\cdot)=0 $ for all $ L>0 $. Thus, the sequence $ \{\tilde{x}_j\}_{j\in\mathbb{N}}=\{x_{n_j}\}_{j\in\mathbb{N}} $ satisfies the statement of Lemma \ref{lemma.osc} concerning the supremum. \\
				
				Using that any solution to \eqref{trav.wave.3.n.R} which attains its infimum at the interior is constant according to Lemma \ref{lemma.no.sup.interior}, we conclude the proof of this lemma repeating the same arguments for the sequence $ f_k= f(\xi_k+\cdot) $, for which a subsequence converges uniformly in every compact set to $ g=\inf\limits_{\R} f $.
			\end{proof}
		\end{lemma}
		Finally, Lemma \ref{lemma.osc} and Theorem \ref{thm.stability} together with the previous results in Lemma \ref{lemma.no.sup.interior}, Lemma \ref{lemma.no.sup.infty} and Lemma \ref{lemma.no.inf.infty} imply that the solution $ f $ to \eqref{trav.wave.3.n.R} is constant.
		\begin{theorem}\label{thm.constant}
			Let $ f $ solve \eqref{trav.wave.3.n.R} for $ 0<\lambda<T_M $. Then $ f $ is constant.
			\begin{proof}
				Let $ f $ solve \eqref{trav.wave.3.n.R}. Let us assume that $ f $ is not constant. By Lemma \ref{lemma.no.sup.interior}, Lemma \ref{lemma.no.sup.infty} and Lemma \ref{lemma.no.inf.infty} there exist $ \{x_n\}_{n\in\mathbb{N}} $ and $ \{\xi_k\}_{k\in\mathbb{N}} $ monotonically decreasing sequences with $ \lim\limits_{n\to\infty}x_n= -\infty $ as well as $  \lim\limits_{k\to\infty}\xi_k= -\infty $ satisfying \[ \lim\limits_{n\to\infty}f(x_n)= \sup\limits_\R f\quad\text{ and }\quad \lim\limits_{k\to\infty}f(\xi_k)= \inf\limits_\R f.\]
				Let also $ \eps<\eps_0 $ be arbitrary and $ L>L_0(\eps) $ for $ \eps_0 $ and $ L_0(\eps) $ as in Theorem \ref{thm.stability}. According to Lemma \ref{lemma.osc} there exists $ N_0(\eps,L) $ such that 
				\[\osc\limits_{[-L,L]} f(x_n+\cdot)<\eps\quad \text{ for all }n\geq N_0.\]
				Since by the translation invariance $ f(x_n+\cdot) $ solve \eqref{trav.wave.3.n.R} with $ \lambda\leq f(x_n+\cdot)\leq T_M $, Theorem \ref{thm.stability} implies that
				\[\osc\limits_{[L,\infty)} f(x_n+\cdot)<3\eps\quad \text{ for all }n\geq N_0.\]
				Thus,
				\[\osc\limits_{[-L,\infty)} f(x_n+\cdot)<4\eps\quad \text{ for all }n\geq N_0.\]
				Similarly, there exists $ K_0(\eps,L)>0 $ such that 
				\[\osc\limits_{[-L,\infty)} f(\xi_k+\cdot)<4\eps\quad \text{ for all }k\geq K_0.\]
				Hence, for any $ n\geq N_0 $ and $ k\geq K_0 $ it is either $ x_n-\xi_k>0 $ or $ \xi_k-x_n>0 $.
				In the first case we estimate $ \left|f(\xi_k)-f(x_n)\right|\leq \osc\limits_{[-L,\infty)}f(\xi_k+\cdot)<4\eps $, while in the latter situation $ \left|f(x_n)-f(\xi_k)\right|\leq \osc\limits_{[-L,\infty)}f(x_n+\cdot)<4\eps $.
				Therefore, \[\left|f(x_j)-f(\xi_j)\right|<4\eps\quad \text{ for all } j\geq \max\{N_0,K_0\}.\]
				Taking now the limit as $ j\to\infty $ we conclude
				\[\sup\limits_\R f-\inf\limits_\R f\leq4\eps.\]
				Since $ \eps<\eps_0 $ was arbitrary, this implies that $\sup\limits_\R f=\inf\limits_\R f $ and hence $ f $ is constant.
			\end{proof}
		\end{theorem}
		\subsection{Existence of a positive limit of the traveling waves as $ y\to\infty $}
		We now finish this section proving that any traveling wave solving \eqref{trav.wave.3.1} for $ y>0 $ has a limit as $ y\to\infty $. We first of all need to show a corollary to the stability result in Theorem \ref{thm.stability}.
		\begin{corollary}\label{cor.stability}
			Let $ f $ solve \eqref{trav.wave.3.1} according to Theorem \ref{thm.existence.all} for $ 0<\lambda\leq f\leq T_M $ and $ c>0 $. Let $ \eps<\eps_0(c,\lambda,T_M) $ and $ L_0(\eps,\lambda,T_M,c) $ be as in Theorem \ref{thm.stability}. Let also $ a>L_0(\eps) $ and $ \tilde{f}(y):=f(a+y) $. Then $ \tilde{f}:[-a,\infty)\to \R_+ $ solves
			\begin{equation}\label{trav.wave.3.a}
				-\tilde{f}''(y)+c\tilde{f}'(y)+\tilde{f}^4(y)-\int_{-a}^\infty E(\eta-y)\tilde{f}^4(\eta)d\eta=0
			\end{equation}
			with $ \tilde{f}(-a)=T_M $ and $ 0<\lambda\leq \tilde{f}\leq T_M $. Moreover, if \[\osc\limits_{[-L,L]}\tilde{f}<\eps\quad \text{ for }L_0(\eps)<L<a\]
			then
			\[\osc\limits_{[L,\infty)}\tilde{f}<3\eps.\]
			\begin{proof}
				It is easy to see that $\tilde{f}$ solves \eqref{trav.wave.3.a}. In order to simplify the reading we use the same notation as in Theorem \ref{thm.stability}. Let hence $ m_L $ and $ M_L $ being the minimum and the maximum of $ \tilde{f} $ on $ [-L,L] $, respectively. Moreover, $ \beta,\;\zeta>0 $, $ \delta<\delta_0 $ and $ \gamma<\gamma_0 $ are defined for $ \tilde{f} $ as in Theorem \ref{thm.stability}. Let finally $ \psi_\delta^L $ as in \eqref{def.psi.delta.1} and $ \psi_\gamma^L $ as in \eqref{def.psi.gamma}. We argue that $ \psi_\delta^L\mathds{1}_{[-a,\infty)} $ and $ \psi_\gamma^L\mathds{1}_{[-a,\infty)} $ are subsolutions and supersolutions for the equation \eqref{trav.wave.3.a} on $ (L,R_\delta) $ and $ (L,R_\gamma) $, respectively. \\
				
				Indeed, by definition $ \psi_\delta^L\mathds{1}_{[-a,\infty)}=\psi_\delta^L $ since $ L_0(\eps)<L<a $ and $ \psi_\delta^L=0 $ for $ y<-L $. This implies that 
				\[\int_{-a}^\infty E(\eta-y) \left(\psi_\delta^L(\eta)\right)^4d\eta=\int_{-L}^\infty E(\eta-y) \left(\psi_\delta^L(\eta)\right)^4d\eta=\int_{-\infty}^\infty E(\eta-y) \left(\psi_\delta^L(\eta)\right)^4d\eta.\]
				Hence, for $ y\in(L,R_\delta) $ the functions $ \psi_\delta^L\mathds{1}_{[-a,\infty)} $ are subsolutions for the equation \eqref{trav.wave.3.a} with $ \psi_\delta^L\mathds{1}_{[-a,\infty)}\leq \tilde{f} $ on $ [-a,\infty)\setminus(L,R_\delta) $, $ \psi_\delta^L\left|_{\{L,R_\delta\}}\right.<\tilde{f}\left|_{\{L,R_\delta\}}\right. $, as well as $ \psi_{\delta_0}^L\mathds{1}_{[-a,\infty)}< \tilde{f} $ on $ [-a,\infty) $.\\
				
				Similarly, $ \psi_\gamma^L\mathds{1}_{[-a,\infty)}\leq\psi_\gamma^L $ since $ L_0(\eps)<L<a $ and $ \psi_\gamma^L=2T_M $ for $ y<-L $. This implies that 
				\[\int_{-a}^\infty E(\eta-y) \left(\psi_\gamma^L(\eta)\right)^4d\eta\leq\int_{-\infty}^\infty E(\eta-y) \left(\psi_\gamma^L(\eta)\right)^4d\eta.\]
				Thus, the functions $ \psi_\gamma^L\mathds{1}_{[-a,\infty)} $ are supersolutions for the equation \eqref{trav.wave.3.a} for $ y\in(L,R_\gamma) $. Moreover, they satisfy  $ \psi_\gamma^L\mathds{1}_{[-a,\infty)}\geq \tilde{f} $ on $ [-a\infty)\setminus(L,R_\gamma) $, $ \psi_\gamma^L\left|_{\{L,R_\gamma\}}\right.>\tilde{f}\left|_{\{L,R_\gamma\}}\right. $, as well as $ \psi_{\gamma}^L\mathds{1}_{[-a,\infty)}> \tilde{f} $ on $ [-a,\infty) $.\\
				
				Hence, as we saw in Theorem \ref{thm.stability}, an application of the maximum principle and of the uniform continuity on compact sets of the families of sub- and supersolutions with respect of $ \delta $ and $ \gamma $, respectively, implies 
				\[\osc\limits_{[L,\infty)} \tilde{f}<3\eps.\]
			\end{proof}
		\end{corollary}
		Finally, we can prove the convergence of the traveling wave to a positive constant as $ y\to\infty $.
		\begin{theorem}\label{thm.convergence}
			Let $ f $ solve \eqref{trav.wave.3.1} according to Theorem \ref{thm.existence.all} for $ T_M>0 $ and $ c>0 $. Then there exists a limit
			\[\lim\limits_{y\to\infty}f(y)=:f_\infty>0.\]
			\begin{proof}
				By Theorem \ref{thm.existence.all}, Lemma \ref{lemma.monotonicity} and Theorem \ref{thm.small.melting} we know that $ f\geq\lambda>0 $ for some $ \lambda>0 $. Let us take $\{x_n\}_{n\in\mathbb{N}} $ and $ \{\xi_n\}_{n\in\mathbb{N}}$ two diverging monotone increasing sequences such that
				\[\lim\limits_{n\to\infty}f(x_n)=\limsup\limits_{y\to\infty}f(y)=:\overline{f_\infty}\quad\text{ and }\quad\lim\limits_{n\to\infty}f(\xi_n)=\liminf\limits_{y\to\infty}f(y)=:\underline{f_\infty} .\]
				We notice that $ \overline{f_\infty},\underline{f_\infty} \in[\lambda,T_M] $.
				Up to subsequences we know that $ f(x_n+\cdot) $ and $ f(\xi_n+\cdot) $ converge to constant functions, as we have proved in Theorem \ref{thm.constant}. We denote these subsequences $ x_n $ and $ \xi_n $. Hence, we have
				\[\lim\limits_{n\to\infty}f(x_n+\cdot)=\overline{f_\infty} \text{ and }\lim\limits_{n\to\infty}f(\xi_n+\cdot)=\underline{f_\infty} \]
				uniformly on compact sets.
				Therefore, for all $ L>0 $ there exists $ N_0(L) $ such that $ x_n,\xi_n>L $ for all $ n\geq N_0(L) $ and such that
				\[\osc\limits_{[-L,L]} f(x_n+\cdot)\to 0 \quad\text{ and }\quad\osc\limits_{[-L,L]} f(\xi_n+\cdot)\to 0 \quad\text{ as }n\to\infty \text{ and }n\geq N_0(L).\]
				Let now $ \eps<\eps_0(c,\lambda,T_M) $ and $ L_0(\eps,c,\lambda,T_M) $ as defined in Theorem \ref{thm.stability} and in Corollary \ref{cor.stability}. Then there exists $ N_1(\eps,L_0(\eps))>0 $ such that $ x_n,\xi_n>L_0(\eps) $ for all $ n\geq N_1 $. Let also $ L\in(L_0(\eps),\min\{x_{N_1},\xi_{N_1}\}) $. Then there exists $ N_2(\eps,L)\geq N_1(\eps) $ such that 
				\[\osc\limits_{[-L,L]} f(x_n+\cdot)<\eps \quad\text{ and }\quad\osc\limits_{[-L,L]} f(\xi_n+\cdot)<\eps \quad\text{ for all }n\geq N_2(\eps,L).\]
				We remark that $ L_0(\eps)<L<\min\{x_n,\xi_n\} $ for all $ n\geq N_2(\eps,L) $. Then by the Corollary \ref{cor.stability} we can conclude that 
				\[\osc\limits_{[L,\infty)} f(x_n+\cdot)<3\eps \quad\text{ and }\quad\osc\limits_{[L,\infty)} f(\xi_n+\cdot)<3\eps \quad\text{ for all }n\geq N_2(\eps,L).\]
				This implies that 
				\begin{equation}\label{convergence1}
					\left|f(x_n)-f(\xi_n)\right|<4\eps \text{ for all }n\geq N_2(\eps,L).
				\end{equation}
				Indeed, let $ n\geq N_2(\eps,L) $. If $ x_n-\xi_n>0 $ we compute
				\[\left|f(\xi_n)-f(x_n)\right|\leq \osc\limits_{[-L,\infty)} f(\xi_n+\cdot)<4\eps,\]
				while if $ \xi_n-x_n>0 $
				\[\left|f(x_n)-f(\xi_n)\right|\leq \osc\limits_{[-L,\infty)} f(x_n+\cdot)<4\eps.\]
				Taking the limit $ n\to\infty $ in \eqref{convergence1} we obtain
				\[0\leq \overline{f_\infty}-\underline{f_\infty}\leq4\eps,\]
				which implies that $ f $ has a limit, since $ \eps<\eps_0 $ is arbitrarily small, i.e.
				\[\overline{f_\infty}=\underline{f_\infty}=f_\infty.\]
			\end{proof}
		\end{theorem}
		\section{Formal description of the long time asymptotic for arbitrary values of $ T(\pm \infty) $}\label{Sec.pict.}
		In this last section we conclude giving the expected behavior of the solution to the Stefan problem \eqref{syst.3} as $ t\to\infty $. We remark that what we present here is formal.
		
		Theorem \ref{thm.trav.wave.0} shows the existence of $ c_{\max}>0 $ such that for any $ c\in(0,c_{\max}) $ there exists traveling waves $ T_1(x+ct)=:T_1^c(y) $ and $ T_2(x+ct)=:T_2^c(y) $ solving the Stefan problem for $ s(t)=-ct $. The first problem we should solve concerns the uniqueness of the traveling waves.
		\begin{definition}
			Prove or disprove that for any $ c\in(0,c_{\max}) $ and $ T>0 $ the traveling waves $ T^c_1,T^c_2 $ solving \eqref{trav.wave.1} are unique.
		\end{definition}
		Notice that it is enough to have the uniqueness of the traveling wave in the solid solving \eqref{trav.wave.3}.
		
		Recall that $ \lim\limits_{y\to\infty}T_2^c>0 $ and $ \lim\limits_{y\to-\infty}T_1^c=T_M-\frac{cL+\partial_yT_2^c(0+)}{LK} $. Moreover, we notice that also for $ c=c_{\max} $ there exist traveling waves. Indeed, $ T_2^{c_{\max}} $ exists by Theorem \ref{thm.existence.all}. By definition $ \partial_yT_2^{c_{\max}}(0^+)=-Lc_{\max} $. Thus, since $ \partial_y T_1^{c_{\max}}(0^-)=0 $, in this case the traveling wave is constant in the liquid part, i.e. $ T_1^{c_{\max}}=T_M $.
		
		Also the existence of a traveling wave $ T_2^0 $ solving \eqref{trav.wave.3} for $ c=0 $ is an important problem that should be considered. 
		\begin{definition}
			Prove or disprove that there exists a unique traveling wave $ T_2^0 $ solving \eqref{trav.wave.3} for $ c=0 $. Moreover, $ T_2^0 $ converges to a positive constant as $ y\to\infty $.	
		\end{definition}
		\begin{remark}
			The existence of $ T_2^0 $ can be proved as follows using an iterative argument. First of all the function \[f_1(y)=\frac{A}{\left(B+y\right)^{\frac{2}{3}}} \text{, where }A=\frac{2}{\sqrt[3]{9}}\text{ and }B=\frac{1}{3}\left(\frac{2}{T_M}\right)^{\frac{2}{3}},\]
			is a solution to $ f_1''-f_1^4=0 $ on $ \R_+ $ with $ f_1(0)=T_M $. Moreover, $ f_1 $ is monotonically decreasing with $ \lim\limits_{y\to\infty}f_1(y)=0 $. It is also possible to show the existence of a monotone sequence $ 0\leq f_1\leq f_2\leq ...\leq f_n\leq f_{n-1}\leq ...\leq T_M $ solving for $ n\geq 2 $ equation \eqref{trav.wave.3.1} for $ c=0 $. In this case though, the variational principle method we used in Proposition \ref{prop.variational} does not work. Nevertheless, knowing for $ n\geq 2 $ the existence of $ f_{n-1}\in C^{0,1/2}(R_+) $ with $ f_1\leq f_{n-1}\leq T_M $, the method of sub- and supersolutions (c.f \cite{evans}) can be implemented in order to find for any $ R>0 $ a solution $ f_n^R\in C^{2,1/2}([0,R]) $ of the boundary value problem
			\[\begin{cases}
				-\left(f_n^R\right)''(y)+\left(f_n^R(y)\right)^4=\int_0^\infty E(y-\eta)f_{n-1}^4(\eta)d\eta&\quad y\in(0,R)\\
				f_n^R(0)=T_M\\
				f_n^R(R)=f_1(R).
			\end{cases}\]
			Indeed, $ f_1 $ and $ T_M $ are sub- and supersolutions of the operator $ L(u)=-\partial_y^2u+4T_M^3u $ and the function $ \lambda\mapsto -\lambda^4+4T_M^3\lambda $ is increasing for $ \lambda\in[0,T_M] $. Moreover, since $ \Arrowvert f_n^R\Arrowvert_\infty\leq T_M $ as well as $ \Arrowvert \partial_y^2f_n^R\Arrowvert_\infty\leq T_M^4 $ we conclude that $ f_n^R\in C^{2,1/2}([0,R]) $ with uniformly bounded norm with respect to $ R $. Hence, taking the limit we prove the existence of a function $ f_n\in C^{2,1/2}(\R_+) $ solution to \eqref{trav.wave.3.1} for $ c=0 $. Since the monotonicity argument in Theorem \ref{thm.existence.all} applies also in this case, such a monotone sequence exists. This implies the existence of a traveling wave solving \eqref{trav.wave.3} for $ c=0 $ and $ y>0 $. However, the uniqueness and the existence of a positive limit are more involved problems. 
		\end{remark} 
		This remark shows that $ T_2^0 $ exists, moreover, $ \partial_y T_2^0(0^+)<0 $ by the Hopf-principle. Hence, in the liquid the traveling wave $ T_1^0 $ solves $ \partial_y^2T_1^0=0 $ with $ T_1^0(0)=T_M $ and $ \partial_yT_1^0(0^-)=\frac{\partial_yT_2^0(0^+)}{K} $. Thus, we obtain \[\lim\limits_{c\to 0}T_1^c(y)=T_M-\frac{\partial_yT_2^0(0^+)}{K}y\text{ with }\lim\limits_{y\to-\infty}T_1^{0}(y)=\infty.\]
		These observations lead to the following open problem.
		\begin{definition}
			Prove or disprove that for any $ T_{-\infty}\in[T_M,\infty] $ there exists a unique $ c\in\left[0,c_{\max}\right] $ such that in the liquid the traveling wave $ T_1^c $ of Theorem \ref{thm.trav.wave.all} satisfies\[\lim\limits_{y\to-\infty}T_1^c(y)=T_{-\infty}.\]
		\end{definition}
		On the contrary, in the solid we already know that there exists $ \theta>0 $ such that for any $ c\geq 0 $ the traveling waves satisfy $ \lim\limits_{y\to\infty}T_2^c=T_{\text{int}}^c\geq \theta $. Therefore, we cannot expect that $ T_{\text{int}}^c $ can attain all the values in $ [0,T_M] $ for $ c\in\left[0,c_{\max}\right] $. Nevertheless, we can reach any value in $ [0,T_M] $ if we include an additional layer in which the radiative transfer equation is approximated using the diffusion approximation. More precisely, we expect to approximate the evolution equation of the temperature by an equation of the form
		\[ T_t=T_{xx}+\left(T^4\right)_{xx}\] in a domain $ x>s(t) $ where $ T $ changes in a length scale much larger than $ 1 $.
		
		We conclude the final picture of the asymptotic of the solution $ (T_1,T_2,s) $ to the Stefan problem \eqref{syst.3} as $ t\to\infty $ with the following claim.\\
		
		Given $ T_{-\infty}\in[T_M,\infty] $ and $ T_\infty\in[0,T_M] $ there exist $ c\in\left[0,c_{\max}\right] $ and functions $ T_1^c,T_2 $ with the following properties: \begin{enumerate}
			\item[(i)] $ s(t)=-ct $;
			\item[(ii)] $ T_1^c $ is the traveling wave of Theorem \eqref{thm.existence.all} for $ y<0 $ with $ \lim\limits_{y\to-\infty}T_1^c(y)=T_{-\infty} $;
			\item[(iii)] $ T_2 $ is given by the traveling wave $ T_2^c $ of Theorem \eqref{thm.existence.all} for $ y>0 $ and by a self-similar profile $ F $ connecting $ T_{\text{int}}^c $ to $ T_\infty $, which solves
			\begin{equation}\label{selfsimilar}
				\begin{cases}
					-\frac{z}{2}F'(z)-F''(z)-\frac{1}{\alpha^2}\left(F^4(z)\right)''=0\\
					F(-\infty)=T_{\text{int}}^c\\
					F(\infty)=T_\infty.
				\end{cases}
			\end{equation}
		\end{enumerate}
		\begin{figure}[H]
			\begin{center}
				\includegraphics[height=5.5cm]{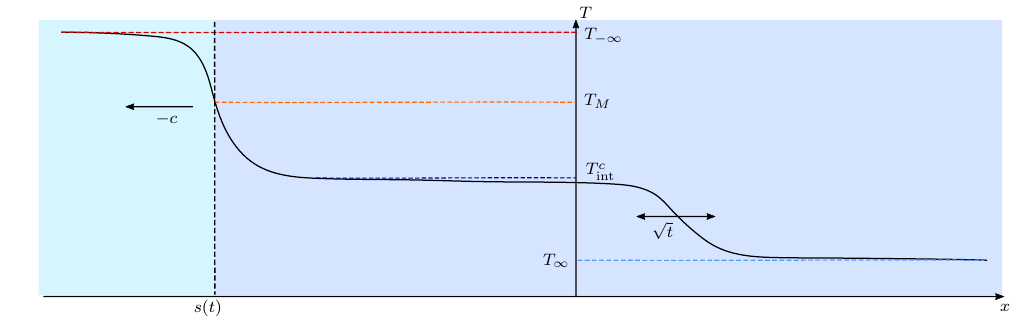}
			\end{center}\caption{Illustration of the expected profile as $ t\to\infty $.}
		\end{figure}
		\begin{remark}
			The self-similar profile $ F $ and equation \eqref{selfsimilar} can be expected due to the diffusion approximation of the radiative transfer equation. Indeed, let us define $ T_2(x,t)=F\left(\frac{x}{\sqrt{t}}\right):=F(z) $ for $ x>-ct $ as $ t\to\infty $. Then, using the H\"older regularity of $ T_2 $ we compute for the radiation term
		\[\begin{split}
			F&\left(\frac{x}{\sqrt{t}}\right)^4-\int_{-ct}^\infty \frac{\alpha E_1(\alpha(\eta-x))}{2}F\left(\frac{\eta}{\sqrt{t}}\right)^4d\eta=F(z)^4- \int_{-ct}^\infty \frac{\alpha E_1(\alpha(\eta-\sqrt{t}z))}{2}F\left(z+\frac{\eta-\sqrt{t}z}{\sqrt{t}}\right)^4d\eta\\
			=&F(z)^4-\int_{-ct-z\sqrt{t}}^\infty\frac{\alpha  E_1(\alpha\eta)}{2}\left[F^4(z)+\partial_z F^4(z)\frac{\eta}{\sqrt{t}}+\frac{\partial_z^2F^4(z)}{2}\frac{\eta^2}{t}+\mathcal{O}\left(\frac{|\eta|^{2+\delta}}{t^{1+\delta/2}}\right)\right]\\
			=&F^4(z)\int_{ct+z\sqrt{t}}^\infty\frac{\alpha  E_1(\alpha\eta)}{2}d\eta-\frac{\partial_zF^4(z)}{\sqrt{t}}\int_{ct+z\sqrt{t}}^\infty\frac{\alpha  E_1(\alpha\eta)}{2}\eta d\eta+\frac{\partial_z^2F^4(z)}{2t}\int_{-ct-z\sqrt{t}}^\infty\frac{\alpha  E_1(\alpha\eta)}{2}\eta^2 d\eta\\&+\mathcal{O}\left(\frac{\int_{-ct-z\sqrt{t}}^\infty\frac{\alpha  E_1(\alpha\eta)}{2}|\eta|^{2+\delta}d\eta}{t^{1+\delta/2}}\right).
		\end{split} \]
		\end{remark}
		Using that
		\[t\int_{\alpha(ct+z\sqrt{t})}^\infty E(\eta)d\eta\sim te^{-\alpha(ct+z\sqrt{t})}\underset{t\to\infty}{\longrightarrow} 0,\quad \sqrt{t}\int_{\alpha(ct+z\sqrt{t})}^\infty E(\eta)\eta d\eta\sim \sqrt{t}(\alpha(ct+z\sqrt{t})+1)e^{-\alpha(ct+z\sqrt{t})}\underset{t\to\infty}{\longrightarrow} 0\]\[\text{ and }\int_{-\infty}^\infty E(\eta)|\eta|^{2+\delta}d\eta<\infty\]
		we conclude multiplying by $ t $ and letting $ t\to\infty $ that
		\[t\left(F\left(\frac{x}{\sqrt{t}}\right)^4-\int_{0}^\infty E(\eta-x)F\left(\frac{\eta}{\sqrt{t}}\right)^4d\eta\right)\underset{t\to\infty}{\longrightarrow}\frac{\partial_z^2F^4(z)}{2\alpha^2}\int_{-\infty}^\infty E(\eta)\eta^2 d\eta=\frac{1}{\alpha^2}\partial_z^2F^4(z).\]
		Finally, we recover \eqref{selfsimilar} observing that $ \partial_t F\left(\frac{x}{\sqrt{t}}\right)=-\frac{z}{t}F'(z) $ and $ \partial_x^2F\left(\frac{x}{\sqrt{t}}\right)=\frac{1}{t}F''(z) $.
		
		\bibliographystyle{siam}
		\bibliography{literature_stefan}
	\end{document}